\documentclass[11pt]{article}
\usepackage{latexsym}
\usepackage{amsfonts,amssymb,amsmath}
\newtheorem{thm}{Theorem}[section]

\newtheorem{lem}[thm]{Lemma}
\newtheorem{pro}[thm]{Proposition}
\newtheorem{defn}[thm]{Definition}

\title{One-way permutations, computational asymmetry and distortion}

\author{ 
   Jean-Camille Birget
         \thanks{Supported by NSF grant CCR-0310793. \ \   
         Some of the results of this paper were presented at the AMS 
         Section Meeting, Oct.\ 21-23, 2005, Lincoln, Nebraska
         (http://www.ams.org/amsmtgs/2117\_program.html),
         and at the conference ``Various Faces of Cryptography'', 
         10 Nov.\ 2006 at City College of CUNY, New York.
                }
       }
\date{\today}
\begin{document}
\maketitle

\begin{abstract}
Computational asymmetry, i.e., the discrepancy between the complexity of
transformations and the complexity of their inverses, is at the core of 
one-way transformations.
We introduce a computational asymmetry function that measures the amount of
one-wayness of permutations. We also introduce the word-length asymmetry
function for groups, which is an algebraic analogue of computational 
asymmetry.
We relate boolean circuits to words in a Thompson monoid, over a fixed 
generating set, in such a way that circuit size is equal to word-length. 
Moreover, boolean circuits have a representation in terms of elements 
of a Thompson group, in such a way that circuit size is polynomially 
equivalent to word-length. 
We show that circuits built with gates that are not constrained
to have fixed-length inputs and outputs, are at most quadratically more 
compact than circuits built from traditional gates (with fixed-length inputs 
and outputs).
  Finally, we show that the computational asymmetry function is closely
related to certain distortion functions:
The computational asymmetry function is polynomially equivalent to the 
distortion of the path length in Schreier graphs of certain Thompson groups, 
compared to the path length in Cayley graphs of certain Thompson monoids. 
We also show that the results of Razborov and others on monotone circuit 
complexity lead to exponential lower bounds on certain distortions. 
\end{abstract}

%%%%%%%%%%%%%%%%%%%%%%%%%%%%%%%%%%%%%%%%%%%%%%%%%%%%%%%%
% Section 1
%%%%%%%%%%%%%%%%%%%%%%%%%%%%%%%%%%%%%%%%%%%%%%%%%%%%%%%%

\section{Introduction}

The existence of one-way functions, i.e., functions that are ``easy to
evaluate'' but ``hard to invert'', is a major open problem. Much of 
cryptography depends on one-way functions; moreover, indirectly, their 
existence is connected to the question whether {\sf P} is different from 
{\sf NP}. In this paper we give some connections between these questions and 
some group-theoretic concepts: \\   
(1) \ We continue the work of \cite{BiCoNP}, \cite{BiFact}, and 
\cite{BiThomMon}, on the relation between combinational circuits, on the one 
hand, and Thompson groups and monoids on the other hand. 
We give a representation of any circuit by a word over the Thompson group,
such that circuit size is polynomially equivalent to word-length. \\    
(2) \ We establish connections between the existence of one-way permutations
and the distortion function in a certain Thompson group. Distortion is an
important concept in metric spaces (e.g., Bourgain \cite{Bourgain}) and in 
combinatorial group theory (e.g., Gromov \cite{Gromov}, Farb \cite{Farb}).

\bigskip

\noindent {\bf Overview:}  \\     
Subsections {\bf 1.1 - 1.6} of the present Section define and motivate the 
concepts used: One-way functions and one-way permutations; computational 
asymmetry; word-length asymmetry; reversible computing; distortion; Thompson 
groups and monoids.
In Section {\bf 2} we show that circuits can be represented by elements of 
Thompson monoids: A boolean circuit is equivalent to a word over a fixed 
generating set of a Thompson monoid, with circuit size being equal (or 
linearly equivalent) to word-length over the generating set. The Thompson
monoids that appear here are monoid generalizations of the Thompson group 
$G_{2,1}$, obtained when bijections are generalized to partial functions 
\cite{BiThomMon}.
Section {\bf 3} shows that computational asymmetry and word-length 
asymmetry (for the Thompson groups and monoids) are linearly related.
In Section {\bf 4} we give a representation of arbitrary (not necessarily 
bijective) circuits by elements of the Thompson group $G_{2,1}$; circuit size 
is polynomially equivalent to word-length over a certain generating set in 
the Thompson group. In Section {\bf 5} we show that the 
computational asymmetry function of permutations is polynomially related to
a certain distortion in a Thompson group. Section {\bf 6} contains 
miscellaneous results, in particular that the work of Razborov and others 
on monotone circuit complexity leads to exponential lower bounds on certain 
distortion functions.

%%%%%%%%%%%%%%%%%%%%%%%%%%%%%%%%%%%%%%%%%%%%%%%%%%%%%%%%%%%%%%%
\subsection{ One-way functions and one-way permutations}

Intuitively, a one-way function is a function $f$ (mapping words to words, 
over a finite alphabet), such that $f$ is ``easy to evaluate'' (i.e., 
given $x_0$ in the domain, it is ``easy'' to compute $f(x_0)$), but ``hard 
to invert'' (i.e., given $y_0$ in the range, it is ``hard'' to find any 
$x_0$ such that $f(x_0) = y_0$). The concept was introduced by Diffie and 
Hellman \cite{DiffieHellman}. 

There are many ways of defining the words ``easy'' and ``hard'', and 
accordingly there exist many different rigorous notions of a one-way 
function, all corresponding to a similar intuition.
It remains an open problem whether one-way functions exist, for any 
``reasonable'' definition. Moreover, for certain definitional choices, this 
problem is a generalization of the famous question whether {\sf P} $\neq$
{\sf NP} \ \cite{GrollSel,Selman,BoppLag}.
 
\medskip

We will base our one-way functions on combinational circuits and their size. 
The size of a circuit will also be called its complexity.
Below, $\{0,1\}^n$ (for any integer $n \geq 0$) 
denotes the set of all bitstrings of length $n$.
A combinational circuit with input-output function
$f: \{0,1\}^m \to \{0,1\}^n$ is an acyclic boolean circuit with $m$ input
wires (or ``input ports'') and $n$ output wires (or ``output ports''). The 
circuit is made from gates of type {\sf not, and, or, fork}, as well as 
wire-crossings or wire-swappings. These gates are very traditional and are 
defined as follows. 

\smallskip

{\sf and}$: (x_1,x_2) \in \{0,1\}^2 \longmapsto y \in \{0,1\}$, \ where
$y = 1$ if $x_1 = x_2 = 1$, and $y = 0$ otherwise. 

\smallskip

{\sf or}$: (x_1,x_2) \in \{0,1\}^2 \longmapsto y \in \{0,1\}$, \ where
$y = 0$ if $x_1 = x_2 = 0$, and $y = 1$ otherwise. 

\smallskip

{\sf not}$: x \in \{0,1\} \longmapsto y \in \{0,1\}$, \ where
$y = 0$ if $x = 1$, $y = 1$ otherwise. 

\smallskip

{\sf fork}$: x \in \{0,1\} \longmapsto (x,x) \in \{0,1\}^2$.

\smallskip

\noindent Another gate that is often used is the exclusive-or gate, 

\smallskip

{\sf xor}$: (x_1,x_2) \in \{0,1\}^2 \longmapsto y \in \{0,1\}$, \ where
$y = 1$ if $x_1 \neq x_2$, and $y = 0$ otherwise. 
 
\smallskip

\noindent The wire-swapping of the $i$th and $j$th wire ($i<j$) is described 
by the {\it bit transposition} (or bit position transposition) \

\smallskip

$\tau_{i,j} : ux_ivx_jw \in \{0,1\}^{\ell} \longmapsto ux_jvx_iw \in
\{0,1\}^{\ell}$, \ 
where \ $|u| = i-1$, \ $|v| = j-i-1$, \ $|w| = \ell -j-1$. 

\smallskip

\noindent The {\sf fork} and wire-swapping operations, although heavily used,
are usually not explicitly called ``gates''; but because of their important 
role we will need to consider them explicitly.
Other notations for the gates: \ {\sf and}$(x_1,x_2) = x_1 \wedge x_2$, 
 \ {\sf or}$(x_1,x_2) = x_1 \vee x_2$, \ {\sf not}$(x) = \overline{x}$, 
 \ {\sf xor}$(x_1,x_2) = x_1 \oplus x_2$.

A combinational circuit for a function $f: \{0,1\}^m \to \{0,1\}^n$ is 
defined by an acyclic directed graph drawn in the plane (with crossing of 
edges allowed). In the circuit drawing, 
the $m$ input ports are vertices lined up in a vertical column on the left 
end of the circuit, and the $n$ output ports are vertices lined up in a 
vertical column on the right end of the circuit. The input and output ports 
and the gates of the circuit (including the {\sf fork} gates, but not the wire 
transpositions) form the vertices of the circuit graph. We often view the 
circuit as cut into vertical {\it slices}. 
A slice can be any collection of gates and wires in the circuit such that 
no gate in a slice is an ancestor of another gate in the same slice, and no 
wire in a slice is an ancestor of another wire in the same slice 
(unless these two wires are an input wire and an output wire of a same gate). 
Two slices do not overlap, and every wire and every gate belongs to some
slice. For more details on combinational circuits, see 
\cite{JESavage, Wegener, BoppLag}. 

The {\it size} of a combinational circuit is defined to be the number of gates
in the circuit, including forks and wire-swappings, as well as the input ports
and the output ports. For a function $f: \{0,1\}^m \to \{0,1\}^n$, the 
{\it circuit complexity} (denoted $C(f)$) is the smallest size of any 
combinational circuit with input-output function $f$. 

A cause of confusion about gates in a circuit is that gates of a certain type
(e.g., {\sf and}) are traditionally considered the same, no matter where they 
occur in the circuit. However, gates applied to different wires in a circuit 
are different functions; e.g., for the {\sf and} gate, 
$(x_1,x_2,x_3) \mapsto (x_1 \wedge x_2, x_3)$ \ is a different 
function than \, $(x_1,x_2,x_3) \mapsto (x_1, x_2 \wedge x_3)$.

%%%%%%%%%%%%%%%%%%%%%%%%%%%%%%%%%%%%%%%%%%%%%%%%%%%%%%%%%%%%%%%
\subsection{Computational Asymmetry}

Computational asymmetry is the core property of one-way functions.
Below we will define computational asymmetry in a quantitative way, and in 
a later Section we will relate it to the group-theoretic notion of distortion.

For the existence of one-way functions, it is mainly the relation 
between the circuit complexity $C(f)$ of $f$ and the circuit complexity 
$C(f^{-1})$ of $f^{-1}$ that matters, not the complexities of $f$ and 
of $f^{-1}$ themselves.
Indeed, a classical {\em padding argument} can be used: If we add $C(f)$ 
``identity wires'' to a circuit for $f$, then the 
resulting circuit has linear size as a function of its number of input
wires; see Proposition \ref{compAsymm_1w} below. (An identity wire is a
wire that goes directly from an input port to an output port, without being
connected to any gate.) 

In \cite{BoppLag} (page 230) Boppana and Lagarias considered
${\rm log} \, C(f') / {\rm log} \, C(f)$ as a measure of one-wayness; here,
$f'$ denotes an inverse of $f$, i.e., any function such that 
$f \circ f' \circ f =f$.
Massey and Hiltgen \cite{Massey,HiltgenDiss} introduced the phrases 
{\em complexity asymmetry} and {\em computational asymmetry} for injective 
functions, in reference to the situation where the circuit complexities 
$C(f)$ and $C(f^{-1})$ are very different. 
The concept of computational asymmetry can be generalized to arbitrary 
(non-injective) functions, with the meaning that for every inverse $f'$ of 
$f$, $C(f)$ and $C(f')$ are very different. 

\smallskip

In \cite{Massey} Massey made the following observation.
For any large-enough fixed $m$ and for almost all permutations $f$ of 
$\{0,1\}^m$, the circuit complexities $C(f)$ and $C(f^{-1})$ are very 
similar: 

\medskip

 \ \ \   \ \ \ \ \   
$\frac{1}{10} \ C(f) \ \leq \ C(f^{-1}) \ \leq \ 10 \ C(f)$ 

\medskip

\noindent 
Massey's proof is adapted from the Shannon lower bound \cite{Shannon} 
and the Lupanov upper bound \cite{Lupanov} (see also \cite{HiltgenDiss}, 
\cite{JESavage}), from which it follows that almost all functions and almost 
all permutations (and their inverses) have circuit complexity close to 
the Shannon bounds.
Massey's observation can be extended to the set of all functions 
$f: \{0,1\}^m \to \{0,1\}^n$, i.e., for almost all $f$ and for every inverse 
$f'$ of $f$, the complexities $C(f)$ and $C(f')$ are within constant factors 
of each other. 

Hence, computationally asymmetric permutations are rare among the boolean 
permutations overall (and similarly for functions). This is an interesting
fact about computational asymmetry, but by itself it does not imply anything 
about the existence or non-existence of one-way functions, not even 
heuristically. Indeed, Massey proved his linear relation 
$C(f) = \Theta(C(f'))$ in the situation where $C(f) = \Theta(2^m)$, and 
then uses the fact that the condition $C(f) = \Theta(2^m)$ holds for almost 
all boolean permutations and for almost all boolean functions.  
But there also exist functions with $C(f) = O(m^k)$, with $k$ a small 
constant. In particular, one-way functions (if they exist) have small 
circuits; by definition, one-way functions violate the condition 
$C(f) = \Theta(2^m)$.  

A well-known candidate for a one-way permutation is the following. For a
large prime number $p$ and a primitive root $r$ modulo $p$, consider
the map \  $x \in \{0, 1, \ldots, p-2 \} \ \longmapsto \ $
$r^x - 1 \ \in \{0, 1, \ldots, p-2 \}$. 
This is a permutation whose inverse, known as the {\it discrete logarithm}, 
is believed to be difficult to compute.

\bigskip

%%%%%%%%%%%%%%%%%%%%%

\noindent {\bf Measuring computational asymmetry:} 

\smallskip

\noindent Let ${\mathfrak S}_{\{0,1\}^m}$ denote the set of all permutations
of $\{0,1\}^m$, i.e., ${\mathfrak S}_{\{0,1\}^m}$ is the symmetric 
group. We will measure the computational asymmetry of all permutations
of $\{0,1\}^m$ (for all $m > 0$) by defining a {\em computational
asymmetry function}, as follows. 
A function $a: {\mathbb N} \to {\mathbb N}$ is an upper bound on the
computational asymmetry function iff for all all $m > 0$ and all permutations
$f$ of $\{0,1\}^m$ we have: \ \ \ $C(f^{-1}) \leq a\big(C(f)\big)$.
The computational asymmetry function $\alpha$ of the boolean permutations
is the least such function $a(.)$. Hence: 

\begin{defn} \label{def_compAsymm} \
The computational asymmetry function $\alpha$ of the boolean permutations 
is defined as follows for all $s \in {\mathbb N}:$
 \ \ \ $\alpha(s) \ = \ {\rm max}\big\{C(f^{-1}) : $
$ \ C(f) \leq s, \ f \in {\mathfrak S}_{\{0,1\}^m}, \ m >0 \big\}.$   
\end{defn}

\noindent
Note that in this definition we look at all combinational circuits,
for all permutations in $\bigcup_{m>0} {\mathfrak S}_{\{0,1\}^m}$;
we don't need to work with non-uniform or uniform families of circuits.

\medskip

\noindent Computational asymmetry is closely related to one-wayness, as 
the next proposition shows.

\begin{pro} \label{compAsymm_1w} \  

\noindent {\bf (1)} \  
For infinitely many $n$ we have: There exists a permutation $f_n$ of 
$\{0,1\}^n$ such that $f_n$ is computed by a circuit of size $\leq 3 \, n$, 
but $f_n^{-1}$ has no circuit of size $< \alpha(n)$.

\noindent {\bf (2)} \  Suppose that $\alpha$ is exponential, i.e., there is 
$k > 1$ such that for all $n$, \ $\alpha(n) \geq k^n$.
Then $k \leq 2$, and there is a constant $c > 1$ such that we have:
For every integer $n \geq 1$ there exists a permutation $F_n$ of $\{0,1\}^n$ 
which is computed by a circuit of size $\leq c \, n$, but $F_n^{-1}$ has no 
circuit of size $< k^n$.
\end{pro}
{\bf Proof. (1)} By the definition of $\alpha$, for every $m>0$ there exists 
a permutation $F$ of $\{0,1\}^m$ such that $F$ is computed by a circuit 
of some size $C_F$, but $F^{-1}$ has no circuit of size $< \alpha(C_F)$.
Let $n = C_F$, and let us consider the function 
$f_n: \{0,1\}^{C_F} \to \{0,1\}^{C_F}$ defined by 
$f_n: (x, w) \longmapsto (F(x), w)$, for all $x \in \{0,1\}^m$ and
$w \in \{0,1\}^{C_F-m}$. 

Then $f_n(x,w)$ is computed by a circuit of size $C_F + 2 \, (C_F - m)$; the 
term ``$2 \, (C_F-m)$'' comes from counting the input-output wires of $w$. 
Hence $f_n$ has a circuit of size $\leq 3n$. On the other hand, 
$(y,w) \longmapsto f_n^{-1}(y,w) = (F^{-1}(y),w)$
is not computed by any circuit of size $< \alpha(C_F)$, so $f_n^{-1}$ has no
circuit of size $< \alpha(n)$.

\smallskip

\noindent {\bf (2)} For every $n \geq 1$ there exists a permutation $F$ of 
$\{0,1\}^n$ such that $F$ is computed by a circuit of some size $C_F$, and 
$F^{-1}$ has a circuit of size $C_{F^{-1}} = \alpha(C_F) \geq k^{C_F}$; 
moreover, $F^{-1}$ has no circuit of size $< \alpha(C_F)$.   Thus, 
$k^{C_F} \leq C_{F^{-1}} \leq 2^n \, (1 + c_o \, \frac{\log n}{n})$, for 
some constant $c_o > 1$; the latter inequality comes from the Lupanov upper 
bound \cite{Lupanov} (or see Theorem 2.13.2  in \cite{JESavage}).   
Hence, $k \leq 2$ and \ $n \leq C_F \leq \frac{1}{\log_2 k} \ n \ + $
$c_1 \, \frac{\log n}{n}$, for some constant $c_1 > 0$. Hence, for all 
$n \geq 1$ there exists a permutation $F$ of $\{0,1\}^n$ with circuit size 
$C_F \in \ $ 
$[n, \ \frac{1}{\log_2 k} \cdot n + c_1 \cdot \frac{\log n}{n}]$, such 
that \ $C_{F^{-1}} = \alpha(C_F) \geq k^{C_F} \geq k^n$.
 \ \ \ $\Box$

\medskip

We will show later that the computational asymmetry function is closely 
related to the distortion of certain groups within certain monoids.
 
\bigskip

\noindent {\bf Remarks:} \ 

\noindent Although in this paper we only use the computational asymmetry
function of the boolean permutations, the concept can be generalized. 
Let \ ${\rm Inj}(\{0,1\}^m, \{0,1\}^n)$ \ denote the set of all 
{\em injective functions} $\{0,1\}^m \to \{0,1\}^n$.
The computational asymmetry function $\alpha_{\rm inj}$ of the injective 
boolean functions is defined by  

\smallskip

$\alpha_{\rm inj}(s) \ = \ {\rm max}\big\{C(f^{-1}) : \ $
$C(f) \leq s, \ f \in {\rm Inj}(\{0,1\}^m, \{0,1\}^n), \ m>0, n >0 \big\}$

\smallskip

\noindent More generally, let \ $(\{0,1\}^n)^{\{0,1\}^m}$ \ denote  
the set of all {\em functions} $\{0,1\}^m \to \{0,1\}^n$. 
The computational asymmetry of all finite boolean functions is defined by

\smallskip

$\alpha_{\rm func}(s) \ = \ {\rm max}\big\{C(f') : \ C(f) \leq s, \ \ $
$ff'f =f, \ \ f, f' \in (\{0,1\}^n)^{\{0,1\}^m}, \ n >0, m >0 \big\}$.

\bigskip

When we compare functions we will be mostly interested in their asymptotic
growth pattern. Hence we will often use the big-O notation, and the 
following definitions.

By definition, two functions $f_1: {\mathbb N} \to {\mathbb N}$
and $f_2: {\mathbb N} \to {\mathbb N}$ are {\em linearly equivalent} iff
there are constants $c_0, c_1, c_2 > 0$ such that for all $n \geq c_0:$ \
$f_1(n) \leq c_1 \, f_2(c_1n)$ \ and \ $f_2(n) \leq c_2\, f_1(c_2n)$.
Notation: \ $f_1 \ \simeq_{{\rm lin}} \ f_2$.

Two functions $f_1$ and $f_2$ (from ${\mathbb N}$ to ${\mathbb N}$) are
called {\em polynomially equivalent} iff there are constants
$c_0, c_1, c_2, d, e > 0$ such that for all $n \geq c_0:$ \
$f_1(n) \leq c_1 \, f_2(c_1n^d)^d$ \ and \ $f_2(n) \leq c_2\, f_1(c_2n^e)^e$.
Notation: \ $f_1 \ \simeq_{{\rm poly}} \ f_2$.

%%%%%%%%

\subsection{Wordlength asymmetry}

We introduce an algebraic notion that looks very similar to computational
asymmetry:

\begin{defn} \label{w_length_asym} \
Let $G$ be a group, let $M$ be a monoid with generating set $\Gamma$ (finite
or infinite), and suppose \ $G \subseteq M$. The {\em word-length 
asymmetry function} of $G$ within $M$ (over $\Gamma$) is 

\smallskip

 \ \ \ \ \  $\lambda(n) \ = \ {\rm max}\{ \, |g^{-1}|_{\Gamma} : \ $
      $|g|_{\Gamma} \leq n, \ g \in G \}$.
\end{defn}
The word-length asymmetry function $\lambda$ depends on $G$, $M$, $\Gamma$,
and the embedding of $G$ in $M$. 

Consider the right Cayley graph of the monoid $M$ with generating
set $\Gamma$; its vertex set is $M$ and the edges have the form
 \ $x \stackrel{\gamma}{\longrightarrow} \gamma x$ \ (for $x \in M$,
$\gamma \in \Gamma$). For $x, y \in M$, the {\em directed distance} $d(x,y)$ 
in the Cayley graph is the shortest length over all paths from $x$ to $y$ in 
the Cayley graph; if no path from $x$ to $y$ exists, the directed distance 
is infinite. By ``path'' we always mean directed path.

\begin{lem} \label{Cayley_g_to_1} \    
Under the above conditions on $G$, $M$, $\Gamma$, we have for every 
$g \in G:$ \ \ $d({\bf 1}, g^{-1}) = d(g, {\bf 1})$ \ and \   
$d({\bf 1}, g) = d(g^{-1}, {\bf 1})$. 
\end{lem}
{\bf Proof.} Let $\eta:  \Gamma^* \to M$ be the map that evaluates generator
sequences in $M$.
If $v \in \Gamma^*$ is the label of a shortest path from {\bf 1} to
$g^{-1}$ in the Cayley graph then $g \cdot \eta(v) = {\bf 1}$ in $M$, hence 
$\eta(v) = g^{-1}$. Therefore, the path starting at $g$ and labeled by $v$
ends at {\bf 1}; hence $d(g, {\bf 1}) \leq |v| = d({\bf 1}, g^{-1})$. 
In a similar way one proves that $d({\bf 1}, g^{-1}) \leq d(g, {\bf 1})$. 
The equality $d({\bf 1}, g) = d(g^{-1}, {\bf 1})$ is also 
proved in a similar way.  \ \ \ $\Box$

\medskip

Since $|g|_{\Gamma}$ is the distance $d({\bf 1},g)$ in the  graph of $M$,
and since $|g^{-1}|_{\Gamma} = d({\bf 1},g^{-1}) = d(g, {\bf 1})$,
the word-length asymmetry also measures the asymmetry of the directed 
distance, to or from the identity element {\bf 1} in the Cayley graph of $M$, 
restricted to vertices in the subgroup $G$.  

For distances to or from the identity element of $M$ it does not matter 
whether we consider the left Caley graph or the right Caley graph.

%%%%%%%%%%%%%%%

\subsection{Computational asymmetry and reversible computing}

Reversible computing deals with the following questions:
If a function $f$ is injective (or bijective) and computable, can $f$ be 
computed in such a way that each elementary computation step is injective 
(respectively bijective)?
And if such injective (or bijective) computations are possible, what is
their complexity, compared to the usual (non-injective) complexity? 

One of the main results is the following (Bennett's theorem 
\cite{Bennett73, Bennett89}, and earlier work of Lecerf \cite{Lec}): 
Let $f$ be an injective function, and assume $f$ and $f^{-1}$ are 
computable by deterministic Turing machines with time complexity $T_f(.)$, 
respectively $T_{f^{-1}}(.)$. Then $f$ (and also $f^{-1}$) is computable by
a {\it reversible Turing machine} (in which every transition is 
deterministic and injective) with time complexity $O(T_f + T_{f^{-1}})$. 
Note that only injectiveness (not bijectiveness) is used here. 

Bennett's theorem has the following important consequence, which relates 
reversible computing to one-way functions: {\em Injective one-way functions 
exist iff there exist injective functions that have efficient traditional 
algorithms but that do not have efficient reversible algorithms.}  

\bigskip

\noindent {\bf Toffoli representation}

\smallskip

Remarkably, it is possible to ``simulate'' any function 
$f: \{0,1\}^m \to \{0,1\}^n$ (injective or not, one-way or not) by a 
bijective circuit; a circuit is called bijective iff the circuit is
made from bijective gates. Here, bijective circuits will be built from 
the wire swapping operations and the following bijective gates: 
{\sf not} (negation), {\sf c-not} (the Controlled Not, also called 
``Feynman gate'') defined by \     
$(x_1,x_2) \in \{0,1\}^2 \longmapsto (x_1, \, x_1 \oplus x_2) \in \{0,1\}^2$, 
and {\sf cc-not} (the Doubly Controlled Not, also called ``Toffoli gate'')
defined by \ $(x_1,x_2,x_3) \in \{0,1\}^3 \longmapsto$
$(x_1, \, x_2, \, (x_1 \wedge x_2) \oplus x_3) \in \{0,1\}^3$.

\begin{thm} {\rm (Toffoli \cite{Toff80Memo}).} \label{ToffRepr} \    
For every boolean function $f: \{0,1\}^m \to \{0,1\}^n$ there exists a 
bijective boolean circuit $\beta_f$ (over the bijective gates 
{\sf not}, {\sf c-not}, {\sf cc-not}, and wire transpositions), 
with input-output function \ $\beta_f:$
$x \, 0^n \ \in \{0,1\}^{m+n} \ \longmapsto \ f(x) \ x \ \in \{0,1\}^{n+m}$. 
\end{thm} 

\noindent
In other words, $f(x)$ consists of the projection onto the first $n$ bits of
$\beta_f(x \, 0^n)$; equivalently, 
 \ $f(.) = {\sf proj}_n \circ \beta_f \circ {\sf concat}_{0^n}(.)$, where
${\sf proj}_n$ projects a string of length $n+m$ to the first $n$ bits,
and ${\sf concat}_{0^n}$ concatenates $0^n$ to the right of a string.
See Theorems 4.1, 5.3 and 5.4 of \cite{Toff80Memo}, and see Fig.\ 1 below. 

\vspace{-1in}

\unitlength=0.90mm \special{em:linewidth 2pt}
\linethickness{0.4pt}

\begin{picture}(150,80)

\put(40,0){\framebox(40,40)}
\put(60,20){\makebox(0,0)[cc]{$\beta_f$}}

%bottom left 0^n
\put(28,10){\vector(1,0){12}}
\put(24,10){\makebox(0,0)[cc]{$0^n$}}
\put(34,7){\makebox(0,0)[cc]{$_n$}}
\put(34,10){\makebox(0,0)[cc]{{\tiny / }}}

%top left x
\put(12,30){\vector(1,0){28}}
\put( 9,30){\makebox(0,0)[cc]{$x$}}
\put(30,27){\makebox(0,0)[cc]{$_m$}}
\put(30,30){\makebox(0,0)[cc]{{\tiny / }}}

%bottom right x
\put(80,10){\vector(1,0){15}}
\put(97,10){\makebox(0,0)[cc]{$ \ x $}}
\put(86,7){\makebox(0,0)[cc]{$_m$}}
\put(86,10){\makebox(0,0)[cc]{{\tiny / }}}

%top right f(x)
\put(80,30){\vector(1,0){28}}
\put(114,30){\makebox(0,0)[cc]{$f(x)$}}
\put(91,27){\makebox(0,0)[cc]{$_n$}}
\put(91,30){\makebox(0,0)[cc]{{\tiny / }}}

\end{picture}

\bigskip

\bigskip

{\sf Fig.\ 1: \ Toffoli representation of the function $f$. }

\vspace{.3in}
 
The Toffoli representation contains two non-bijective actions:
The projection at the output, and the forced setting of the value of some 
of the input wires. 

\smallskip

%There is a special Toffoli representation of boolean permutations, by 
%bijective circuits that have ``storage wires'', but no forced input settings
%and no projections at the output: \ 
%For every boolean permutation $g: \{0,1\}^m \to \{0,1\}^m$ there
%exists a bijective boolean circuit $\pi_g$ (over the bijective gates
%seen above), with input-output function $\{0,1\}^{4m} \to \{0,1\}^{4m}$, 
%such that for all $x \in \{0,1\}^m$, and for all $z \in \{0,1\}^{3m}$,
%$\pi_g(x z) = g(x) \, z$. (See Theorems 5.3 and 5.4 in \cite{Toff80Memo}.)
%The input substring $z$ re-appears unchanged in the output, on the same wires, 
%but the circuit uses these wires between input and output. 
%
%\smallskip

Toffoli's proofs and constructions are based on truth tables, and he does not
prove anything about the circuit size of $\beta_f$ (counting the bijective 
gates),
compared to the circuit size of $f$. The following gives a polynomial bound
on the size of the bijective circuit, at the expense of a large number of 
input- and output-wires.

\begin{thm} {\rm (E.~Fredkin, T.~Toffoli \cite{FredToff}).} 
 \label{FredToffRepr} \
For every boolean function $f: \{0,1\}^m \to \{0,1\}^n$ with circuit size
$C(f)$ there exists a bijective boolean circuit $B_f$ (over a bounded 
collection of bijective gates, e.g., {\sf not}, {\sf c-not}, {\sf cc-not}, 
and wire transpositions), with input-output function 

\smallskip

$B_f: \ x \ 0^{n+C(f)} \ \in \ \{0,1\}^{m+n+C(f)} \ \ \longmapsto \ \ $
$f(x) \ z(x) \ \in \ \{0,1\}^{m+n+C(f)}$ \ \ \ 

\smallskip

\noindent for some $z(x) \in \{0,1\}^{m+C(f)}$. 

If $g: \{0,1\}^m \to \{0,1\}^m$ is a permutation then there exists a bijective
boolean circuit $U_g$ (over bijective gates), with input-output function

\smallskip

$U_g: \ x \ 1^m \ 0^{m+C} \ \in \ \{0,1\}^{3m+C} \ \ \longmapsto \ \ $
$g(x) \ {\overline {g(x)}} \ x \ 0^C \ \in \ \{0,1\}^{3m+C}$    

\smallskip

\noindent where \ $C = {\rm max}\{C(g), \, C(g^{-1})\}$, and ${\overline
{g(x)}}$ is the bitwise complement of $g(x)$. 
\end{thm} 

Later we will introduce another reversible representation of boolean 
functions by bijective gates; we will need only one 0-wire, but the gates 
will be taken from the Thompson group $G_{2,1}$, i.e., we will also use 
non-length-preserving transformations of bitstrings (Theorems 
\ref{repres_M_by_G} and \ref{repres_G_circuits_by_G} below).

%%%%%%%%%%%%%%%%%%%%%%%%%%%%%%%%%%%%%%%%%%%%%%%%%%%

\subsection{Distortion}   % distorsion, Verzerrung 

We will prove later (Theorem \ref{polyn_rel_alpha_delta}) that computational 
asymmetry has a lot to do with distortion, a concept introduced into group 
theory by Gromov \cite{Gromov} and Farb \cite{Farb}. Distortion is already 
known to have connections with isoperimetric functions (see \cite{Olsh}, 
\cite{OlSap}, \cite{MargMeakSunik}).
A somewhat different problem about distortion (for finite metric spaces) 
was tackled by Bourgain \cite{Bourgain}.

We will use a slightly more general notion of distortion, based on (possibly
directed) countably infinite rooted graphs, and their (directed) path metric. 

A weighted directed graph is a structure $(V, E, \omega)$ where $V$ is a set
(called the vertex set), $E \subseteq V \times V$ (called the edge set), and
$\omega: E \longmapsto {\mathbb R}_{> 0}$ is a function (called the weight 
function); note
that every edge has a strictly positive weight. It is sometimes convenient 
to define $\omega(u,v) = \infty$ when $(u,v) \in V \times V - E$. 
A path in $(V, E)$ is a sequence of edges $(u_i,v_i)$ ($1 \leq i \leq n$) 
such that $u_{i+1} = v_i$ for all $i<n$, and such that all elements in 
$\{u_i : 1 \leq i \leq n\} \cup \{v_n\}$ are distinct;
$u_1$ is called the start vertex of this path, and $v_n$ is called the end
vertex of this path; the sum of weights \ $\sum_{i=1}^n \omega(u_i,v_i)$ \ 
over the edges in the path is called the {\em length} of the path. 
Here we do not consider any paths with infinitely many edges; but we allow 
$V$ and $E$ to be countably infinite.
A vertex $w_2$ is said to be reachable from a vertex $w_1$ in $(V, E)$ iff
there exists a path with start vertex $w_1$ and end vertex $w_2$.  
If $w_2$ is reachable from $w_1$ then the minimum length over all paths from
$w_1$ to $w_2$ is called the {\em directed distance} from $w_1$ to $w_2$,
denoted $d(w_1,w_2)$; since we only consider finite paths here, this minimum 
exists. If $w_2$ is not reachable from $w_1$ then we define
$d(w_1,w_2)$ to be $\infty$. Clearly we have $w_1 = w_2$ iff $d(w_1,w_2) =0$,
and for all $u,v,w \in V$, $d(u,w) \leq d(u,v) + d(v,w)$. In a directed
graph, the function $d(.,.)$ need not be symmetric. The function 
$d: V \times V \to {\mathbb R}_{\geq 0} \cup \{ \infty\}$ is called the 
directed path metric of $(V, E, \omega)$.
A rooted directed weighted graph is a structure $(V, E, \omega, r)$ where
$(V, E, \omega)$ is a directed weighted graph, $r \in V$, and all vertices
in $V$ are reachable from $r$.

A set $M$ with a function 
$d: M \times M \to {\mathbb R}_{\geq 0} \cup \{ \infty\}$, satisfying 
the two axioms $w_1 = w_2$ iff $d(w_1,w_2) =0$, and $d(u,w) \leq d(u,v) +
d(v,w)$, will be called {\em directed metric space} (a.k.a.\ quasi-metric
space).

Any subset $G$ embedded in a directed metric space $M$ becomes a directed 
metric space by using the directed distance of $M$.  
We call this the {\it directed distance on $G$ inherited from} $M$.

If $G \subseteq V$ for a rooted directed weighted graph $(V, E, \omega, r)$,
we consider the function 
 \ $\ell: g \in G \longmapsto d(r,g) \in {\mathbb R}_{\geq 0}$,
which we call the {\em directed length function} on $G$ inherited from 
$(V, E, \omega, r)$. (The value $\infty$ will not appear here since all of
$G$ is reachable from $r$.)   

\smallskip
 
We now define distortion in a very general way. Intuitively, distortion in
a set is a quantitative comparison between two (directed) length functions 
that are defined on the same set. 

\begin{defn} \label{generalDistortion} \ 
Let $G$ be a set, and let $\ell_1$ and $\ell_2$ be two functions 
$G \to {\mathbb R}_{\geq 0}$.
The {\bf distortion} of $\ell_1$ with respect to $\ell_2$ is the function
$\delta_{\ell_1,\ell_2} :  {\mathbb R}_{\geq 0} \to {\mathbb R}_{\geq 0}$ 
defined by

\medskip

\hspace{1in}
$\delta_{\ell_1,\ell_2}(n) \ = \ {\rm max}\{ \ell_1(g) : g \in G, $
                                 $ \ \ell_2(g) \leq n \}$. 
\end{defn} 
We will also use the notation $\delta[\ell_1,\ell_2](.)$ for
$\delta_{\ell_1,\ell_2}(.)$.
When we consider a distortion $\delta_{\ell_1,\ell_2}(.)$ we often assume 
that $\ell_2 \leq \ell_1$ or $\ell_2 \leq O(\ell_1)$; this insures that the 
distortion is at least linear, i.e., \,  
$\delta_{\ell_1,\ell_2}(n) \geq c \, n$,
for some constant $c>0$. We will only deal with functions obtained from the
lengths of finite paths in countable directed graphs, so in that case the
functions $\ell_i$ are discrete, and the distortion function exists. 
The next Lemma generalizes the distortion result of Prop.\ 4.2 of \cite{Farb}.  

\begin{lem} \label{distor_composition} \
Let $G$ be a set and consider three functions $\ell_3, \ell_2, \ell_1:$
$G \to {\mathbb R}_{\geq 0}$ such that 
$\ell_1(.) \geq \ell_2(.) \geq \ell_3(.)$. 
Then the corresponding distortions satisfy: \ \ \  
$\delta_{\ell_1, \ell_3}(.) \ \leq \ \delta_{\ell_1, \ell_2} \, \circ \, $
        $\delta_{\ell_2, \ell_3}(.)$.
\end{lem}
{\bf Proof.} The inequalities $\ell_1(.) \geq \ell_2(.) \geq \ell_3(.)$
guarantee that the three distortions $\delta_{\ell_1, \ell_3}$, \ 
$\delta_{\ell_1, \ell_2}$, and $\delta_{\ell_2, \ell_3}$ are at least as 
large as the identity map.  By definition, 

\smallskip

$\delta_{\ell_1,\ell_2} \big( \delta_{\ell_2,\ell_3}(n) \big) \ = \ $
${\rm max}\{\ell_1(x): x \in G, \ $
$ \ell_2(x) \leq \delta_{\ell_2,\ell_3}(n)\}$

\smallskip

$ = \ {\rm max}\big\{\ell_1(x): x \in G, \ \ \ell_2(x) \leq $
      ${\rm max}\{\ell_2(z): z \in G, \ \ell_3(z) \leq n \} \big\}$

\smallskip

$= \ {\rm max}\big\{\ell_1(x): x \in G, \ $
$(\exists z \in G)\big(\ell_2(x) \leq \ell_2(z)$ 
$ \ {\rm and} \ $
$\ell_3(z) \leq n \big) \big\}$

\smallskip

$\geq \ {\rm max} \{\ell_1(x): x \in G, \ \ell_3(x) \leq n \} \ = \ $
$\delta_{\ell_1, \ell_3}(n)$.

\smallskip

\noindent The last inequality follows from the fact that if
$\ell_3(x) \leq n$ then for some $z$ (e.g., for $z = x$): \,
$\ell_2(x) \leq \ell_2(z)$ and $\ell_3(z) \leq n$.  \ \ \ $\Box$

\bigskip

\bigskip

\noindent {\large \bf Examples of distortion:}

\medskip

\noindent Distortion and asymmetry are unifying concepts that apply to 
many fields. 

\medskip

\noindent
{\bf 1. Gromov distortion:} \ Let $G$ be a subgroup of a group $H$, with 
generating sets $\Gamma_G$, respectively $\Gamma_H$, such that 
$\Gamma_G \subseteq \Gamma_H$, and such that $\Gamma_G = \Gamma_G^{-1}$ and
$\Gamma_H = \Gamma_H^{-1}$.  This determines a Cayley graph for $G$ and a
Cayley graph for $H$. Now we have two distance functions on $G$, one 
obtained from the Cayley graph of $G$ itself (based on $\Gamma_G$), and the 
other inherited from the embedding of $G$ in $H$. See \cite{Gromov},
\cite{Bourgain}, and \cite{Farb}.   

\smallskip

The Gromov distortion function is a natural measure of the difficulty
of the {\em generalized word problem}.
A very important case is when both $\Gamma_G$ and $\Gamma_H$ are finite.
Here are some results for that case:

Theorem of Ol$'$shanskii and Sapir \cite{OlSap} (making precise and proving 
the outline on pp.\ 66-67 in \cite{Gromov}): All {\em Dehn functions} of
finitely presented groups (and ``approximately all'' {\em time complexity} 
functions of {\em nondeterministic} Turing machines) are Gromov distortion 
functions of finitely generated subgroups of \ FG$_2 \times $FG$_2$; here,
FG$_2$ denotes the 2-generated free group. 
Moreover, in \cite{BiThomps} it was proved that FG$_2 \times $FG$_2$ is
embeddable with {\em linear} distortion in the Thompson group $G_{2,1}$. 
So the theorem of Ol$'$shanskii and Sapir also holds for the finitely 
generated subgroups of $G_{2,1}$.

\smallskip

Actually, Gromov \cite{Gromov} and Bourgain \cite{Bourgain} defined the 
distortion to be 
 \ $\frac{1}{n} \cdot {\rm max}\{ |g|_{\Gamma_G} : $
           $ |g|_{\Gamma_H} \leq n, \ g \in G\}$, i.e., they use an extra 
factor $\frac{1}{n}$. 
%Compare this to the induced operator 
%norm, where for a linear operator $A$ on a space ${\mathcal H}$ we have
% \ $\|A\| = {\rm sup}\{ \|Ax\| : \|x\| =1, \ x \in {\mathcal H} \} = $
%$\frac{1}{n} \cdot {\rm sup}\{\|Ax\| : \|x\| \leq n, \ x \in {\mathcal H}\}$.
However, the connections between distortion, the generalized word problem, 
and complexity (as we just saw, and will further see in the present paper) 
are more direct without the factor $\frac{1}{n}$.

\bigskip

\noindent
{\bf 2. Bourgain's distortion theorem:} Given a finite metric space $G$ with 
$n$ elements, the aim is to find embeddings of $G$ into a finite-dimensional 
euclidean space. The two distances of $G$ are its given distance and the 
inherited euclidean distance. In this problem the goal is to have small 
distortion, as a function of the cardinality of $G$, while also
keeping the dimension of the euclidean space small. Bourgain \cite{Bourgain} 
found a bound $O(n \log n)$ for the distortion (or ``$O(\log n)$'' in 
Bourgain's and Gromov's terminology). This is an important result. See also 
\cite{IndykMatousek}, \cite{ArzhGubaSapir}, \cite{ArzhDrutuSapir}.

\medskip

\noindent
{\bf 3. Generator distortion:} A variant of Gromov's distortion is obtained 
when $G = H$, but $\Gamma_G \subsetneqq \Gamma_H$. So here we look at the 
distorting effect of a change of generators in a given group. 
When $\Gamma_G$ and $\Gamma_H$ are both finite the generator distortion is 
linear; however, when $\Gamma_G$ is finite and $\Gamma_H$ is infinite the 
distortion becomes interesting.
E.g., for the Thompson group $G_{2,1}$ let us take $\Gamma_G$ to be any 
finite generating set, and for $\Gamma_H$ let us take 
$\Gamma_G \cup \{\tau_{i,j}: 1 \leq i < j\}$; here $\tau_{i,j}$ is the 
position transposition defined earlier. Then the generator distortion is 
exponential (see \cite{BiCoNP}). Also, the word problem of $G_{2,1}$ over any
finite generating set $\Gamma_G$ is in {\sf P}, but the word problem of 
$G_{2,1}$ over $\Gamma_G \cup \{\tau_{i,j}: 1 \leq i < j\}$ is 
{\sf coNP}-complete (see \cite{BiCoNP} and \cite{BiFact}).

\medskip

\noindent
{\bf 4. Monoids and directed distance:} \ Gromov's distortion and the 
generator distortion can be generalized to monoids. 
We repeat what we said about Gromov distortion, but $G$ and $H$ are now 
monoids, and $\Gamma_G$, respectively $\Gamma_H$, are monoid generating 
sets which are used to define monoid Cayley graphs. We will use the left 
Cayley graphs.
We assume $\Gamma_G \subseteq \Gamma_H$. In each Cayley graph there is a 
directed distance, defined by the lengths of directed paths.
The monoid $G$ now has two directed distance functions, the distance in the 
Cayley graph of $G$ itself, and the directed distance that $G$ inherits 
from its embedding into the Cayley graph of $H$. We denote the word-length
of $g \in G$ over $\Gamma_G$ by $|g|_G$; this is the minimum length of all
words over $\Gamma_G$ that represent $g$; it is also the length of a shortest
path from the identity to $g$ in the Cayley graph of $G$. Similarly, we 
denote the word-length of $h \in H$ over $\Gamma_H$ by $|h|_H$. 
The definition of the distortion becomes: \ \ \    
$\delta(n) \ = \ {\rm max}\{ |g|_G : g \in G, \ |g|_H \leq n\}$.

\medskip

\noindent {\bf 5. Schreier graphs:} Let $G$, $H$, and $F$ be groups, where
$F$ is a subgroup of $H$. Let $\Gamma_H$ be a generating set of $\Gamma_H$, 
and assume $\Gamma_H = \Gamma_H^{-1}$.
We can define the Schreier left coset graph of $H/F$ over the generating
set $\Gamma_H$, and the distance function $d_{H/F}(.,.)$ in this coset graph.
By definition, this Schreier graph has vertex set $H/F$ (i.e., the left 
cosets, of the form $h \cdot F$ with $h \in H$), and it has directed
edges of the form
 \ $h \cdot F \stackrel{\gamma}{\longrightarrow} \gamma g \cdot F$, \ for 
$h \in H$, $\gamma \in \Gamma_H$.
The graph is symmetric; for every edge as above there is an opposite edge
 \ $\gamma h \cdot F \stackrel{\gamma^{-1}}{\longrightarrow} h \cdot F$.
Because of symmetry the Schreier graph has a (symmetric) distance function
based on path length, \ $d_{H/F}(.,.): H/F \times H/F \to {\mathbb N}$. 
 
Next, assume that $G$ is embedded into $H/F$ by some injective function
$G \hookrightarrow H/F$. Such an embedding happens, e.g., if $G$ and $F$
are subgroups of $H$ such that $G \cap F = \{{\bf 1}\}$.
Indeed, in that case each coset in $H/F$ contains at most one element of $G$
(since $g_1F = g_2F$ implies $g_2^{-1}g_1 \in F \cap G = \{{\bf 1}\}$). 

The group $G$ now inherits a distance function from the path length in 
the Schreier graph of $H/F$. 
Comparing this distance with other distances in $G$ leads to distortion
functions. E.g., if the group $G$ is also embedded in a monoid $M$ with
monoid generating set $\Gamma_M$, this leads to the following distortion 
function: \ \ \     
$\delta_G(n) = {\rm max}\{d_{H/F}(F, gF): g \in G, \ |g|_M \leq n \}$.

It will turn out that for appropriate choices of $G, F, H$, 
$\Gamma_H$, and $\Gamma_M$, this last distortion is polynomially related to
the computational asymmetry function $\alpha$ of boolean permutations 
(Theorem \ref{polyn_rel_alpha_delta}). 

\medskip

\noindent {\bf 6. Asymmetry functions:} We already saw the computational 
asymmetry function of combinational circuits, and the word-length asymmetry 
function of a group embedded in a monoid. More generally, in any 
quasi-metric space $(S,d)$, where $d(.,.)$ is a directed distance function, 
an asymmetry function $A: {\mathbb R}_{\geq 0} \to {\mathbb R}_{\geq 0}$ can 
be defined by \ \  $A(n) = $
${\rm max}\{ d(x_2,x_1) : \, x_1,x_2 \in S, \ d(x_1,x_2) \leq n \}$.

This asymmetry function can also be viewed as the distortion of
$d^{{\rm rev}}$ with respect to $d$ in $S$; here $d^{{\rm rev}}$ denotes
the reverse directed distance, defined by 
$d^{{\rm rev}}(x_1,x_2) = d(x_2,x_1)$.

\medskip

\noindent {\bf 7. Other distortions:} 

\noindent
- Distortion can compare lengths of proofs (or lengths of expressions) 
in various, more or less powerful proof systems (respectively description 
languages). Distortion can also compare the duration of computations or of 
rewriting processes in various models of computation. 
Hence, many (perhaps all) notions of complexity are examples of distortion. 
Distortion is an algebraic or geometric representation (or cause) of
complexity.

\noindent
- Instead of length and distance, other measures (e.g., volumes in higher
dimension, energy, action, entropy, etc.) could be used.

%%%%%%%%%%%%%%%%%%%%%%%%%%%%

\subsection{Thompson-Higman groups and monoids}

The Thompson groups, introduced by Richard J.\ Thompson \cite{Th0,McKTh,Th},
are finitely presented infinite groups that act as bijections between certain
subsets of $\{0,1\}^*$. So, the elements of the Thompson groups are
transformations of bitstrings, and hence they are related to input-output 
maps of boolean circuits. In this subsection we define the Thompson group 
$G_{2,1}$ (also known as ``$V$''), as well as its generalization (by Graham 
Higman \cite{Hig74}) to the group $G_{k,1}$ that partially acts on $A^*$, for 
any finite alphabet $A$ of size $k \geq 2$.
We will follow the presentation of \cite{BiThomps} (see also \cite{BiFact}
and \cite{BiCoNP}); another reference is \cite{Scott}, which is also based
on string transformations but with a different terminology; the classical
references \cite{Th0,McKTh,Th,Hig74,CFP} do not describe the Thompson groups 
by transformations of finite strings.
Because of our interest in strings and in circuits, we also use 
generalizations of the Thompson groups to monoids, as introduced in
\cite{BiThomMon}.

Some preliminary definitions, all fairly standard, are needed in order to 
define the Thompson-Higman group $G_{k,1}$. First, we pick any alphabet 
$A$ of cardinality $|A| = k$. By $A^*$ we denote the set of all finite words
(or ``strings'') over $A$; the empty word $\varepsilon$ is also in $A^*$.
We denote the length of $w \in A^*$ by $|w|$ and we let $A^n$ denote the set 
of words of length $n$. We denote the concatenation of two words 
$u,v \in A^*$ by $uv$ or by $u \cdot v$; the concatenation of two subsets 
$B, C \subseteq A^*$ is defined by $BC = \{uv : u \in B, v \in C\}$.
A {\it right ideal} of $A^*$ is a subset $R \subseteq A^*$ such that
$RA^* \subseteq R$. A generating set of a right ideal $R$ is, by definition,
a set $C$ such that $R$ is equal to the intersection of all right ideals that 
contain $C$; equivalently, $C$ generates $R$ (as a right ideal) iff 
$R = CA^*$. A right ideal $R$ is called {\it essential} iff $R$ has a 
non-empty intersection with every right ideal of $A^*$.
For $u,v \in A^*$, we call $u$ a {\it prefix} of $v$ iff there exists 
$z \in A^*$ such that $uz = v$. A {\it prefix code} is a subset
$C \subseteq A^*$ such that no element of $C$ is a prefix of another element
of $C$. A prefix code $C$ over $A$ is maximal iff $C$ is not a strict 
subset of any other prefix code over $A$.
It is easy to prove that a right ideal $R$ has a unique minimal (under
inclusion) generating set $C_R$, and that $C_R$ is a prefix code; moreover, 
$C_R$ is a maximal prefix code iff $R$ is an essential right ideal.

For a partial function $f: A^* \to A^*$ we denote the domain by Dom$(f)$ 
and the image (range) by Im$(f)$. A restriction of $f$ is any partial 
function $f_1: A^* \to A^*$ such that Dom$(f_1) \subseteq$ Dom$(f)$, and such 
that $f_1(x) = f(x)$ for all $x \in {\rm Dom}(f_1)$. An extension of $f$ is 
any partial function of which $f$ is a restriction.
An {\it isomorphism} between right ideals $R_1, R_2$ of $A^*$ is a bijection
$\varphi: R_1 \to R_2$ such that for all $r_1 \in R_1$ and all $z \in A^*$:
 \ $\varphi(r_1z) = \varphi(r_1) \cdot z$. The isomorphism $\varphi$ is
uniquely determined by a bijection between the prefix codes that minimally 
generate $R_1$, respectively $R_2$.  One can prove \cite{Th,Scott,BiThomps} 
that every isomorphism $\varphi$ between essential right ideals has a 
{\it unique maximal extension} (within the category of isomorphisms between 
essential right ideals of $A^*$); we denote this unique maximal extension by 
max$(\varphi)$. 

Now, finally, we define the {\it Thompson-Higman group} $G_{k,1}$: It 
consists of all maximally extended isomorphisms between finitely generated 
essential right ideals of $A^*$.
The multiplication consists of composition followed by maximum extension:
 \ $\varphi \cdot \psi = $ max$(\varphi \circ \psi)$.
Note that $G_{k,1}$ acts partially and faithfully on $A^*$ on the {\it left}.

Every element $\varphi \in G_{k,1}$ can be described by a bijection between
two finite maximal prefix codes; this bijection can be described concretely
by a finite function {\it table}. When $\varphi$ is described
by a maximally extended isomorphism between essential right ideals,
$\varphi: R_1 \to R_2$, we call the minimum generating set of $R_1$ the
domain code of $\varphi$, and denote it by domC$(\varphi)$; similarly, the 
minimum generating set of $R_2$ is called the image code of $\varphi$, 
denoted by imC$(\varphi)$.

Thompson and Higman proved that $G_{k,1}$ is finitely presented. Also, when
$k$ is even $G_{k,1}$ is a simple group, and when $k$ is odd $G_{k,1}$ has a 
simple normal subgroup of index 2. In \cite{BiThomps} it was proved that the 
word problem of $G_{k,1}$ over any finite generating set is in {\sf P} (in 
fact, more strongly, in the parallel complexity class {\sf AC}$_1$).
In \cite{BiFact,BiCoNP} it was proved that the word problem of $G_{k,1}$
over $\Gamma \cup \{\tau_{i,j}: 1 \leq i < j\}$ is {\sf coNP}-complete,
where $\Gamma$ is any finite generating set of $G_{k,1}$, and where 
$\tau_{i,j}$ is the position transposition introduced in Subsection 1.1.

Because of connections with circuits we consider the subgroup ${\it
lp}G_{k,1}$ of all length-preserving elements of $G_{k,1}$; more precisely,
 \  ${\it lp}G_{k,1} = \{ \varphi \in G_{k,1} : $
$\forall x \in {\rm Dom}(\varphi), \, |x| = |\varphi(x)| \}$.
See \cite{BiFact} for a study of ${\it lp}G_{k,1}$ and some of its properties.
In particular, it was proved that ${\it lp}G_{k,1}$ is a direct limit of finite 
alternating groups, and that ${\it lp}G_{2,1}$ is generated by the set 
 \ $\{N, C, T\} \, \cup \, \{\tau_{i,i+1} : 1 \leq i \}$, where 
 \ $N: x_1w \mapsto \overline{x}_1w$, \
$C: x_1x_2w \mapsto x_1 \, (x_2 \oplus x_1) \, w$, \ and \
$T: x_1x_2x_3w \mapsto x_1x_2 \, (x_3 \oplus (x_2 \wedge x_1)) \, w$ \ (for 
$x_1, x_2, x_3 \in \{0,1\}$ and $w \in \{0,1\}^*$). Thus (recalling 
Subsection 1.4), $N, C, T$ are the {\sf not, c-not, cc-not} gates, applied 
to the first (left-most) bits of a binary string. It is known that the 
gates {\sf not, c-not, cc-not}, together with the wire-swappings, form a 
complete set of gates for bijective circuits (see 
\cite{Shende, Toff80Memo,FredToff}); hence, ${\it lp}G_{2,1}$ is closely 
related to the field of reversible computing.

\medskip

It is natural to generalize the bijections between finite maximal prefix 
codes to functions between finite prefix codes. Following \cite{BiThomMon} 
we will define below the {\it Thompson-Higman monoids} $M_{k,1}$. 
First, some preliminary definitions. A {\it right-ideal homomorphism} of 
$A^*$ is a total function $\varphi: R_1 \to A^*$ such that $R_1$ is a right 
ideal, and such that for all $r_1 \in R_1$ and all $z \in A^*$: 
 \ $\varphi(r_1z) = \varphi(r_1) \cdot z$. It is easy to prove that 
Im$(\varphi)$ is then also a right ideal of $A^*$. From now on we will write
a right-ideal homomorphism as a total surjective function 
$\varphi: R_1 \to R_2$, where both $R_1$ and $R_2$ are right ideals.  The 
homomorphism $\varphi$ is uniquely determined by a total surjective function 
$f: P_1 \to S_2$, with $P_1, S_2 \subset A^*$ where $P_1$ is the prefix code 
(not necessarily maximal) that generates $R_1$ as a right ideal, and where 
$S_2$ is a set (not necessarily a prefix code) that generates $R_2$ as a right
ideal; $f$ can be described by a finite function table. 

For two sets $X, Y$, we say that $X$ and $Y$ ``intersect'' iff 
$X \cap Y \neq \varnothing$.  We say that a right ideal $R'_1$ is 
{\em essential in} a right ideal $R_1$ iff
$R'_1$ intersects every right ideal that $R_1$ intersects. 
An {\em essential restriction} of a right-ideal homomorphism
$\varphi: R_1 \to R_2$ is a right ideal-homomorphism $\Phi: R'_1 \to R'_2$
such that $R'_1$ is essential in $R_1$, and for all $x'_1 \in R'_1$: 
$\varphi(x'_1) = \Phi(x'_1)$.
In that case we also say that $\varphi$ is an {\em essential extension} of 
$\Phi$. If $\Phi$ is an essential restriction of $\varphi$ then
$R'_2 = {\rm Im}(\Phi)$ will automatically be essential in
$R_2 = {\rm Im}(\varphi)$. Indeed, if $I$ is any no-empty right subideal of 
$R_1$ then $I \cap R'_1 \neq \varnothing$, hence \  
$\varnothing \neq \Phi(I \cap R'_1)$
$\subseteq \Phi(I) \, \cap \, \Phi(R'_1)$ $= \Phi(I) \cap R'_2$;
moreover, any non-empty right subideal $J$ of $R_2$ is of the form
$J = \Phi(I)$, where $I = \Phi^{-1}(J)$ is a non-empty right 
subideal of $R_1$; hence, for any non-empty right subideal $J$ of $R_2$, 
$\varnothing \neq J \cap R'_2$.

The free monoid $A^*$ can be pictured by its right Cayley graph, which is 
easily seen to be the infinite regular $k$-ary tree with vertex set $A^*$ 
and edge set $\{(v,va): v \in A^*, a \in A\}$. We simply call this the 
{\em tree of} $A^*$. It is a directed, rooted tree, with all paths directed 
away from the root $\varepsilon$ (the empty word); by ``path'' we will 
always mean a directed path. Many of the previously defined concepts can
be reformulated more intuitively in the context of the tree of $A^*$:  
A word $v$ is a prefix of a word $w$ iff $v$ is an ancestor of $w$ in the 
tree.
A set $P$ is a prefix code iff no two elements of $P$ are on a common path.
A set $R$ is a right ideal iff any path that starts in $R$ has all its
vertices in $R$.
The prefix code that generates $R$ consists of the elements of $R$ that are
maximal (within $R$) in the prefix order, i.e., maximally close (along paths)
to the root $\varepsilon$.
A finitely generated right ideal $R$ is essential iff every infinite path
eventually reaches $R$ (and then stays in it from there on). Similarly, a
finite prefix code $P$ is maximal iff any infinite path starting at the root
eventually intersects $P$.
For two finitely generated right ideals $R', R$ with $R' \subset R$ we have:
$R'$ is essential in $R$ iff any infinite path starting in $R$ eventually 
reaches  $R'$ (and then stays in it from there on).

Assume now that a total order $a_1 < a_2 < \ldots < a_k$ has been chosen for 
the alphabet $A$; this means that the tree of $A^*$ is now an {\em oriented} 
rooted tree, i.e., the children of each vertex $v$ have a total order
$va_1 < va_2 < \ldots < va_k$. The following can be proved (see 
\cite{BiThomMon}, Prop.\ 1.4(1)): 
$\Phi$ is an essential restriction of $\varphi$
iff $\Phi$ can be obtained from $\varphi$ by starting from the table of
$\varphi$ and applying a finite number of {\em restriction steps} of the
following form: {\it ``replace $(x,y)$ in a table by \   
$\{(xa_1,ya_1), \ldots, (xa_k,ya_k)\}$''.} 
In the tree of $A^*$ this means that $x$ and $y$ are replaced by their
children $xa_1, \ldots, xa_k$, respectively $ya_1, \ldots, ya_k$, paired 
according to the order on the children.
One can also prove (see \cite{BiThomMon}, Remark after Prop.\ 1.4):
Every right ideal homomorphism $\varphi$ with table $P \to S$ has an
essential restriction $\varphi'$ that has a table $P' \to Q'$ such that both
$P'$ and $Q'$ are prefix codes.

An important fact is the following (see \cite{BiThomMon}, Prop.\ 1.4(2)):
Every homomorphism between finitely generated right ideals 
of $A^*$ has a {\em unique maximal} essential extension; we call it the 
maximum essential extension of $\Phi$ and denote it by max($\Phi$). 

Finally here is the definition of the {\it Thompson-Higman monoid}:
 \ $M_{k,1}$ consists of all
maximum essential extensions of homomorphisms between finitely generated
right ideals of $A^*$. The multiplication is composition followed by maximum
essential extension.

One can prove the following, which implies associativity: 
For all right ideal homomorphisms $\varphi_1, \varphi_2 : $
 \ \ ${\rm max}(\varphi_2 \circ \varphi_1) = $
   ${\rm max}({\rm max}(\varphi_2) \circ \varphi_1) = $
   ${\rm max}(\varphi_2 \circ {\rm max}(\varphi_1))$.

\medskip

\noindent 
In \cite{BiThomMon} the following are proved about the Thompson-Higman
monoid $M_{k,1}$: \\  
$\bullet$ \ The Thompson-Higman group $G_{k,1}$ is the group of invertible
 elements of the monoid $M_{k,1}$. \\    
$\bullet$ \ $M_{k,1}$ is finitely generated. \\  
$\bullet$ \ The word problem of $M_{k,1}$ over any finite generating set is
  in {\sf P}.  \\  
$\bullet$ \ The word problem of $M_{k,1}$ over a generating set 
  $\Gamma \cup \{\tau_{i,j} : 1 \leq i < j \}$, where $\Gamma$ is any finite
  generating set of $M_{k,1}$, is {\sf coNP}-complete.

%%%%%%%%%%%%%%%%%%%%%%%%%%%%%%%%%%%%%%%%%%%%%%%%%%%
% Section 2
%%%%%%%%%%%%%%%%%%%%%%%%%%%%%%%%%%%%%%%%%%%%%%%%%%%

\section{Boolean functions as elements of Thompson monoids}

The input-output functions of digital circuits map bitstrings of some fixed 
length to bitstrings of a fixed length (possibly different from the input 
length). In other words, circuits have input-output maps that are total
functions of the form $f: \{0,1\}^m \to \{0,1\}^n$ for some $m,n >0$.   
The Thompson-Higman monoid $M_{k,1}$ has an interesting submonoid that 
corresponds to fixed-length maps, defined as follows.

\begin{defn} \label{lepM} {\bf (the submonoid ${\it lep}M_{k,1}$).} \
Let $\varphi: PA^* \to QA^*$ be a right-ideal homomorphism, where
$P, Q \subset A^*$ are finite prefix codes, and where $P$ is a {\em maximal}
prefix code. Then $\varphi$ is called {\em length equality preserving} iff 
for all $x_1, x_2 \in {\rm Dom}(\varphi):$ \ \ $|x_1| = |x_2|$ implies 
$|\varphi(x_1)| = |\varphi(x_2)|$.

The submonoid ${\it lep}M_{k,1}$ of $M_{k,1}$ consists of those elements of 
$M_{k,1}$ that can be represented by length-equality preserving
right-ideal homomorphisms.
\end{defn}
It is easy to check that an essential restriction of an element of
${\it lep}M_{k,1}$ is again in ${\it lep}M_{k,1}$, so ${\it lep}M_{k,1}$ is
well defined as a subset of $M_{k,1}$; moreover, one can easily check that
${\it lep}M_{k,1}$ is closed under composition, so ${\it lep}M_{k,1}$ is
indeed a submonoid of $M_{k,1}$.

For $\varphi \in M_{k,1}$ we have $\varphi \in {\it lep}M_{k,1}$ iff there 
exist $m > 0$ and $n>0$ such that $A^m \subset {\rm Dom}(\varphi)$ and 
$\varphi(A^m) \subseteq A^n$. So (by means of an essential restriction, if 
necessary), $\varphi$ can be represented by a function table 
$A^m \to Q \subseteq A^n$ with a {\em fixed input length and a fixed output 
length} (but the input and output lengths can be different).

\smallskip

The motivation for studying the monoid ${\it lep}M_{k,1}$ is the following.
Every boolean function $f: \{0,1\}^m \to \{0,1\}^n$ (for any $m,n >0$)
determines an element of ${\it lep}M_{k,1}$, and conversely, this element of
${\it lep}M_{k,1}$ determines $f$ when restricted to $\{0,1\}^m$. By 
considering all boolean functions as elements of ${\it lep}M_{k,1}$ we gain
the ability to compose arbitrary boolean functions, even if their domain and
range ``do not match''. Moreover, in ${\it lep}M_{k,1}$ we are able to 
generate all boolean functions from gates by using ordinary {\em functional 
composition} (instead of graph-based circuit lay-outs). The following 
remains open: 
 
\medskip

{\bf Question:} Is ${\it lep}M_{k,1}$ finitely generated?

\medskip

\noindent However we can find nice infinite generating sets, in connection
with circuits. 

\begin{pro} {\bf (Generators of ${\it lep}M_{k,1}$).} \label{gen_tau_lepM} \
The monoid ${\it lep}M_{k,1}$ has a generating set of the form
 \ $\Gamma \cup \, \{\tau_{i,i+1} : 1 \leq i \}$, for some finite subset
 \ $\Gamma \subset {\it lep}M_{k,1}$.
\end{pro}
{\bf Proof.} We only prove the result for $k=2$; a similar reasoning works 
for all $k$ (using $k$-ary logic). 

It is a classical fact that any function $f: \{0,1\}^m \to \{0,1\}^n$ can be 
implemented by a combinational circuit that uses copies of {\sf and}, 
{\sf or}, {\sf not}, {\sf fork} and wire-crossings. So all we need to do is 
to express theses gates, at any place in the circuit, by a finite subset of 
${\it lep}M_{2,1}$ and by positions transpositions $\tau_{i,i+1}$.  
For each gate $g \in$ \{{\sf and, or}\} we define an element 
$\gamma_g \in {\it lep}M_{k,1}$ by

\smallskip 

$\gamma_g: \ x_1x_2 w \ \in \ \{0,1\}^m \ \longmapsto \ $
  $g(x_1,x_2) \ w \ \in \ \{0,1\}^{m-1}$.

\smallskip

\noindent Similarly we define 
 $\gamma_{\sf not}, \gamma_{\sf fork} \in {\it lep}M_{k,1}$ by

\smallskip 

$\gamma_{\sf not}: \ x_1w \ \in \ \{0,1\}^m \ \longmapsto \ $
  $\overline{x_1} \ w \ \in \ \{0,1\}^m$,

\smallskip

$\gamma_{\sf fork}: \ x_1w \ \in \ \{0,1\}^m \ \longmapsto \ $
  $x_1 \, x_1 \, w \ \in \ \{0,1\}^{m+1}$. 
  
\smallskip 

\noindent For each $g \in \{ {\sf and, or, not, fork} \}$, \ $\gamma_g$ 
transforms only the first one or two boolean variables, and leaves the 
other boolean variables unchanged. We also need to simulate the effect of a 
gate $g$ on any variable $x_i$ or pair of variables $x_ix_{i+1}$, i.e., 
we need to construct the map 

\smallskip

$ux_ix_{i+1}v \ \in \ \{0,1\}^m \ \longmapsto \ $
 $u \ g(x_i,x_{i+1}) \ v \ \in \ \{0,1\}^{m-1}$ 

\smallskip

\noindent (and similarly in case where $g$ is {\sf not} or {\sf fork}). 
For this, we apply wire-transpositions to move $x_ix_{i+1}$ to the 
wire-positions 1 and 2, then we apply $\gamma_g$, then we apply more 
wire-transpositions in order to move $g(x_1,x_2)$ back to position $i$. 
Thus the effect of any gate anywhere in the circuit can be expressed as a 
composition of $\gamma_g$ and position transpositions in 
$\{\tau_{i,i+1}: 1 \leq i\}$. 
 \ \ \ $\Box$

\begin{pro} {\bf (Change of generators of ${\it lep}M_{k,1}$).} 
\label{(Change_gen_tau_lepM} \  
Let $\{\tau_{i,i+1} : 1 \leq i \}$ be denoted by $\tau$.
If $\Gamma, \Gamma' \subset {\it lep}M_{k,1}$ are two finite sets such that
$\Gamma \cup \tau$ and $\Gamma' \cup \tau$ generate ${\it lep}M_{k,1}$,
then the word-length over $\Gamma \cup \tau$ is linearly related to the 
word-length over $\Gamma' \cup \tau$. In other words, there are
constants $c' \geq c \geq 1$ such that for all $m \in {\it lep}M_{k,1}:$
 \ \ \    
$|m|_{\Gamma \cup \tau} \ \leq \ c \cdot |m|_{\Gamma' \cup \tau} \ \leq \ $
$c' \cdot |m|_{\Gamma \cup \tau}$.   
\end{pro}
{\bf Proof.} Since $\Gamma$ is finite, the elements of $\Gamma$ can be 
expressed by a finite set of words of bounded length ($ \leq c$) over 
$\Gamma' \cup \tau$. Thus, every word of length $n$ over $\Gamma \cup \tau$ 
is equivalent to a word of length $\leq c \, n$ over $\Gamma' \cup \tau$.
This proves the first inequality.
A similar reasoning proves the second inequality. 
 \ \ \ $\Box$

%%%% Proposition
\begin{pro} {\bf (Circuit size vs.\ ${\it lep}M_{2,1}$ word-length).} 
\label{circ_lepM} \\     
Let \ $\Gamma_{{\it lep}M_{2,1}} \cup \, \{\tau_{i,j} : 1 \leq i < j \}$ 
 \ be a generating set of 
${\it lep}M_{2,1}$ with $\Gamma_{{\it lep}M_{2,1}}$ finite.  
Let $f: \{0,1\}^m \to Q \ \ (\subseteq \{0,1\}^n)$ be a function
defining an element of ${\it lep}M_{2,1}$,
and let $|f|_{{\it lep}M_{2,1}}$ the word-length of $f$ over the generating
set \ $\Gamma_{{\it lep}M_{2,1}} \cup \, \{\tau_{i,j} : 1 \leq i < j \}$. 
Let $|C_f|$ be the circuit size of $f$ (using any finite universal set of 
gates and wire-swappings).
Then $|f| _{{\it lep}M_{2,1}}$ and $|C_f|$ are linearly related. 
More precisely, for some constants $c_1 \geq c_o \geq 1:$ 

\smallskip

 \ \ \  \ \ \ \ \ $|C_f| \ \ \leq \ \ c_o \cdot |f| _{{\it lep}M_{2,1}}$
  $ \ \ \leq \ \ c_1 \cdot |C_f|$.
\end{pro}
{\bf Proof.} For the proof we assume that the set of gates for circuits
(not counting the wire-transpositions) is $\Gamma_{{\it lep}M_{2,1}}$. 
If we make a different choice for the universal set of gates for circuits, 
and a different choice for the finite portion $\Gamma_{{\it lep}M_{2,1}}$ of 
the generating set of ${\it lep}M_{2,1}$ then the inequalities remain the 
same, except for the constants $c_1, c_o$. 

\smallskip

The inequality \ $|C_f| \leq |f| _{{\it lep}M_{2,1}}$ \ is  
obvious, since a word $w$ over 
 \ $\Gamma_{{\it lep}M_{2,1}} \cup \, \{\tau_{i,j} : 1 \leq i < j\}$ \ is
automatically a circuit of size $|w|$. 

\smallskip

For the other inequality, we want to simulate each gate of the circuit $C_f$
by a word over 
$\Gamma_{{\it lep}M_{2,1}} \cup \, \{\tau_{i,j} : 1 \leq i < j\}$. The
reasoning is the same for every gate, so let us just focus on an {\sf or}
gate. The essential difference between circuit gates and elements of 
${\it lep}M_{2,1}$ is that in a circuit, a gate (with 2 input wires,
for example) can be applied to any two wires in the circuit; on the other 
hand, the functions in ${\it lep}M_{2,1}$ are applied to the first 
few wires. However, the circuit gate {\sf or}, applied to $(i, i+1)$
can be simulated by an element of $\Gamma_{{\it lep}M_{2,1}}$ and a few wire 
transpositions, since we have: \ \  ${\sf or}_{i,i+1}(.) = $
$\gamma_{{\sf or}} \circ \tau_{2,i+1} \circ \tau_{1,i}(.)$.  

The output wire of ${\sf or}_{i,i+1}(.)$ is wire number $i$, whereas the
output wire of $\gamma_{{\rm or}} \circ \tau_{2,i+1} \circ \tau_{1,i}(.)$
is wire number 1. However, instead of permuting all the wires in order 
to place the output of 
 \ $\gamma_{{\sf or}} \, \tau_{2,i+1} \, \tau_{1,i}(.)$ \ on wire $i$, we
just leave the output of 
$\gamma_{{\sf or}} \, \tau_{2,i+1} \, \tau_{1,i}(.)$ \ on wire 1 for now.
The simulation of the next gate will then use appropriate transpositions
$\tau_{2,j} \cdot \tau_{1,k}$ for fetch the correct input wires for the next 
gate. Thus, each gate of $C_f$ is simulated by one function in 
$\Gamma_{{\it lep}M_{2,1}}$ and a bounded number of wire-transpositions in \   
$\{\tau_{i,j} : 1 \leq i < j\}$. 

At the output end of the circuit, a permutation of the $n$ output wires is 
needed in order to send the outputs to the correct wires; any permutation of
$n$ elements can be realized with $< n$ ($\leq |C_f|$) transpositions. 
(The inequality $n \leq |C_f|$ holds because since we count the output ports 
in the circuit size.)
 \ \ \ $\Box$

\bigskip

\noindent
{\bf Remark.} The above Proposition motivates our choice of generating set
of the form $\Gamma \cup \{\tau_{i,j} : 1 \leq i < j\}$ (with $\Gamma$
finite) for ${\it lep}M_{k,1}$; in particular, it motivates the inclusion of 
all the position transpositions $\tau_{i,j}$ in the generating set. The
Proposition also motivates the definition of word-length in which 
$\tau_{i,j}$ has word-length 1 for all $j > i \geq 1$. 

\bigskip

\noindent
Next we will study the {\bf distortion} of ${\it lep}M_{k,1}$ in $M_{k,1}$. 
We first need some Lemmas. 

\begin{lem}
\label{idealInter}  {\rm (Lemma 3.3 in \cite{BiThomps}).} \
If $P, Q, R \subseteq A^*$ are such that $PA^* \cap QA^* = RA^*$ and
$R$ is a prefix code, then $R \subseteq P \cup Q$.
\end{lem}
{\bf Proof.} For any $r \in R$ there are $p \in P, q \in Q$ and  
$v, w \in A^*$ such that $r = pv = qw$. Hence $p$ is a prefix of $q$ or $q$
is a prefix of $p$. Let us assume $p$ is a prefix of $q = px$, for some 
$x \in A^*$ (the other case is similar) 
Hence $q = px \in PA^* \cap QA^* = RA^*$, and $q$ is a
prefix of $r = qw$. Since $R$ is a prefix code, $r = q$, hence $r \in Q$.
 \ \ \ $\Box$

\begin{lem} \label{InvIm_of_rightideal} \
Let $P,Q \subset A^*$ be finite prefix codes, and let
$\theta: PA^* \to QA^*$ be a right-ideal homomorphism with domain $PA^*$
and image $QA^*$.  Let $S$ be a prefix code with $S \subset QA^*$. 
Then $\theta^{-1}(S)$ is a prefix code and 
 \ $\theta^{-1}(SA^*) = \theta^{-1}(S) \ A^*$.
\end{lem}
{\bf Proof.}  First, $\theta^{-1}(S)$ is a prefix code. Indeed, if we had
$x_1 = x_2u$ for some $x_1, x_2 \in \theta^{-1}(S)$ with $u$ non-empty, then
$\theta(x_1) = \theta(x_2) \ u$. This would contradict the assumption that
$S$ is a prefix code.

Second, $\theta^{-1}(S) \subset \theta^{-1}(SA^*)$, hence
 \ $\theta^{-1}(S) \ A^* \subseteq \theta^{-1}(SA^*)$, since 
$\theta^{-1}(SA^*)$ is a right ideal. (Recall that the inverse image of a
right ideal under a right-ideal homomorphism is a right ideal.)

We also want to show that $\theta^{-1}(SA^*) \subseteq \theta^{-1}(S) \ A^*$. 
Let  $x \in \theta^{-1}(SA^*)$. So, $\theta(x) = sv$ for some $s \in S$,
$v \in A^*$, and $s = qu$ for some $q \in Q$, $u \in A^*$. Since 
$\theta(x) = quv$, we have $x = puv$ for some $p \in P$ with 
$\theta(p) = q$. Hence $\theta(pu) = qu = s$. Therefore, $x = puv$ with
$pu \in \theta^{-1}(s) \subseteq \theta^{-1}(S)$, hence 
$x \in \theta^{-1}(S) \ A^*$.   
  \ \ \ $\Box$

\bigskip

\noindent {\bf Notation:}
For a right-ideal homomorphism 
$\varphi: {\rm Dom}(\varphi) = PA^* \to {\rm Im}(\varphi) = QA^*$, where
$P,Q \subset A^*$ are finite prefix codes, we define

\smallskip

 \ $\ell(\varphi) \ = \ {\rm max}\{|z| : z \in P \cup Q\}$,

\smallskip

\noindent For any finite prefix code $C \subset A^*$ we define

\smallskip

 \ $\ell(C) \ = \ {\rm max}\{|z| : z \in C\}$.

\begin{lem} \label{max_length_phi_P} \   
Let $\varphi: {\rm Dom}(\varphi) = PA^* \to {\rm Im}(\varphi) = QA^*$ be a 
right-ideal homomorphism, where $P$ and $Q$ are finite prefix codes. Let 
$R \subset A^*$ be any finite prefix code. Then we have:

\smallskip

\noindent {\rm (1)} 
 \ \ \ \ \ $\ell(\varphi^{-1}(R)) < \ell(\varphi) + \ell(R)$, 

\smallskip

\noindent {\rm (2)} \ \ \ \ \ $\ell(\varphi(R)) < \ell(\varphi) + \ell(R)$. 
\end{lem}
{\bf Proof.} (1) \ Let $r \in R \cap {\rm Im}(\varphi)$. Then every element 
of $\varphi^{-1}(r)$ has the form $p_1w$ for some $p_1 \in P$ and 
$w \in A^*$ such that $r = q_1w$ for some $q_1 \in Q$ 
(with $\varphi(p_1) = q_1$).
Hence $|p_1w| = |p_1| + |r| - |q_1| = |r| + |p_1| - |q_1|$. Moreover,
$|r| \leq \ell(R)$ and $|p_1| - |q_1| < \ell(\varphi)$, so 
$|p_1 w| < \ell(R) + \ell(\varphi)$.

(2) \ If $r \in R \cap {\rm Dom}(\varphi)$ then $\varphi(r)$ has the form 
$q_1v$ for some $q_1 \in Q$ and $v \in A^*$ such that $r = p_1w$ for some 
$p_1 \in P$ (with $\varphi(p_1) = q_1$). Hence $|q_1v| = |q_1| + |r| - |p_1|$
$= |r| + |q_1| - |p_1|$. Moreover, $|r| \leq \ell(R)$  and 
$|q_1| - |p_1| < \ell(\varphi)$, so $|q_1 w| < \ell(R) + \ell(\varphi)$.
  \ \ \ $\Box$

\bigskip

\noindent For any right-ideal homomorphisms $\varphi_i$ 
(with $i = 1, \ldots, N$), the composite map
 \ $\varphi_N \circ \ldots \circ \varphi_1(.)$ \ is a right-ideal
homomorphism. We say that right-ideal homomorphisms $\Phi_i$ (with 
$i = 1, \ldots, N$) are {\em directly composable} iff 
${\rm Dom}(\Phi_{i+1}) = {\rm Im}(\Phi_i)$, for $i = 1, \ldots, N-1$.
The next Lemma shows that we can replace composition by direct composition.

%%%%%%%%%%%%%%%%%
\begin{lem} \label{max_length_factors} \  
Let $\varphi_i: {\rm Dom}(\varphi_i) = P_iA^* \to $
 ${\rm Im}(\varphi_i) = Q_iA^*$ be a right-ideal homomorphism (for 
$i = 1, \ldots, N$), where $P_i$ and $Q_i$ are finite prefix codes. 
Then each $\varphi_i$ has a (not necessarily essential) restriction to a 
right-ideal homomorphism $\Phi_i$ with the following properties:

\smallskip

\noindent $\bullet$ \ \ $\Phi_N \circ \ldots \circ \Phi_1(.) \ = \ $
          $\varphi_N \circ \ldots \circ \varphi_1(.)$; 

\smallskip

\noindent $\bullet$ \ \ ${\rm Dom}(\Phi_{i+1}) = {\rm Im}(\Phi_i)$, 
     \ for $i = 1, \ldots, N-1$;

\smallskip

\noindent $\bullet$ 
 \ \ $\ell(\Phi_i) \ \leq \ \sum_{j=1}^N \ell(\varphi_j)$ \ for every 
$i = 1, \ldots, N$. 
\end{lem}
{\bf Proof.} We use induction on $N$. For $N=1$ there is nothing to prove.
So we let $N>1$ and we assume that the Lemma holds for 
$\varphi_i: P_iA^* \to Q_iA^*$ with $i = 2, \ldots, N$, i.e., we assume 
that each $\varphi_i$ (for $i = 2, \ldots, N$) has a restriction
$\varphi'_i: P'_i A^* \to Q'_i A^*$ such that 
 \ $\varphi'_N \circ \ldots \circ \varphi'_2 \ = \ $
$ \varphi_N \circ \ldots \circ \varphi_2$, \   
$P'_{i+1} = Q'_i$ (for $i = 2, \ldots, N-1)$, 
 \ and \ $\ell(\varphi'_i) \ \leq \ \sum_{j=2}^N \ell(\varphi_j)$ \ for every
$i = 2, \ldots, N$. From $P'_{i+1} = Q'_i$ (for $i = 2, \ldots, N-1)$
it follows that \ $\ell(\varphi'_N \circ \ldots \circ \varphi'_2) \ \leq \ $
${\rm max}\{ \ell(\varphi'_i) : i = 2, \ldots, N\}$
$ \ \leq \ \sum_{j=2}^N \ell(\varphi_j)$.

Using the notation $\varphi'_{[N,2]}$ for 
$\varphi'_N \circ \ldots \circ \varphi'_2$ we have 
${\rm Dom}(\varphi'_{[N,2]}) =  P_2A^*$ and 
${\rm Im}(\varphi'_{[N,2]}) = Q_NA^*$. When we compose $\varphi_1$ and
$\varphi'_{[N,2]}$ we obtain 

\smallskip

 \ \ \ \ \ 
$\varphi_1^{-1}(Q_1A^* \cap P_2A^*) \ \ \stackrel{\Phi_1}{\longrightarrow}$
$ \ \ Q_1A^* \cap P_2A^* \ \stackrel{\Phi'_{[N,2]}}{\longrightarrow} $
$ \ \ \varphi'_{[N,2]}(Q_1A^* \cap P_2A^*)$.

\smallskip

\noindent In this diagram, $\Phi_1$ is the restriction of $\varphi_1$
to the domain \ $\varphi_1^{-1}(Q_1A^* \cap P_2A^*)$ \ and image
\ $Q_1A^* \cap P_2A^*$; and $\Phi'_{[N,2]}$ is the restriction of 
$\varphi'_{[N,2]}$ to the domain \ $Q_1A^* \cap P_2A^*$ \ and image
 \ $\varphi'_{[N,2]}(Q_1A^* \cap P_2A^*)$. Hence, 
 \ $\Phi'_{[N,2]} \circ \Phi_1 = \varphi'_{[N,2]} \circ \varphi_1$, \ and
 \ ${\rm Dom}(\Phi'_{[N,2]}) = {\rm Im}(\Phi_1)$ \ ($= Q_1A^* \cap P_2A^*$).  
So $\Phi_1$ and $\Phi'_{[N,2]}$ are directly composable.

By Lemma \ref{idealInter} there is a prefix code $S \subset A^*$
such that \ $SA^* = Q_1A^* \cap P_2A^*$ \ and \ $S \subseteq Q_1 \cup P_2$.
Hence, $\ell(S) \leq {\rm max}\{\ell(Q_1), \ell(P_2)\}$
$\leq {\rm max}\{\ell(\varphi_1), \ell(\varphi'_2)\} \leq $
${\rm max}\{\ell(\varphi_1), \sum_{j=2}^N \ell(\varphi_j) \}$
$\leq \sum_{j=1}^N \ell(\varphi_j)$. 

It follows also that 
 \ $\varphi_1^{-1}(Q_1A^* \cap P_2A^*) = \varphi_1^{-1}(SA^*) = $
$\varphi_1^{-1}(S) \ A^*$ (the latter equality is from Lemma 
\ref{InvIm_of_rightideal}). Since $S \subseteq Q_1 \cup P_2$ implies
$\varphi_1^{-1}(S) \subseteq \varphi_1^{-1}(Q_1) \cup \varphi_1^{-1}(P_2)$
$= P_1 \cup \varphi_1^{-1}(P_2)$, we have \ $\ell(\varphi_1^{-1}(S)) \leq $
 $ {\rm max}\{\ell(P_1), \ \ell(\varphi_1^{-1}(P_2))\}$.
Obviously, $\ell(P_1) \leq \ell(\varphi_1)$. 
Moreover, by Lemma \ref{max_length_phi_P}, 
 \ $\ell(\varphi_1^{-1}(P_2)) \leq \ell(\varphi_1) + \ell(P_2)$. 
Since $\ell(P_2) \leq \ell(\varphi'_2)$ 
$\leq \sum_{j=2}^N \ell(\varphi_j)$ (the latter ``$\leq$'' by induction),
we have \ $\ell(\varphi_1^{-1}(S)) \leq \ell(\varphi_1) + $
   $\sum_{j=2}^N \ell(\varphi_j) \ = \ \sum_{j=1}^N \ell(\varphi_j)$.

Since the domain code of $\Phi_1$ is $\varphi_1^{-1}(S)$ and its image code 
is $S$, we conclude that
 \ \ $\ell(\Phi_1) \leq  \sum_{j=1}^N \ell(\varphi_j)$.

\smallskip

Let us now consider any $\Phi'_{[i,2]}$, for $i = 1, \ldots, N$. 
By definition, $\Phi'_{[i,2]}$ is the restriction of 
$\varphi'_i \circ \ldots \circ \varphi'_2$ to the domain $SA^*$. So the domain 
code of $\Phi'_{[i,2]}$ is $S$, and we just proved that 
 \ $\ell(S) \leq \ \sum_{j=1}^N \ell(\varphi_j)$.  The image code of 
$\Phi'_{[i,2]}$ is \ $\varphi'_i \circ \ldots \circ \varphi'_2(S)$.
Since $S \subseteq Q_1 \cup P_2$ we have \ 

\smallskip

$\varphi'_i \circ \ldots \circ \varphi'_2(S) \ \ \subseteq \ \ $
$\varphi'_i \circ \ldots \circ \varphi'_2(Q_1)$  $ \ \cup \ $
$\varphi'_i \circ \ldots \circ \varphi'_2(P_2) \ \ = \ \ $ 
$\varphi'_i \circ \ldots \circ \varphi'_2(Q_1) \ \cup \ Q'_i$. 

\smallskip

\noindent Therefore:
 \ $\ell(\varphi'_i \circ \ldots \circ \varphi'_2(S)) \ \leq \ $
${\rm max}\{ \ell(\varphi'_i \circ \ldots \circ \varphi'_2(Q_1)), $
$ \ \ell(Q'_i) \}$.

We have \ $\ell(Q'_i) \leq \ell(\varphi'_i) \leq \sum_{j=2}^N \ell(\varphi_j)$
 \ (the last ``$\leq$'' by induction).

By Lemma \ref{max_length_phi_P}, 
 \ $\ell(\varphi'_i \circ \ldots \circ \varphi'_2(Q_1)) \ \leq \ $
   $\ell(\varphi'_i \circ \ldots \circ \varphi'_2) + \ell(Q_1) \ \leq \ $
   $\ell(\varphi'_i \circ \ldots \circ \varphi'_2) + \ell(\varphi_1)$. 
And \ $\ell(\varphi'_i \circ \ldots \circ \varphi'_2) \ \leq \ $
  ${\rm max}\{ \ell(\varphi'_j) : j = 2, \ldots, i\}$, because 
${\rm Dom}(\varphi'_{r+1}) = {\rm Im}(\varphi'_r)$ for all 
$r = 2, \ldots, N-1$.
And by induction, \ $\ell(\varphi'_j) \leq \sum_{j=2}^N \ell(\varphi_j)$.
Hence, $\ell(\varphi'_i \circ \ldots \circ \varphi'_2(Q_1)) \ \leq \ $
$\sum_{j=1}^N \ell(\varphi_j)$.

Thus, \ $\ell(\Phi'_{[i,2]}) \ \leq \ \sum_{j=1}^N \ell(\varphi_j)$ 
 \ for every $i = 2, \ldots, N$.

\smallskip

Finally, we factor $\Phi'_{[N,2]}$ as 
 \ $\Phi'_{[N,2]} \ = \ \Phi_N \circ \ldots \circ \Phi_2$,
where $\Phi_i$ (for $i = 2, \ldots, N$) is defined to be the restriction of
$\varphi'_i$ to the domain 
 \ $\varphi'_{i-1} \circ \ldots \circ \varphi'_2(SA^*)$ 
 \ ($= \Phi'_{[i-1,2]}(SA^*)$).
Since \ ${\rm Dom}(\varphi'_{r+1}) = {\rm Im}(\varphi'_r)$ \ (for 
all $r = 2, \ldots, N-1$), the domain of $\varphi'_i$ is equal to the image 
of \ $\varphi'_{i-1} \circ \ldots \circ \varphi'_2$.
So, the domain code of $\Phi_i$ is 
 \ $\varphi'_{i-1} \circ \ldots \circ \varphi'_2(S)$, and its image code is
 \ $\varphi'_i \circ \varphi'_{i-1} \circ \ldots \circ \varphi'_2(S)$. 
Since we already proved that 
 \ $\ell(\varphi'_i \circ \ldots \circ \varphi'_2(S)) \ \leq \ $
 $\sum_{j=1}^N \ell(\varphi_j)$ \ (for all $i$), it follows that
 \ $\ell(\Phi_i) \ \leq \  \sum_{j=1}^N \ell(\varphi_j)$.
 \ \ \ $\Box$

\bigskip

In the next theorem we show that the distortion of ${\it lep}M_{k,1}$ in
$M_{k,1}$ is at most quadratic (over the generators considered so far, which
include the bit position transpositions). Combined with Proposition 
\ref{circ_lepM}, this means the following: 

{\it  Assume circuits are built with gates that are not constrained to have 
fixed-length inputs and outputs, but assume the input-output function has
fixed-length inputs and outputs. Then the resulting circuits are not much 
more compact than conventional circuits, built from gates that have
fixed-length inputs and outputs (we gain at most a square-root in size). }

% From an engineering point of view, non-fixed-length inputs and outputs 
%make sense, e.g., in the case of sequential input and output.
% (compare with the phone number system in many 
%countries, where phone numbers do not have a fixed length but form a prefix 
%code). 
 
%%%%%%%%  

\begin{thm} {\bf (Distortion of ${\it lep}M_{k,1}$ in $M_{k,1}$).} 
 \label{distor_lepM_M} \  
The word-length (or Cayley graph) distortion of ${\it lep}M_{k,1}$ in 
$M_{k,1}$ has a quadratic upper bound; in other words, for all 
$x \in {\it lep}M_{k,1}$:

\smallskip
 
\hspace{2in} $|x|_{{\it lep}M_{k,1}} \ \leq \ c \cdot (|x|_{M_{k,1}})^2$

\smallskip

\noindent where $c\geq 1$ is a constant. Here the generating sets used are    
 \ $\Gamma_{M_{k,1}} \cup \ \{\tau_{i,j}: 1 \leq i<j\}$ \ for $M_{k,1}$, 
and \ $\Gamma_{{\it lep}M_{k,1}} \cup \ \{\tau_{i,j} : 1 \leq i<j\}$ 
 \ for ${\it lep}M_{k,1}$,
where $\Gamma_{M_{k,1}}$ and $\Gamma_{{\it lep}M_{k,1}}$ are finite.
By $|x|_{M_{k,1}}$ and $|x|_{{\it lep}M_{k,1}}$ we denote the word-length of 
$x$ over $\Gamma_{M_{k,1}} \cup \ \{\tau_{i,j}: 1 \leq i<j\}$,
respectively $\Gamma_{{\it lep}M_{k,1}} \cup \ \{\tau_{i,j} : 1 \leq i<j\}$.
\end{thm}
{\bf Proof.} We only prove the result for $k=2$; a similar proof applies 
for any $k$. 
We abbreviate the set $\{\tau_{i,j} : 1 \leq i <j\}$ by $\tau$. The choice 
of the finite sets $\Gamma_{M_{k,1}}$ and $\Gamma_{{\it lep}M_{k,1}}$ does 
not matter (it only affects the constant $c$ in the Theorem.  
By Corollary 3.6 in \cite{BiThomMon} we can choose $\Gamma_{M_{k,1}}$ so 
that each $\gamma \in \Gamma_{M_{k,1}}$ satisfies the following (recall that
$\ell(S)$ denotes the length of the longest words in a set $S$):

\smallskip
 
 \ \ \ \ \   
 $\ell\big({\rm domC}(\gamma) \ \cup \ {\rm imC}(\gamma)\big) \ \leq \ 2$, 
 \ \ \ and
 
\smallskip

 \ \ \ \ \ $\big| |\gamma(x)| - |x| \big| \leq 1$ \ 
for all $x \in {\rm Dom}(\gamma)$.

\smallskip
 
Let $\varphi \in {\it lep}M_{k,1}$, and let
 \ $w = \alpha_N \ldots \alpha_1$ \ be a shortest word over the generating
set  $\Gamma_{M_{k,1}} \cup \tau$ of $M_{k,1}$, representing $\varphi$. 
So $N = |\varphi|_{M_{k,1}}$. 
We restrict each partial function $\alpha_i$ to a partial function 
$\alpha'_i$ such that 
 \ ${\rm imC}(\alpha'_i) = {\rm domC}(\alpha'_{i+1})$ \ for 
$i = 1, \ldots, N-1$, according to Lemma \ref{max_length_factors}. Hence,
 \ $\alpha_N \circ \ldots \circ \alpha_1(.) \ = \ $
$\alpha'_N \circ \ldots \circ \alpha'_1(.)$, and 
 \ $\ell(\alpha'_i) \leq \sum_{j=1}^N \ell(\alpha_j)$ \ for every 
$i = 1, \ldots, N$.  Then 
 \ $\alpha_N \circ \ldots \circ \alpha_1(.)$ \ is a function 
 \ $\{0,1\}^m \, \{0,1\}^* \to Q \, \{0,1\}^*$, representing $\varphi$, and 
we will identify \ $\alpha_N \circ \ldots \circ \alpha_1(.)$ \ with $\varphi$. 
It follows that \ ${\rm domC}(\alpha'_1) = {\rm domC}(\varphi) = \{0,1\}^m$,
and  \ ${\rm imC}(\alpha'_N) = {\rm imC}(\varphi) = Q \subseteq \{0,1\}^n$. 
More generally, it follows that  
 \ ${\rm imC}(\alpha'_i \circ \ldots \circ \alpha'_1) = {\rm imC}(\alpha'_i)$, 
and 
 \ ${\rm domC}(\alpha'_N \circ \ldots \circ \alpha'_i) = $
 ${\rm domC}(\alpha'_i)$.

Since $\ell(\alpha'_i) \leq \sum_{j=1}^N \ell(\alpha_j)$, and 
$\ell(\alpha_j) \leq 2$ for all $j$, we have for every $i = 1, \ldots, N$:
 \ \ \ $\ell(\alpha'_i) \leq 2 \, N$. 

\medskip

  From here on we will simply denote $\ell(\alpha'_i)$ by $\ell_i$. 
Now, we will replace each $\alpha'_i \in M_{k,1}$ by 
$\beta_i \in {\it lep}M_{k,1}$, such that
 \ ${\rm domC}(\beta_i) = \{0,1\}^{\ell_i}$, and  
 \ ${\rm imC}(\beta_i) \subseteq \{0,1\}^{\ell_{i+1}}$; so $\beta_i$ is
length-equality preserving. 
This will be done by artificially lengthening those words in
${\rm domC}(\alpha'_i)$ that have length $< \ell_i$ and those words in
${\rm imC}(\alpha'_i)$ that have length $< \ell_{i+1}$. Moreover, we make 
$\beta_i$ defined on {\it all} of $\{0,1\}^{\ell_i}$. In detail, $\beta_i$ 
is defined as follows:

\medskip

\noindent $\bullet$ \ If $\ell_i \leq \ell_{i+1}$ :

\smallskip

$\beta_i(u \, z) \ = \ v \ z \ 0^{\ell_{i+1}-\ell_i -|v| + |u|}$ 
 \ \ \ for all \ $u \in {\rm domC}(\alpha'_i)$, 
 \ and \ $z \in \{0,1\}^{\ell_i - |u|}$; \ here \ $v = \alpha'_i(u)$;

\smallskip

$\beta_i(x) \ = \ x \ 0^{\ell_{i+1}-\ell_i}$         
 \ \ \ for all \ $x \not\in {\rm Dom}(\alpha'_i)$, \ $|x| = \ell_i$. 

\medskip

\noindent $\bullet$ \ If $\ell_i > \ell_{i+1}$ :

\smallskip

$\beta_i(u \, z_1 \, z_2) \ = \ v \, z_1$ \ \ \ for 
 all \ $u \in {\rm domC}(\alpha'_i)$ \ and all \  $z_1, z_2 \in \{0,1\}^*$
 \ with 

 \hspace{1in}  $|z_1| = \ell_{i+1} - |v|$, 
 \ \ $|z_2| = \ell_i - \ell_{i+1} + |v| - |u|$; \ here, $v = \alpha'_i(u)$;

\smallskip

$\beta_i(x_1 \, x_2) \ = \ x_1$ \ \ \ for all \ $x_1, x_2 \in \{0,1\}^*$ \
 such that \ $x_1 x_2 \not\in {\rm Dom}(\alpha'_i)$, \ with 

 \hspace{1in} $|x_1| = \ell_{i+1}$, \ \ $|x_2| = \ell_i - \ell_{i+1}$.

\bigskip

\noindent {\sf Claim.} 
 \ \ \ $\beta_N \circ \ldots \circ \beta_1(.) = \varphi$.  

\smallskip

\noindent Proof of the Claim: We observe first that \   
${\rm domC}(\beta_1) = {\rm domC}(\alpha'_1)$ ($ = {\rm domC}(\varphi)$
$ = \{0,1\}^m$).   
Next, assume by induction that for every $x \in \{0,1\}^m:$ 
 \ $\alpha'_{i-1} \circ \ldots \circ \alpha'_1(x) = u$ \ is a prefix of
 \ $\beta_{i-1} \circ \ldots \circ \beta_1(x) = u \, z$. 
Then \ $\beta_i (u \, z) = v \, z \, 0^{\ell_{i+1}-\ell_i -|v| + |u|}$ \ (if
$\ell_i \leq \ell_{i+1}$); or \ $\beta_i (u \, z) = v \, z_1$ \ (if
$\ell_i \geq \ell_{i+1}$, with $|z_1| = \ell_{i+1} - |v|$ and $z = z_1z_2$). 
In either case we find that 
 \ $\alpha'_i (\alpha'_{i-1} \circ \ldots \circ \alpha'_1(x)) = v$ \ is a 
prefix of \ $\beta_i (\beta_{i-1} \circ \ldots \circ \beta_1(x))$
$ = \beta_i (u \, z)$.

Hence, when $i = N$ we obtain for any $x \in \{0,1\}^m$:
 \ $\beta_N \circ \ldots \circ \beta_1(x) = y \, s$ \ is a prefix of
 \ $\alpha'_N \circ \ldots \circ \alpha'_1(x) = \varphi(x) = y$ \ for some 
$y$ and $s$ with $|y \, s| = \ell_N = n$. Since 
$y \in {\rm imC}(\varphi) \subseteq \{0,1\}^n$ we conclude that $s$ is
empty, hence \ $\beta_N \circ \ldots \circ \beta_1(x) = $
 $\alpha'_N \circ \ldots \circ \alpha'_1(x)$. 
 \ \ \ {\sf [End, proof of Claim.]}
 
\medskip

At this point we have expressed $\varphi$ as a product of $N$ elements 
$\beta_i \in {\it lep}M_{k,1}$, where $N = |\varphi|_{M_{k,1}}$.
We now want to find the word-length of each $\beta_i$ over
$\Gamma_{{\it lep}M_{k,1}} \cup \tau$, in order to find an upper bound on
the total word-length of $\varphi$ over $\Gamma_{{\it lep}M_{k,1}} \cup \tau$. 
As we saw above, \ $\ell_i \leq 2 \, N$ \ for every $i = 1, \ldots, N$.
 
\smallskip

We examine each generator in $\Gamma_{M_{k,1}} \cup \tau$.

\smallskip
 
If $\alpha_i \in \tau$ then $\beta_i \in \tau$, so in this case  
 \ $|\beta_i|_{{\it lep}M_{k,1}} = 1$. 
 
\medskip

Suppose now that $\alpha_i \in \Gamma_{M_{k,1}}$. 
By Proposition \ref{circ_lepM} it is sufficient to construct a circuit that
computes $\beta_i$; the circuit can then be immediately translated into a
word over $\Gamma_{{\it lep}M_{k,1}} \cup \tau$ with linear increase in 
length. 

Since ${\rm domC}(\alpha_i) \subseteq \{0,1\}^{\leq 2}$, we can restrict 
$\alpha_i$ so that its domain code becomes a subset of $\{0,1\}^2$;
next, we extend $\alpha_i$ to a map $\alpha''_i$ that acts as the identity 
map on $\{0,1\}^2$ where $\alpha_i$ was undefined. The image code of
$\alpha''_i$ is a subset of $\{0,1\}^{\leq 3}$. In order to compute $\beta_i$ 
we first introduce a circuit $C(\alpha''_i)$ that computes $\alpha''_i$. 
A difficulty here is that $\alpha''_i$ does not produce fixed-length 
outputs in general, whereas $C(\alpha''_i)$ has to work with fixed-length 
inputs and outputs; so the output of $C(\alpha''_i)$ represents the output 
of $\alpha''_i$ indirectly, as follows: 

The circuit $C(\alpha''_i)$ 
has two input bits $u = u_1u_2 \in \{0,1\}^2$, and 5 output bits:
First there are 3 output bits $0^{3 - |v|} \, v \in \{0,1\}^3$, where 
$v = \alpha''_i(u)$; second, there are two more output bits, 
$c_1c_2 \in \{0,1\}^2$, defined by \ $c_1c_2 = {\sf bin}(3 - |v|)$ \ (the 
binary representation of the non-negative integer $3 - |v|$). Hence,
 \ $c_1c_2 = 00$ if $|v| = 3$, \ $c_1c_2 = 01$ if $|v| = 2$, \ $c_1c_2 = 10$ 
if $|v| = 1$; since $|v| > 0$, the value $c_1c_2 = 11$ will not occur.
Thus $c_1c_2 \, 0^{3 - |v|} \, v$ contains the same information as $v$,
but has the advantage of having a fixed length (always 5).
The circuit $C(\alpha''_i)$ can be built with a small constant number of 
{\sf and, or, not, fork} gates, and we will not need to know the details.

\smallskip

We now build a circuit for $\beta_i$.

\smallskip

\noindent $\bullet$ Circuit for $\beta_i$ if  $\ell_i \leq \ell_{i+1}$: \\  
On input $u \, z \in \{0,1\}^{\ell_i}$ (with $u \in \{0,1\}^2$), we want to
produce the output \ $v \ z \ 0^{\ell_{i+1}-\ell_i -|v| + |u|}$, where
$v = \alpha''_i(u)$. 

We first apply the circuit $C(\alpha''_i)$, thus obtaining
 \  $c_1c_2 \, 0^{3-|v|} \, v \, z$. 
Then we apply two {\sf fork} operations (always to the last bit in $z$) 
to produce \ $c_1c_2 \, 0^{3-|v|} \, v \, z \ b \, b$, where $b$ is the last 
bit of $z$. Applying a negation to the first $b$ and an {\sf and} operation, 
we obtain \ $c_1c_2 \, 0^{3-|v|} \, v \, z \ 0$. 
Applying \ $\ell_{i+1} - \ell_i - 1$ \ more {\sf fork} operations to the last
0 yields \ $c_1c_2 \, 0^{3-|v|} \, v \, z \ 0^{\ell_{i+1} - \ell_i - 1}$.

Next, we want to move $0^{3-|v|}$ to the right of the output, in order to
obtain \ $c_1c_2 \ v \, z \ 0^{3-|v| +\ell_{i+1} - \ell_i - 1}$.
For this effect we introduce a {\it controlled cycle}. Let 
 \ $\kappa: x_1x_2x_3 \in \{0,1\}^3 \longmapsto x_3x_1x_2$ \ be the 
usual cyclic permutations of 3 bit positions.
The controlled cycle acts as the identity map when $c_1c_2 = 00$ or $11$,
$\tau_{1,2}$ when $c_1c_2 = 01$, and $\kappa$ when $c_1c_2 = 10$. 
More precisely,

 \ \ \ \ \ $\kappa_c: \ c_1c_2 \, x_1x_2x_3 \in \{0,1\}^5 \ \ \longmapsto$
$ \ \ \left\{ \begin{array}{ll}
         c_1c_2 \, x_1x_2x_3 \ &  \mbox{if \ $c_1c_2 = 00$ \ or \ $11$,} \\   
         c_1c_2 \, x_2x_1x_3 \ &  \mbox{if \ $c_1c_2 = 01$,} \\   
         c_1c_2 \, x_3x_1x_2 \ &  \mbox{if \ $c_1c_2 = 10$}.
         \end{array}   \right. $

\noindent We apply $\ell_i$ copies of $\kappa_c(c_1,c_2,.,.,.)$ (all 
controlled by the same value of $c_1c_2$) to 
 \ $0^{3-|v|} \, v \, z$. The first $\kappa_c(c_1,c_2,.,.,.)$ is applied to
the 3 bits $0^{3-|v|} \, v$, producing 3 bits $y_1y_2y_3$; the second 
$\kappa_c(c_1,c_2,.,.,.)$ is applied to $y_2y_3$ and the first bit of $z$, 
producing 3 bits $y'_1y'_2y'_3$; the third $\kappa_c(c_1,c_2,.,.,.)$ is 
applied to $y'_2y'_3$ and the second bit of $z$,  etc.
So, each one of the $\ell_i$ copies of $\kappa_c$ acts one bit further 
down than the previous copy of $\kappa_c$.
This will yield \ $c_1c_2 \, v \, z \ 0^{3-|v| +\ell_{i+1} - \ell_i - 1}$.
Finally, to make $c_1c_2$ disappear, we apply two {\sf fork} operations to 
$c_1$, then a negation and an {\sf and}, to make a 0 appear. We combine this 0
with $c_1$ and $c_2$ by {\sf and} gates, thus transforming $0c_1c_2$  into 0.
Finally, an {\sf or} operation between this 0 and the first bit of $v$ makes
this 0 disappear. 

The number of gates used to compute $\beta_i$ is $O(\ell_{i+1} + \ell_i)$,
which is $\leq O(N)$. 

\smallskip

\noindent $\bullet$ Circuit for $\beta_i$ if $\ell_i > \ell_{i+1}$: \\
On input $u \, z \in \{0,1\}^{\ell_i}$ (with $u \in \{0,1\}^2$), we want to
produce the output $v \, z_1$, where $v = \alpha''_i(u)$.  

We first apply the circuit $C(\alpha''_i)$, which yields the output 
 \ $c_1c_2 \, 0^{3-|v|} \, v \, z$.
Now we want to erase the $\ell_i - \ell_{i+1} +1$ last bits of $z$. 
For this we apply two {\sf fork} operations to the last bit of $z$ 
(let's call it $b$),
then a negation and an {\sf and}, to make a 0 appear. We combine this 0 with
the last $\ell_i - \ell_{i+1}$ bits of $z$, using that many {\sf and} gates,
turning all these bits into a single 0; finally, an {\sf or} operation between
this 0 and the bit of the remainder of $z$ makes this 0 disappear.
At this point, the output is \ $c_1c_2 \, 0^{3-|v|} \, v \, Z_1$, where 
$Z_1$ is the prefix of length $\ell_{i+1} -1$ of $z$. 
 
Next, we apply $O(\ell_{i+1})$ position transpositions to $Z_1$ in order move 
the two last bits of $Z_1$ to the front of $Z_1$. Let $b_1b_2$ be the last two
bits of $Z_1$; so, $Z_1 = z_0b_1b_2$ (where $z_0$ is the prefix of length 
$\ell_{i+1} -3$ of $z$); at this point, the output of the circuit is 
 \ $c_1c_2 \, 0^{3-|v|} \, v \, b_1b_2 \, z_0$. 

We now introduce a fixed small circuit with 7 input bits and 5 output bits,
defined by the following input-output map:

\smallskip

 \ \ \ \ \ 
$\omega_c: \ c_1c_2 \, x_1x_2x_3 \, b_1b_2 \in \{0,1\}^7 \ \ \longmapsto$
$ \ \ \left\{ \begin{array}{ll}
   c_1c_2 \, x_1x_2x_3 \ &  \mbox{if \ $c_1c_2 = 00$ \ or \ $11$,} \\
   c_1c_2 \, x_1x_2 \, b_1 \ &  \mbox{if \ $c_1c_2 = 01$,} \\
   c_1c_2 \, x_3 \, b_1b_2 \ &  \mbox{if \ $c_1c_2 = 10$}.
   \end{array}   \right. $

\smallskip

\noindent When this map is applied to $c_1c_2 \, 0^{3-|v|} \, v \, b_1b_2$ 
the output is therefore given by 

\smallskip

 \ \ \ \ \ 
$\omega_c: \ c_1c_2 \, 0^{3-|v|} \, v \, b_1b_2 \in \{0,1\}^7 \ \ \longmapsto$
$ \ \ \left\{ \begin{array}{ll}
   c_1c_2 \, v       \ &  \mbox{if \ $|v| = 3$,} \\
   c_1c_2 \, v \, b_1 \ &  \mbox{if \ $|v| = 2$,} \\
   c_1c_2 \, v \, b_1b_2 \ &  \mbox{if \ $|v| = 1$}.
   \end{array}   \right. $

\smallskip

\noindent A circuit for $\omega_c$ can be built with a small fixed number of
{\sf and, or, not, fork} gates, and we will not need to know the details.

\smallskip

After applying $\omega_c$ to 
 \ $c_1c_2 \, 0^{3-|v|} \, v \, b_1b_2 \, z_0$ \ the output has length 
$\ell_{i+1} + 2$; the ``$+2$'' comes from $c_1c_2$. The output is 
 \ $c_1c_2 \, v \, z_0$, or \ $c_1c_2 \, v \, b_1 \, z_0$, or \   
$c_1c_2 \, v \, b_1b_2 \, z_0$, depending on whether $|v| = 3, 2$, or $1$. 
 
We need to move $b_1b_2$ or $b_1$ (or nothing) back to the 
right-most positions of $z_0$. We do this by applying $\ell_{i+1}$   
copies of the {\it controlled cycle} $\kappa_c(c_1,c_2,.,.,.)$ (all copies
controlled by the same value of $c_1c_2$). We proceed in the same way as 
when we used $\kappa_c$ in the previous case, and we obtain the output 
 \ $c_1c_2 \, v \, z_0$ \ (if $|v| = 3$), or 
 \ $c_1c_2 \, v \, z_0 \, b_1$ \ (if $|v| = 2$), or  
 \ $c_1c_2 \, v \, z_0 \, b_1b_2$ \ (if $|v| = 1$).  

Finally, we erase $c_1c_2$ in the same way as in the previous case, thus 
obtaining the final output. 
The number of gates used to compute $\beta_i$ is $O(\ell_{i+1} + \ell_i)$
$\leq O(N)$.

\smallskip

This completes the constuction of a circuit for $\beta_i$.
Through this circuit, $\beta_i: \{0,1\}^{\ell_i} \to \{0,1\}^{\ell_{i+1}}$ 
is expressed as a word over the generating set
$\Gamma_{{\it lep}M_{k,1}} \cup \tau$, of length 
$\leq O(\ell_{i+1} + \ell_i)$ $\leq O(N)$.

Since we have described $\varphi$ as a product of $N = |\varphi|_{M_{k,1}}$ 
elements $\beta_i \in {\it lep}M_{k,1}$, each of word-length $O(N)$, 
we conclude that $\varphi$ has word-length $\leq O(N^2)$ over the generating 
set  $\Gamma_{{\it lep}M_{k,1}} \cup \tau$ of ${\it lep}M_{k,1}$.
 \ \ \ $\Box$

\bigskip

\noindent {\bf Question:} Does the distortion of ${\it lep}M_{k,1}$ in 
$M_{k,1}$ (over the generators of Theorem \ref{distor_lepM_M}) have an upper 
bound that is less than quadratic?

%%%%%%%%%%%%%%%%%%%%%%%%%%%%%%%%%%%%%%%%%%%%%%%%%%%%%%%%
% Section  3
%%%%%%%%%%%%%%%%%%%%%%%%%%%%%%%%%%%%%%%%%%%%%%%%%%%%%%%%

\section{Wordlength asymmetry vs.\ computational asymmetry}

\begin{pro} \label{wl_asym_vs_compasym} \
The word-length asymmetry function $\lambda$ of the Thompson group
${\it lp}G_{2,1}$ within the Thompson monoid ${\it lep}M_{2,1}$ is 
linearly equivalent to the computational asymmetry function $\alpha$: 

\smallskip

 \hspace{.5in}   $\alpha \ \simeq_{{\rm lin}} \ \lambda$.

\smallskip

\noindent Here the generating set used for ${\it lep}M_{2,1}$ is
 \ $\Gamma_{{\it lep}M_{2,1}} \, \cup \, \{\tau_{i,j} : 0\leq i < j\}$,  
where $\Gamma_{{\it lep}M_{2,1}}$ is finite. The gates used for circuits are 
any finite universal set of gates, together with the wire-swapping operations 
$\{\tau_{i,j} : 0\leq i < j\}$.

We can choose $\Gamma_{{\it lep}M_{2,1}}$ to consist  exactly of
the gates used in the circuits; then $\alpha = \lambda$.
\end{pro}
{\bf Proof.} For any $g \in {\it lp}G_{2,1}$ we have

\smallskip

$C(g^{-1}) \ \leq \ c_0 \cdot |g^{-1}|_{{\it lep}M_{2,1}} \ \leq \ $
$c_0 \cdot \lambda(|g|_{{\it lep}M_{2,1}}) \ \leq \ $ 
$c_0 \cdot \lambda(c_1 \cdot C(g))$.

\smallskip

\noindent The first and last ``$\leq$'' come from Prop.\ \ref{circ_lepM}
(since ${\it lp}G_{2,1} \subset {\it lep}M_{2,1}$),
and the middle ``$\leq$'' comes from the definition of $\lambda$; 
$c_0$ and $c_1$ are positive constants.  Hence,  

\smallskip

$\alpha(n) \ \leq \ c_0 \cdot \lambda(c_1 \, n)$ \ for all $n$.

\smallskip

\noindent In a very similar way we prove that 
 \ $\lambda(n) \leq \ c'_0 \cdot \alpha(c'_1 \, n)$ \ for some positive
constants $c'_0, c'_1$. 
 \ \ \ $\Box$

\begin{pro} \label{wl_asymM_vs_wl_asymLepM} \
The word-length asymmetry function $\lambda_{M_{2,1}}$ of the Thompson group
${\it lp}G_{2,1}$ within the Thompson monoid $M_{2,1}$ is polynomially 
equivalent to the word-length asymmetry function $\lambda_{{\it lep}M_{2,1}}$
of ${\it lp}G_{2,1}$ within the Thompson monoid ${\it lep}M_{2,1}$. More
precisely we have for all $n:$

\smallskip

$\lambda_{M_{2,1}}(n) \ \leq \ $
  $c_0 \cdot \lambda_{{\it lep}M_{2,1}}(c_1 \, n^2)$, 
 
\smallskip

$\lambda_{{\it lep}M_{2,1}}(n) \ \leq \ $
              $c'_0 \cdot \big(\lambda_{M_{2,1}}(c'_1 \, n)\big)^2$, 

\smallskip

\noindent where $c_0, c_1, c'_0, c'_1$ are positive constants.
Here the generating set used for ${\it lep}M_{2,1}$ is
 \ $\Gamma_{{\it lep}M_{2,1}} \, \cup \, \{\tau_{i,j} : 0\leq i < j\}$,
where $\Gamma_{{\it lep}M_{2,1}}$ is finite. The generating set used for
$M_{2,1}$ is
 \ $\Gamma_{M_{2,1}} \, \cup \, \{\tau_{i,j} : 0\leq i < j\}$,
where $\Gamma_{M_{2,1}}$ is a finite generating set of $M_{2,1}$.
\end{pro}
{\bf Proof.} For any $g \in {\it lp}G_{2,1}$ we have

\smallskip

$|g^{-1}|_{M_{2,1}} \ \leq \ c_0 \cdot |g^{-1}|_{{\it lep}M_{2,1}} \ \leq$
 \ $c_0 \cdot \lambda_{{\it lep}M_{2,1}}(|g|_{{\it lep}M_{2,1}}) \ \leq$
$ \ c_0 \cdot \lambda_{{\it lep}M_{2,1}}(c_1 \cdot |g|_{M_{2,1}}^2)$. 

\smallskip

\noindent The first ``$\leq$'' holds because 
${\it lp}G_{2,1} \subset {\it lep}M_{2,1} \subset M_{2,1}$ and because of the 
choice of the generating sets.
The second ``$\leq$'' holds by the definition of $\lambda_{{\it lep}M_{2,1}}$.
The third ``$\leq$'' comes from the quadratic distortion of ${\it lep}M_{2,1}$
in $M_{2,1}$ (Theorem \ref{distor_lepM_M}). 
For the same reasons we also have the following:

\smallskip

$|g^{-1}|_{{\it lep}M_{2,1}} \ \leq \ c'_0 \cdot |g^{-1}|_{M_{2,1}}^2 $
$ \ \leq \ c'_0 \cdot (\lambda_{M_{2,1}}(|g|_{M_{2,1}}))^2 \ \leq \ $ 
$ c'_0 \cdot (\lambda_{M_{2,1}}(c_1 \cdot |g|_{{\it lep}M_{2,1}}))^2$ 

\noindent where $c'_0, c'_1$ are positive constants.
 \ \ \ $\Box$

%%%%%%%%%%%%%%%%%%%%%%%%%%%%%%%%%%%%%%%%%%%%%%%%%%%%%%%%
% Section  4
%%%%%%%%%%%%%%%%%%%%%%%%%%%%%%%%%%%%%%%%%%%%%%%%%%%%%%%%

\section{Reversible representation over the Thompson groups}

Theorems \ref{repres_M_by_G} and \ref{repres_G_circuits_by_G} below 
introduce a representation of elements of the Thompson monoid
${\it lep}M_{2,1}$ by elements of the Thompson group $G_{2,1}$, in
analogy with the Toffoli representation (Theorem \ref{ToffRepr} above), 
and the Fredkin representation (Theorem \ref{FredToffRepr} above). 
Our representation preserves complexity, up to a polynomial change, and 
uses only {\em one} constant-0 input. Note that although the 
functions and circuits considered here use fixed-length inputs and outputs, 
the representations is over the Thompson group $G_{2,1}$, which includes 
functions with variable-length inputs and outputs.

In the Theorem below, $\Gamma_{G_{2,1}}$ is any finite generating set of
$G_{2,1}$. We denote the length of a word $w$ by $|w|$, and we denote the
size of a circuit $C$ by $|C|$. The gates {\sf and, or, not} will also be
denoted respectively by $\wedge, \vee, \neg$.
We distinguish between a word $W_f$ (over a generating set of
$G_{2,1}$) and the element $w_f$ of $G_{2,1}$ represented by $W_f$.

%%%
\begin{thm} {\bf (Representation of boolean functions by the Thompson 
group).} \label{repres_M_by_G} \
Let $f : \{0,1\}^m \to \{0,1\}^n$ be any total function and
let $C_f$ be a minimum-size circuit (made of $\wedge, \vee, \neg$,
{\sf fork}-gates and wire-swappings $\tau_{i,j}$) that computes $f$.
Then there exists a word $W_f$ over the generating set
 \ $\Gamma_{G_{2,1}} \cup \{\tau_{i,i+1} : 1 \leq i \}$ \ of $G_{2,1}$
such that:

\smallskip

\noindent $\bullet$ \ For all $x \in \{0,1\}^m$:
                    \ $w_f(0 \, x) = 0 \ f(x) \ x$, \ \ where $w_f$ is the
element of $G_{2,1}$ represented by $W_f$.

\smallskip

\noindent $\bullet$ \ The length of the word $W_f$ is bounded by 
 \ \ $|W_f| \leq O(|C_f|^4)$.

\smallskip

\noindent $\bullet$ \  The largest subscript of any transposition
$\tau_{i,i+1}$ occurring in $W_f$ has an upper bound \
$\leq |C_f|^2 + 2$.
\end{thm}
{\bf Proof.} Wire-swappings in circuits are represented by the position
transpositions $\tau_{i,i+1} \in G_{2,1}$.
The gates {\sf not}, {\sf or}, and {\sf and} of
circuits are represented by the following elements of $G_{2,1}$:

\medskip

\noindent
$ \varphi_{\neg} = \left[ \! \! \begin{array}{cc}
0 & 1 \\
1 & 0
\end{array} \! \! \right]
$, \ \  
$ \varphi_{\vee} = \left[ \! \! \begin{array}{cc}
0x_1x_2             & 1x_1x_2  \\
(x_1\vee x_2)x_1x_2 & (\, {\overline{x_1\vee x_2}} \, )x_1x_2
\end{array}  \!\!  \right]
$, \ \  
$ \varphi_{\wedge} = \left[ \! \! \begin{array}{cc}
0x_1x_2               & 1x_1x_2  \\
(x_1\wedge x_2)x_1x_2 & (\, {\overline{x_1\wedge x_2}} \, )x_1x_2
\end{array}  \! \! \right]
$,

\medskip

\noindent
where $x_1, x_2$ range over $\{0,1\}$. Hence the domain and image codes of
$\varphi_{\vee}$ and $\varphi_{\wedge}$ are all equal to \
$\{0,1\}^3$.

To represent {\sf fork} we use the following element, in which we recognize
$\sigma \in F_{2,1}$, one of the commonly used generators of the 
Thompson group $F_{2,1}$:

\smallskip

$ \sigma \ = \ \left[ \! \! \begin{array}{ccc}
0  & 10 & 11   \\
00 & 01 & 1
\end{array}   \! \!    \right] \ \ = \ \ $  
$   \left[ \! \! \begin{array}{rrcc}
    00  & 01  & 10 & 11   \\   
    000 & 001 & 01 & 1
\end{array}   \! \!    \right] $.

\smallskip

\noindent Note that $\sigma$ agrees with {\sf fork} only on input 0, but 
that is all we will need. By its very essense, the forking operation 
cannot be represented by a length-equality preserving element of $G_{2,1}$,
because $G_{2,1} \cap {\it lep}M_{2,1} = {\it lp}G_{2,1}$ (the group of 
length-preserving elements of $G_{2,1}$). 
A small remark: In \cite{BiThomps, BiCoNP, BiFact}, what we call 
``$\sigma$'' here, was called ``$\sigma^{-1}$''.

\medskip

\noindent We will occasionally use the wire-swapping $\tau_{i,j}$ 
($1 \leq i < j$); note that $\tau_{i,j}$ can be expressed in terms of 
transpositions of neighboring wires as follows: 

\smallskip

$\tau_{i,j}(.) \ = \ $
 $ \tau_{i,i+1} \ \tau_{i+1,i+2} \ \ldots \ \tau_{j-2,j-1} \ \tau_{j-1,j}$
 $ \ \tau_{j-2,j-1} \ \ldots \ \tau_{i+1,i+2} \ \tau_{i,i+1}(.)$

\smallskip

\noindent so the word-length of $\tau_{i,j}$ over
 $\{\tau_{\ell,\ell+1} : 1 \leq \ell\}$ is \ $\leq 2(j-i) -1$.

\bigskip

For \ $x = x_1 \ldots x_m \in \{0,1\}^m$ \ and \
$f(x) = y = y_1 \ldots y_n \in \{0,1\}^n$, we will construct a word $W_f$
over the generators $\Gamma_{G_{2,1}} \cup \{\tau_{i,i+1}: 1 \leq i\}$ \ of
$G_{2,1}$, such that $W_f$ defines the map
 \ $w_f(.): 0 \, x \mapsto 0 \ f(x) \ x$.

The circuit $C_f$ is partitioned into {\it slices} $c_{\ell}$ ($\ell = 1,
\ldots, L$). Two gates $g_1$ and $g_2$ are in the same slice iff the length 
of the {\it longest} path from $g_1$ to any input port is the same as the 
length of the longest path from $g_2$ to any input port.  
We assume that $C_f$ is {\it strictly layered}, i.e., each gate in slice 
$c_{\ell}$ only has in-wires coming from slice $c_{\ell-1}$, and out-wires 
going toward slice $c_{\ell+1}$, for all $\ell$. 
To make a circuit $C$ strictly layered we need to add at most $|C|^2$ 
identity gates (see p.\ 52 in \cite{BiCoNP}).
The input-output map of slice $c_{\ell}$ has the form

\smallskip

$c_{\ell}(.): \ y^{(\ell-1)} \ = \ $
$y_1^{(\ell-1)} \ldots y_{n_{\ell-1}}^{(\ell-1)} \ \in \ \{0,1\}^{n_{\ell-1}}$
$ \ \longmapsto \ y^{(\ell)} \ = \ $
$y_1^{(\ell)} \ldots y_{n_{\ell}}^{(\ell)} \ \in \ \{0,1\}^{n_{\ell}}$.

\smallskip

\noindent Then $y^{(0)} = x$ and $y^{(L)} = y$, where $x \in \{0,1\}^m$ is
the input and $y \in \{0,1\}^n$ is the output of $C_f$.
Each slice is a circuit of depth 1.

\bigskip

Before studying in more detail how $C_f$ is built from slices, let us see how
a slice is built from gates (inductively, one gate at a time).

Let $C$ be a depth-1 circuit with $k+1$ gates, obtained by adding one gate
to  a depth-1 circuit $K$ with $k$ gates.
Let \ $K(.): x_1 \ldots x_m \longmapsto y_1 \ldots y_n$ \ be the
input-output map of the circuit $K$. Assume by induction that $K$ is
represented by a word $W_K$ over the generating set \
$\Gamma_{G_{2,1}} \cup \{\tau_{i,i+1} : 1 \leq i\}$ \ of $G_{2,1}$.
The input-output map of $W_K$ is, by induction hypothesis, \

\smallskip

$w_K(.): 0 \, x_1 \ldots x_m \ \longmapsto \ $
$ 0 \, y_1 \ldots y_n \, x_1 \ldots x_m$.

\smallskip

\noindent
The word $W_C$ that represents $C$ over $G_{2,1}$ is obtained as follows
from $W_K$; there are several cases, depending on the gate that is added to
$K$ to obtain $C$.

\bigskip

\noindent {\sc Case 1}: An identity-gate (or a {\sf not}-gate) is added to 
$K$ to form $C$, i.e.,

\smallskip

$C(.): x_1 \ldots x_m x_{m+1} \ \longmapsto \ y_1 \ldots y_n x_{m+1}$

\smallskip

(or, \ $C(.):$
$x_1 \ldots x_m x_{m+1} \ \longmapsto \ y_1 \ldots y_n {\overline x_{m+1}}$).

\smallskip

\noindent Then $W_C$ is given by

\smallskip

$w_C: \ \ 0\, x_1 x_2 \ldots x_m x_{m+1} \ \stackrel{\sigma}{\longmapsto} \ $
$00 \, x_1 x_2 \ldots x_m x_{m+1} \ \stackrel{\tau_{3,m+3}}{\longmapsto} \ $
$00 \, x_{m+1}\, x_2 \ldots x_m x_1 \ $
  $\stackrel{\varphi_{\vee}}{\longmapsto} \ $

\smallskip

$x_{m+1}0 \, x_{m+1}\, x_2 \ldots x_m x_1 \ $
 $ \stackrel{\tau_{3,m+3}}{\longmapsto} \ $
$x_{m+1} \, 0 \, x_1 x_2 \ldots x_m x_{m+1} \ \stackrel{\pi}{\longmapsto} $
$0 \, x_1 x_2 \ldots x_m x_{m+1} x_{m+1} \ \stackrel{w_K}{\longmapsto} $

\smallskip

$0 \, y_1 \ldots y_n \, x_1 \ldots x_m x_{m+1} x_{m+1} \ $
$\stackrel{\pi'}{\longmapsto} $
$0 \, y_1 \ldots y_n \, x_{m+1} \, x_1 \ldots x_m x_{m+1}$ ,

\smallskip

\noindent where
 \ $\pi(.) = \ \tau_{m+1,m+2} \ \ldots \ \tau_{2,3} \ \tau_{1,2}(.)$
 \ shifts $x_{m+1}$ from position 1 to position $m+2$, while shifting
$0x_1 \ldots x_m$ one position to the left; and \  $\pi'(.) = \ $
$\tau_{m+2,m+3} \ \ldots \ \tau_{n+m+1,n+m+2} \ \tau_{n+m+2,n+m+3}(.)$
 \ shifts $x_{m+1}$ from position $n+m+3$ to position $n+2$, while shifting
$x_1 \ldots x_m$ one position to the right.

So, \ $W_C \ = \ $
$\pi' \ W_K \ \pi \ \tau_{3,m+3} \ \varphi_{\vee} \ \tau_{3,m+3} \ \sigma$,
noting that functions act on the left. Thus, $|W_C| = |W_K| +m+n+5$ if we
use all of $\{\tau_{i,j}: 1\leq i<j\}$ in the generating set; over 
$\{\tau_{i,i+1}: 1\leq i\}$,
$\tau_{3,m+3}$ has length $\leq 2m-1$, hence \ $|W_C| \leq 3m +n +4$.
If we denote the maximum index in the transpositions occurring in $W_C$ by
$J_C$ then we have \ $J_C = {\rm max}\{ J_K, \ n+m+3\}$.

\smallskip

In case a {\sf not}-gate is added (instead of an identity gate), 
$\varphi_{\vee}$ is replaced by $\varphi_{\neg} \, \varphi_{\vee}$ in $W_C$, 
and the result is similar.

\bigskip

\noindent {\sc Case 2}: An {\sf and}-gate (or an {\sf or}-gate) is added 
to $K$ to form $C$, i.e.,

\smallskip

$C(.): x_1 \ldots x_m x_{m+1} x_{m+2} \ \longmapsto \ $
   $y_1 \ldots y_n \ (x_{m+1} \wedge x_{m+2})$

\smallskip

(or, \ $C(.): x_1 \ldots x_m x_{m+1} x_{m+2} \ \longmapsto \ $
 $y_1 \ldots y_n \ (x_{m+1} \vee x_{m+2})$).

\smallskip

\noindent Then $W_C$ is given by

\smallskip

$w_C: \ \ 0 \, x_1 x_2 \ldots x_m x_{m+1} x_{m+2} \   $
$\stackrel{\sigma}{\longmapsto} \ $
$00 \, x_1 x_2 \ldots x_m x_{m+1} x_{m+2} \   $
$\stackrel{\tau_{2,m+3}}{\longmapsto} \ \ $
    $\stackrel{\tau_{3,m+4}}{\longmapsto} \ $

\smallskip

$0 \, x_{m+1} x_{m+2} \, x_2 \ldots x_m 0 \, x_1 \ $
$\stackrel{\varphi_{\wedge}}{\longmapsto} \ $
$(x_{m+1} \wedge x_{m+2}) \ x_{m+1} x_{m+2} \, x_2 \ldots x_m 0 \, x_1 \ $
$\stackrel{\tau_{2,m+3}}{\longmapsto} \ \ $
$\stackrel{\tau_{3,m+4}}{\longmapsto} \ $

\smallskip

$(x_{m+1} \wedge x_{m+2}) \ 0 \, x_1 x_2 \ldots x_m x_{m+1} x_{m+2} \ $
$\stackrel{\pi}{\longmapsto} \ $
$0 \, x_1 x_2 \ldots x_m \ (x_{m+1} \wedge x_{m+2}) \ x_{m+1} x_{m+2}$
$ \ \stackrel{w_K}{\longmapsto} \ $

\smallskip

$0 \, y_1 \ldots y_n \, x_1 x_2 \ldots x_m \ (x_{m+1} \wedge x_{m+2}) \ $
             $x_{m+1} x_{m+2} \ $
$\stackrel{\pi'}{\longmapsto} \ $

\smallskip

$0 \, y_1 \ldots y_n \ (x_{m+1} \wedge x_{m+2}) \ $
                                   $x_1 x_2 \ldots x_m x_{m+1} x_{m+2}$ ,

\smallskip

\noindent where \ $\pi = \ \tau_{m+1,m+2} \ \ldots \ \tau_{2,3} \ \tau_{1,2}$
 \ shifts $(x_{m+1} \wedge x_{m+2})$ from position 1 to position $m+2$,
while shifting \ $0 \, x_1 x_2 \ldots x_m$ \ one position to the left;
and \ $\pi' = \ \tau_{m+2,m+3} \ \ldots \ \tau_{m+n+1,m+n+2} \ $
shifts $(x_{m+1} \wedge x_{m+2})$ from position $n+m+2$ to position
$m+2$, while shifting \ $x_1 \ldots x_m$ \ one position to the right.

So, \ $W_C = \pi' \ W_K \ \pi \ \tau_{3,m+4} \ \tau_{2,m+3} \ $
$\varphi_{\wedge} \ \tau_{3,m+4} \ \tau_{2,m+3} \ \sigma$, hence \
$|W_C| = |W_K| + n +m + 7$ if all of $\{\tau_{i,j}: 1\leq i<j\}$ is used in
the generating set;
over $\{\tau_{i,i+1}: 1\leq i\}$, $\tau_{3,m+4}$ and $\tau_{2,m+3}$ have
length $\leq 2(m+1) -1$, so \ $|W_C| \leq |W_K| + 5m+ n +9$. Moreover,
 \ $J_C = {\rm max}\{ J_K, \ m+n+2\}$.

\bigskip

\noindent {\sc Case 3}: A {\sf fork}-gate is added to $K$ to
form $C$, i.e.,

\smallskip

$C(.): x_1 \ldots x_m x_{m+1} \ \longmapsto \ $
   $y_1 \ldots y_n \, x_{m+1} x_{m+1}$.

\smallskip

\noindent Then $W_C$ is given by

\smallskip

$w_C: \ \ 0 \, x_1 x_2 \ldots x_m x_{m+1} \   $
$\stackrel{\sigma^2}{\longmapsto} \ $
$000 \, x_1 x_2 \ldots x_m x_{m+1} \   $
$\stackrel{\tau_{3,m+4}}{\longmapsto} \ $
$00 \, x_{m+1} x_1 x_2 \ldots x_m 0 \ $
$\stackrel{\varphi_{\vee}}{\longmapsto} \ $

\smallskip

$x_{m+1} 0 \, x_{m+1} x_1 x_2 \ldots x_m 0 \ $
$\stackrel{\tau_{1,m+4}}{\longmapsto} \ $
$00 \, x_{m+1} x_1 x_2 \ldots x_m x_{m+1} \ $
$\stackrel{\varphi_{\vee}}{\longmapsto} \ $
$x_{m+1}0 \, x_{m+1} x_1 x_2 \ldots x_m x_{m+1} \ $
$\stackrel{\pi}{\longmapsto} \ $

\smallskip

$0 \, x_1 x_2 \ldots x_m x_{m+1} x_{m+1} x_{m+1} \ $
$\stackrel{w_K}{\longmapsto} \ $
$0 \, y_1 \ldots y_n \, x_1 x_2 \ldots x_m \, x_{m+1} x_{m+1} x_{m+1} \ $
$\stackrel{\pi'}{\longmapsto} \ $

\smallskip

$0 \, y_1 \ldots y_n \, x_{m+1} x_{m+1} \, x_1 x_2 \ldots x_m x_{m+1}$ ,

\smallskip

\noindent where \ $\pi = \ $
$\tau_{m+3,m+4} \ \ldots \ \tau_{1,2} \ \tau_{m+3,m+4} \ \ldots \ \tau_{3,4}$
 \ shifts the two copies of $x_{m+1}$ at the left end from positions 1 and 3
to positions $m+3$ and $m+4$, while shifting 0 to position 1 and shifting
 \ $x_1 \ldots x_m$ \ two positions to the left; and
 \ $\pi' = \ \tau_{m+3,m+4} \ \ldots \ \tau_{m+n+2,m+n+3} \ $
$\tau_{m+2,m+3} \ \ldots \ \tau_{m+n+1,m+n+2}$ \
shifts $x_{m+1} x_{m+1}$ from positions $m+n+2$ and $m+n+3$ to positions
$m+2$ and $m+3$, while shifting \ $x_1 \ldots x_m$ \ two positions to the
right.

So, \ $W_C = \pi' \ W_K \ \pi \ \varphi_{\vee} \ \tau_{1,m+4} \ $
$\varphi_{\vee} \ \tau_{3,m+4} \ \sigma^2$, hence \
$|W_C| = |W_K| + 2m+ n + 10$, if all of $\{\tau_{i,j}: 1\leq i<j\}$ is used
in the generating set;
over $\{\tau_{i,i+1}: 1\leq i\}$, $\tau_{1,m+4}$ has length $\leq 2(m+3) -1$
and $\tau_{3,m+4}$ has length $\leq 2m-1$. Hence, \
$|W_C| \leq |W_K| + 6m +n +14$.
Moreover, \ $J_C = {\rm max}\{ J_K, \ m+n+3\}$.

\bigskip

In all cases, \ $|W_C| \leq |W_K| + c \cdot (m+n+1)$ \ (for some constant
$c > 1$), and \ $J_C \leq {\rm max}\{J_K, n+m+3\}$.
Thus, each slice $c_{\ell}$, with input-output map \
$c_{\ell}(.) : y^{(\ell-1)} \ \longmapsto \ y^{(\ell)}$,
is represented by a word $W_{c_{\ell}}$ with map \
$w_{c_{\ell}}(.): 0 \, y^{(\ell-1)} \ \longmapsto \ $
  $0 \, y^{(\ell)} \, y^{(\ell-1)}$, such that \
$|W_{c_{\ell}}| \ \leq \ c \cdot ( n_{\ell-1}^2 + n_{\ell}^2)$ \ (for some
constant $c>1$), and
$J_{c_{\ell}} \ \leq \ n_{\ell-1} + n_{\ell} + c$.

\bigskip

Regarding wire-crossings, we do not include them into other slices; we
put the wire-crossings into pure wire-crossing slices. So we
consider two kinds of slices: Slices entirely made of wire-crossings and
identities, slices without any wire-crossings. Wire-crossings in
circuits are identical to the group elements $\tau_{i,i+1}$.

\bigskip

We now construct the word $W_f$ from the words $W_{c_{\ell}}$
($\ell = 1, \ldots, L$). First observe that since the map $w_{c_{\ell}}(.)$
is a right-ideal isomorphism (being an element of $G_{2,1}$), we not only
have

\smallskip

$w_{c_{\ell}}(.) : \ 0 \, y^{(\ell-1)} \ \longmapsto \ $
                  $0 \, y^{(\ell)} y^{(\ell-1)}$

\smallskip

\noindent but also

\smallskip

$w_{c_{\ell}}(.): \ 0 \, y^{(\ell-1)} y^{(\ell-2)} \ldots y^{(1)} y^{(0)}$
$ \ \longmapsto \ 0 \, y^{(\ell)} y^{(\ell-1)} y^{(\ell-2)} \ldots $
$y^{(1)} y^{(0)}$.

\smallskip

\noindent Then, by concatenating all $W_{c_{\ell}}$ (and by recalling that
$y = y^{(L)}$ and $x = y^{(0)}$) we obtain

\smallskip

$w_{c_L} \ w_{c_{L-1}} \ \ldots \ w_{c_2} \ w_{c_1}(.): \ \ 0 \, x \ $
$\longmapsto \ $
$0 \, y \, y^{(L-1)} \ \ldots \ y^{(2)} \, y^{(1)} \, x$.

\smallskip

\noindent Let $\pi_{C_f}$ be the position permutation that shifts $y$ right
to the positions just right of $x$:

\smallskip

$\pi_{C_f}: \ \ 0 \, y \, y^{(L-1)} \ \ldots \ y^{(2)} \, y^{(1)} \, x \ $
$\longmapsto \ $
$0 \, y^{(L-1)} \ \ldots \ y^{(2)} \, y^{(1)} \, x \, y$.

\smallskip

\noindent Observe that for \ $(W_{c_{L-1}} \ \ldots \ W_{c_2} \ W_{c_1})^{-1}$
 \ we have

\smallskip

$(w_{c_{L-1}} \ \ldots \ w_{c_2} \ w_{c_1})^{-1}(.): \ \ $
$0 \, y^{(L-1)} \ \ldots \ y^{(2)} \, y^{(1)} \, x \, y \ $
$\longmapsto \ $
$0 \, x \, y$.

\smallskip

\noindent Then we have:

\smallskip

$w_{c_L} \ w_{c_{L-1}} \ \ldots \ w_{c_2} \ w_{c_1}$ $ \ \pi_{C_f} \ $
$(w_{c_{L-1}} \ \ldots \ w_{c_2} \ w_{c_1})^{-1}(.): $
$ \ \ 0 \, x \ \longmapsto \ 0 \, x \, y$ .

\smallskip

\noindent By using the position permutation 
 \ $\pi_{m,n}: 0 \, x \, y \ \longmapsto \ 0 \, y \, x$,
we now see how to define $W_f$:

\[ W_f \ = \ \pi_{m,n} \ W_{c_L} \ W_{c_{L-1}} \ \ldots \ W_{c_2} \ W_{c_1}
 \ \pi_{C_f} \ (W_{c_{L-1}} \ \ldots \ W_{c_2} \ W_{c_1})^{-1}. \]

\noindent Then we have:

\smallskip

 \hspace{1in}   $w_f(.): \ 0 \, x \ \longmapsto \ 0 \, y \, x$,

\smallskip

\noindent where $y = f(x)$.

\bigskip

 Finally, we need to examine the length of the word $W_f$ in terms
of the size of the circuit $C_f$ that computes $f: \{0,1\}^m \to \{0,1\}^n$.

The position permutation $\pi_{m,n}$ shifts the $n = |y|$ letters of $y$
to the left over the $m = |x|$ positions of $x$. So, $\pi_{m,n}$ can be
written as the product of $nm$ transpositions in 
$\{\tau_{i,i+1}: 1 \leq i\}$, with maximum subscript 
$J_{\pi_{m,n}} \leq m+n+1$.

The position permutation $\pi_{C_f}$ shifts $y$ to the right from positions
in the interval $[2, n+1]$ within the string
 \ $0 \, y \, y^{(L-1)} \ \ldots \ y^{(2)} \, y^{(1)} \, x$ \ to positions
in the interval \ $[2 + \sum_{i=0}^{L-1} n_i , \ 2 + \sum_{i=0}^L n_i]$ 
 \ within the string
  \ $0 \, y^{(L-1)} \ \ldots \ y^{(2)} \, y^{(1)} \, x \, y$.
Note that \ $\sum_{i=0}^L n_i = |C_f|$ \ (the size of the circuit $C_f$), 
and $n_L = |y| = n$, \ $n_0 = |x| = m$. We shift $y$ starting with the 
right-most letters of $y$.  This takes 
 \ $n \ \sum_{i=0}^{L-1} n_i \, = \, n \, (|C_f| - n)$ \ transpositions
in $\{\tau_{i,i+1} : 1 \leq i\}$, with maximum subscript 
$J_{\pi_{C_f}} = |C_f| + 2$.

We saw already that \ $|W_{c_{\ell}}| \leq c \ (n_{\ell-1}^2 + n_{\ell}^2)$,
and \ $J_{c_{\ell}} \leq n_{\ell-1} + n_{\ell} +c$, for some constant $c>1$.
Note that \ $\sum_{i=0}^L n_i^2 \ \leq \ (\sum_{i=0}^L n_i)^2 \ = \ |C_f|^2$.
Hence we have:
 \ $|W_f| \ \leq \ c_o \ |C_f|^2$, \ for some constant $c_o > 1$.
Moreover, the largest subscript in any transposition occurring in
$W_f$ is  \ $J_{W_f} \ \leq \ |C_f| + 2$.

Recall that we assumed that our circuit $C_f$ was {\it strictly layered}, 
and that the circuit size has to be squared (at most) in order to make the 
circuit strictly layered. Thus, if $C_f$ was originally not strictly layered, 
our bounds become
 \ $|W_f| \ \leq \ c_o \ |C_f|^4$, \ and \ $J_{W_f} \ \leq \ |C_f|^2 +2$.
 \ \ \ \ \ $\Box$

%%%%%%%%

\bigskip

The next theorem gives a representation of a boolean permutation by an
element of the Thompson group $G_{2,1}$; the main point of the theorem is
the polynomial bound on the {\it word-length} in terms of {\it circuit size}.  

\begin{thm} {\bf (Representation of permutations by the Thompson group).} 
 \label{repres_G_circuits_by_G} \
Let $g: \{0,1\}^m \to \{0,1\}^m$ be any permutation and let $C_g$ and 
$C_{g^{-1}}$ be minimum-size circuits that compute $g$, respectively $g^{-1}$.
Then there exists a word $W_{(g,g^{-1})}$ over the generating set
 \ $\Gamma_{G_{2,1}} \cup \{\tau_{i,i+1} : 1 \leq i\}$ \ of $G_{2,1}$,
representing an element $w_{(g,g^{-1})} \in G_{2,1}$ such that:

\smallskip

\noindent $\bullet$ \ For all $x \in {\rm Dom}(g)$ and all 
$y \in {\rm Im}(g)$:
 $$w_{(g,g^{-1})}( 0 \, x) \ = \ 0 \ g(x), \ \ \
 {\it and} \ \ \ \ (w_{(g,g^{-1})})^{-1}( 0 \, y) \ = \ 0 \ g^{-1}(y),$$  
 
where $(w_{(g,g^{-1})})^{-1} \in G_{2,1}$ is represented by the
free-group inverse $(W_{(g,g^{-1})})^{-1}$ of the word $W_{(g,g^{-1})}$. 

\smallskip

\noindent $\bullet$ \ $w_{(g,g^{-1})}(.)$ \ and \ $(w_{(g,g^{-1})})^{-1}$ 
 \ stabilize both $0 \, \{0,1\}^*$ and $1 \, \{0,1\}^*$.

\smallskip

\noindent $\bullet$ \ We have a length upper bound 
 \ \ \ $|W_{(g,g^{-1})}| \ = \ $
 $|(W_{(g,g^{-1})})^{-1}| \ \ \leq \ O(|C_g|^4 + |C_{g^{-1}}|^4)$.

\smallskip

\noindent $\bullet$ \ The largest subscript of transpositions
$\tau_{i,i+1}$ occurring in $W_{(g,g^{-1})}$ is
$ \ \ \leq \ {\rm max}\{|C_g|^2, \  |C_{g^{-1}}|^2\} \ + 2$.
\end{thm}
Note that we distinguish between the word $W_{(g,g^{-1})}$ 
(over a generating set of $G_{2,1}$) and the element $w_{(g,g^{-1})}$ of 
$G_{2,1}$ represented by $W_{(g,g^{-1})}$. Also, note that although
$g$ is length-preserving ($g \in {\it lp}G_{2,1}$), 
$w_{(g,g^{-1})} \in G_{2,1}$ is not length-preserving.

\medskip

\noindent {\bf Proof.} Consider the position permutation 
 \ $\pi: 0 \, y \, x \ \longmapsto \ 0 \, x \, y$, for all
$x, y \in \{0,1\}^m$; \ we express $\pi$ as a composition of $\leq m^2$ 
position transpositions of the form $\tau_{i,i+1}$. 
Let $W_g$ be the word constructed in Theorem \ref{repres_M_by_G} for $g$, 
and let $W_{g^{-1}}$ be the word constructed for $g^{-1}$. We define
$W_{(g,g^{-1})}$ by 
$$W_{(g,g^{-1})} \ \ = \ \ (W_{g^{-1}})^{-1} \ \pi \ W_g . $$

\noindent Then for all $x \in {\rm Dom}(g)$ we have:  
 \ $w_{(g,g^{-1})}: \ 0 \, x \ \longmapsto \ 0 \, y$, \ where $y = g(x)$.
More precisely, for all $x \in {\rm domC}(g)$,  

\smallskip

$ 0 \, x \ \xrightarrow{w_g} \ 0 \ g(x) \ x \ = 0 \, y \, x \ $
$\xrightarrow{\pi} \ 0 \ x \, y \ = \ 0 \ g^{-1}(y) \ $ 
$\xrightarrow{(w_{g^{-1}})^{-1}} \ 0 \, y \ = \ 0 \ g(x)$.

\smallskip

\noindent Since ${\rm domC}(g)$ is a maximal prefix code, $w_{(g,g^{-1})}$ 
maps $0 \, \{0,1\}^*$ into $0 \, \{0,1\}^*$ (where defined). 

Similarly, for all $y \in {\rm Im}(g)$ $ = {\rm Dom}(g^{-1})$ we have:
 \ $(w_{(g,g^{-1})})^{-1}: \ 0 \, y \ \longmapsto \ 0 \, x$, 
 \ where $x = g^{-1}(y)$, $y = g(x)$.
Since ${\rm domC}(g^{-1})$ is a maximal prefix code,
$(w_{(g,g^{-1})})^{-1}$ maps $0 \, \{0,1\}^*$ into $0 \, \{0,1\}^*$
(where defined). 
Hence, elements of $0 \, \{0,1\}^*$ are never images of $1 \, \{0,1\}^*$.
Thus, $1 \, \{0,1\}^*$ is also stabilized by $w_{(g,g^{-1})}$ and by
$(w_{(g,g^{-1})})^{-1}$.

The length of the word $W_{(g,g^{-1})}$ is bounded as follows:
We have \ $|W_g| \, \leq \, c_o \ |C_g|^4$, and
 \ $|(W_{g^{-1}})^{-1}| \ = \ |W_{g^{-1}}| \ \leq \ c_o \ |C_{g^{-1}}|^4$,
by Theorem \ref{repres_M_by_G}. Moreover, $\pi$ can be expressed as the 
composition of \ $\leq m^2$ ($< |C_g|^2$) transpositions in
$\{\tau_{i,i+1} : 1 \leq i\}$.

The bound on the subscripts also follows from Theorem \ref{repres_M_by_G}.
 \ \ \ \ \ $\Box$

%%%%%%%%%%%%%%%%%%%%%%%%%%%%%%%%%%%%%%%%%%%%%%%%%%%%%%%%
% Section 5
%%%%%%%%%%%%%%%%%%%%%%%%%%%%%%%%%%%%%%%%%%%%%%%%%%%%%%%%

\section{Distortion vs.\ computational asymmetry}

We show in this Section that the computational asymmetry function $\alpha(.)$
is polynomially related to a certain distortion of the group 
${\it lp}G_{2,1}$. 

By Theorem \ref{repres_G_circuits_by_G}, for every element 
$g \in {\it lp}G_{2,1}$ there is an element $w_{(g,g^{-1})} \in G_{2,1}$ 
which agrees with $g$ on $0 \, \{0,1\}^*$, and which stabilizes 
$0 \, \{0,1\}^*$ and $1 \, \{0,1\}^*$.  
The main property of $W_{(g,g^{-1})}$ is that its length is polynomially 
bounded by the circuit sizes of $g$ and $g^{-1}$; that fact will be crucial 
later. First we want to study how $w_{(g,g^{-1})}$ is related to $g$. 
Recall that we distinguish between the word $W_{(g,g^{-1})}$
(over a generating set of $G_{2,1}$) and the element $w_{(g,g^{-1})}$ of
$G_{2,1}$ represented by $W_{(g,g^{-1})}$.

Theorem \ref{repres_G_circuits_by_G} inspires the following concepts.

\begin{defn} \label{stab_fix_0lowering} \    
Let $G$ be a subgroup of $G_{2,1}$.
For any prefix codes $P_1, \ldots, P_k \subset \{0,1\}^*$, the {\em joint 
stabilizer} (in $G$) of the right ideals $P_1\{0,1\}^*,$ $\ldots,$ 
$P_k \{0,1\}^*$ is defined by

\smallskip

${\sf Stab}_{G}(P_1, \ldots, P_k) \ = \  \big\{ g \in G : \ $
$g(P_i\{0,1\}^*) \subseteq P_i\{0,1\}^*$ \ for every 
$i = 1, \ldots, k \big\}$.

\smallskip

\noindent The {\em fixator} (in $G$) of $P_1\{0,1\}^*$ is defined by

\smallskip

${\sf Fix}_{G}(P_1) \ = \ \big\{ g \in G : \ $
   $g(x) = x$ \ for all $x \in P_1\{0,1\}^*) \big\}$.

\smallskip

\noindent The fixator is also called ``point-wise stabilizer''.
\end{defn}

The following is an easy consequence of the definition: 
${\sf Fix}_{G}(P_i)$ is a subgroup of $G$ ($\subseteq G_{2,1}$), for 
$i=1, \ldots, k$.
If the prefix codes $P_1, \ldots, P_k$ are such that the right ideals
$P_1\{0,1\}^*, \ \ldots, \ P_k \{0,1\}^*$ are two-by-two disjoint, and such 
that $P_1 \cup \ldots  \cup P_k$ is a maximal prefix code, then  
${\sf Stab}_{G}(P_1, \ldots, P_k)$ is closed under inverse. Hence in this 
case ${\sf Stab}_{G}(P_1, \ldots, P_k)$ is a subgroup of $G$. 

\smallskip

\noindent In particular, we will consider the following groups:  

\medskip

\noindent $\bullet$ \ The joint stabilizer of $0 \, \{0,1\}^*$ and 
$ 1 \, \{0,1\}^*$,

\smallskip

 \ \ \ \ \ ${\sf Stab}_{G}(0,1) \ = \ \big\{ g \in G : \ $
    $g(0 \, \{0,1\}^*) \subseteq 0 \, \{0,1\}^*$ \ {\it and} \   
    $g(1 \, \{0,1\}^*) \subseteq 1 \, \{0,1\}^* \big\}$.

\medskip

\noindent $\bullet$ \ The fixator of $0 \, \{0,1\}^*$,

\smallskip

 \ \ \ \ \ ${\sf Fix}_{G}(0) \ = \ \{ g \in G : \ g(x) = x$ 
 \ {\it for all} $x \in 0 \, \{0,1\}^* \}$.

\medskip

\noindent $\bullet$ \ The fixator of $ 1 \, \{0,1\}^*$,

\smallskip

 \ \ \ \ \ ${\sf Fix}_{G}(1) \ = \ \{ g \in G : \ g(x) = x$
 \ {\it for all} $x \in 1\, \{0,1\}^* \}$.

\medskip

\noindent Clearly, ${\sf Fix}_{G}(0)$ and ${\sf Fix}_{G}(1)$ are subgroups
of ${\sf Stab}_{G}(0,1)$. 

%% [Notation:  use 
%% $\partial_0, \partial$ (partial) instead of $\Lambda_0, \Lambda$ ?]

\begin{lem} {\bf (Self-embeddings of $G_{2,1}$).} \label{isom_of_G} \   
Let $G$ be a subgroup of $G_{2,1}$. Then $G$ is isomorphic to 
${\sf Fix}_{G}(1)$ and to ${\sf Fix}_{G}(0)$ by the following isomorphisms: 

\smallskip
 
$\Lambda_0: \ g \in G \ \longmapsto (g)_0 \in {\sf Fix}_{G}(1)$

\smallskip
 
$\Lambda_1: \ g \in G \ \longmapsto (g)_1 \in {\sf Fix}_{G}(0)$

\smallskip 

\noindent where $(g)_0$ and $(g)_1$ defined as follows for any 
$g \in G_{2,1}$:

\smallskip

$(g)_0 : \left\{ \begin{array}{lll}
   0 \, x \in 0 \, \{0,1\}^* & \longmapsto & 0 \ g(x) \\
   1 \, x \in 1 \, \{0,1\}^* & \longmapsto & 1 \, x 
   \end{array}   \right. $         \hspace{1in}
$(g)_1 : \left\{ \begin{array}{lll}
   1 \, x \in 1 \, \{0,1\}^* & \longmapsto & 1\ g(x) \\
   0 \, x \in 0 \, \{0,1\}^* & \longmapsto & 0 \, x 
   \end{array}   \right. $
\end{lem}
{\bf Proof.} It is straightforward to verify that $\Lambda_0$ and
$\Lambda_1$ are injective homomorphisms. That $\Lambda_0$ is onto 
${\sf Fix}_{G}(1)$ can be seen from the fact that every element of
${\sf Fix}_{G}(1)$ has a table of the form 

\medskip

$   \left[ \! \! \begin{array}{cccc}
    0x_1 & \ldots & 0x_n & 1   \\
    0y_1 & \ldots & 0y_n & 1
\end{array}   \! \!    \right] $

\smallskip

\noindent where $\{x_1, \ldots, x_n\}$ and $\{y_1, \ldots, y_n\}$ are
two maximal prefix codes, and 
 \ $\left[ \! \! \begin{array}{ccc}
 x_1 & \ldots & x_n \\
 y_1 & \ldots & y_n 
\end{array}   \! \!    \right] $
 \ is an arbitrary element of $G$.
 \ \ \ $\Box$

\begin{lem} \label{isom_of_GxG} \  
Let $G$ be a subgroup of $G_{2,1}$. Then the direct product $G \times G$ is 
isomorphic to ${\sf Stab}_{G}(0,1)$ by the isomorphism

\smallskip

 \ \ \ \ \  $\Lambda: \ (f,g) \in G \times G \ \ \longmapsto \ \ $
$\big(0 \, x \mapsto 0 \ f(x), \ 1 \, x \mapsto 1 \ g(x) \big) $
$ \ \in \ {\sf Stab}_{G}(0,1)$.
\end{lem}
{\bf Proof.} It is straightforward to verify that $\Lambda$ is a 
homomorphism. That $\Lambda$ is onto ${\sf Stab}_{G}(0,1)$ and injective 
follows from the fact that every element of ${\sf Stab}_{G}(0,1)$ has a table
of the form

\medskip 

 \ \ \ \ \ 
$   \left[ \! \! \begin{array}{ccc ccc}
    0x_1 & \ldots & 0x_m & 1x_1' & \ldots & 1x'_n   \\
    0y_1 & \ldots & 0y_m & 1y'_1 & \ldots & 1y'_n
\end{array}   \! \!    \right] $

\medskip

\noindent where $\{x_1, \ldots, x_m\}$, \ $\{y_1, \ldots, y_m\}$, 
 \ $\{x'_1, \ldots, x'_n\}$, and $\{y'_1, \ldots, y'_n\}$, are maximal prefix
codes, and 
 \ $\left[ \! \! \begin{array}{ccc}
 x_1 & \ldots & x_m \\
 y_1 & \ldots & y_m 
\end{array}   \! \!    \right] $ \ and 
 \ $\left[ \! \! \begin{array}{ccc}
 x'_1 & \ldots & x'_n \\
 y'_1 & \ldots & y'_n 
\end{array}   \! \!    \right] $
 \ are arbitrary elements of $G$ ($\subseteq G_{2,1}$).
 \ \ \ $\Box$

\bigskip

\noindent Lemmas \ref{isom_of_G} and \ref{isom_of_GxG} reveal certain
{\em self-similarity} properties of the Thompson group $G_{2,1}$. 
(Self-similarity of groups with total action on an infinite tree is an 
important subject, see \cite{Nekr_book}. However, the action of $G_{2,1}$ 
is partial, so much of the known theory does not apply directly.)  

The stabilizer and the fixators above have some interesting 
properties. 

\begin{lem} \label{properties_Fix_Stab} \!.  

\noindent 
{\rm (1)} \ \ For all $f, g \in G$: \ \ \ $(f)_0 \, (g)_1 = (g)_1 \, (f)_0$

 \ \ (i.e., the commutator of ${\sf Fix}_G(0)$ and ${\sf Fix}_G(1)$ is the 
 identity).

\smallskip 

\noindent
{\rm (2)} \ \ ${\sf Fix}_G(0) \cdot {\sf Fix}_G(1) \ = \ {\sf Stab}_G(0,1)$ 
 \ \ and \ \ ${\sf Fix}_G(0) \cap {\sf Fix}_G(1) = {\bf 1}$; 

\smallskip 

\noindent
{\rm (3)} \ \ ${\sf Stab}_G(0,1)$ is the internal direct product of 
${\sf Fix}_G(0)$ and ${\sf Fix}_G(1)$.

 \ \ (This is equivalent to the combination of {\rm (1)} and {\rm (2)}.)

\smallskip

\noindent
{\rm (4)} \ \ For all $f, g \in G$: 
 \ \ \ $\Lambda(f,g) = \Lambda_0(f) \cdot \Lambda_1(g)$, 
 \ \ \ $\Lambda_0(f) = \Lambda(f, {\bf 1})$, \ and 
 \ \ $\Lambda_1(g) = \Lambda({\bf 1}, g)$.
\smallskip

 \ \ Moreover, \ ${\sf Fix}_G(0) \ = \ \Lambda_1(G)$, \ \ 
 \ \ ${\sf Fix}_G(1) \ = \ \Lambda_0(G)$, \ \ and 
 \ \ ${\sf Stab}_G(0,1) \ = \ \Lambda(G \times G)$.

%\smallskip 
%
%\noindent
%{\rm (5)} \ \  For all $a_1 \ldots a_m \in \{0,1\}^*$ we have: 
% \ \ \ $\Lambda_{a_m} \circ \ldots \circ \Lambda_{a_1}(.) \ = \ $
%       $\Lambda_{a_1 \ldots a_m}(.)$, 
%
%\smallskip
%
% \ \ where $\Lambda_{a_1 \ldots a_m}$ is defined by \ \ \     
% $\Lambda_{a_1 \ldots a_m}: \ g \in G \ \longmapsto \ (g)_{a_1 \ldots a_m}$
%  $ \in {\sf Fix}_G\big( \, \overline{\{a_1 \ldots a_m\}} \, \big)$. \   
% 
%\smallskip
%
% \ \ Here $\overline{\{a_1 \ldots a_m\}}$ is a finite complementary prefix
% code of $\{a_1 \ldots a_m\}$, and $(g)_{a_1 \ldots a_m}$ is defined by
%
%\smallskip
% 
%\hspace{.5in} $(g)_{a_1 \ldots a_m} : \left\{ \begin{array}{lll}
%   a_1 \ldots a_m \ z & \longmapsto & a_1 \ldots a_m \ g(z),  \\
%   w & \longmapsto & w \ \ \ \ \ {\rm if} \ \ a_1 \ldots a_m \ \ {\rm is 
%        \ not \ a \ prefix \ of} \ w. 
%   \end{array}   \right. $ 
\end{lem}
{\bf Proof.} The proof is a straightforward verification. 
%For a prefix code $P \subset \{0,1\}^*$, a {\it complementary prefix code} 
%of $P$ is (by definition) any prefix code $P' \subset \{0,1\}^*$, such 
%that $P \cup P'$ is a maximal prefix code and such that 
%$P \, \{0,1\}^* \cap P' \, \{0,1\}^* = \varnothing$. 
%It is proved in \cite{BiCoNP} that every finite prefix code has a finite 
%complementary prefix code.
 \ \ \ $\Box$

\begin{lem} \label{tau_0} \   
For every position transposition $\tau_{i,j}$, with $1 \leq i < j$, we have

\smallskip

 \ \ \ \ \  $(\tau_{i,j})_0 \ = \ \tau_{2,i+1} \circ \tau_{3,j+1} \circ $
$(\tau_{1,2})_0 \circ \tau_{3,j+1} \circ \tau_{2,i+1}$.

\smallskip

\noindent Hence, assuming $(\tau_{1,2})_0 \in \Gamma_{G_{2,1}}$, and
abbreviating $\{ \tau_{i,j} : 0 < i < j\}$ by $\tau$, we have: 

\smallskip

 \ \ \ \ \  $|(\tau_{i,j})_0|_{_{\Gamma_{G_{2,1}} \cup \tau}} \ \leq 5$.
\end{lem}
{\bf Proof.} Recall that for $(\tau_{1,2})_0$ we have, by definition,  \  
$(\tau_{1,2})_0(1 \, w) = 1 \, w$, and \ 
$(\tau_{1,2})_0(0 \, x_2x_3 w) = 0 \, x_3x_2w$, for all $w \in \{0,1\}^*$
and $x_2, x_3 \in \{0,1\}$. The proof of the Lemma is a straightforward 
verification.  \ \ \ $\Box$

\bigskip

% \noindent
% In passing we define another, related, injective homomorphism on $G_{2,1}$,
% called the {\it lowering map}:
% 
% \smallskip
%
% $g \in G_{2,1} \longmapsto (g)' \in {\sf Stab}_{G}(0,1)$,
%
% \smallskip
%
% \noindent where $(g)'$ is defined by
%
%\medskip
%
% $(g)': \ x_1 w \longmapsto  x_1 \ g(w)$, 
% \ \ for any $x_1 \in \{0,1\}$ and $w \in \{0,1\}^*$.
%
% \medskip
%
% \noindent We have the following fact: \ \ $(g)' = \Lambda(g,g)$.
%
% \bigskip

\noindent Now we arrive at the relation between $w_{(g,g^{-1})}$ and $g$.

\begin{lem} \label{relation_w_g} \  
For all $g \in {\it lp}G_{2,1}$ the following relation holds
between $g$ and $w_{(g,g^{-1})}$ : 

\smallskip

 \ \ \ \ \ $w_{(g,g^{-1})} \cdot (g)_0^{-1}$ , \ \   
$(g)_0^{-1} \cdot w_{(g,g^{-1})} \ \ \in \ {\sf Fix}_{{\it lp}G_{2,1}}(0)$.

\medskip

\noindent Equivalently,     

\medskip

 \ \ \ \ \ $(g)_0 \cdot {\sf Fix}_{{\it lp}G_{2,1}}(0) \ = \ $
$w_{(g,g^{-1})} \cdot {\sf Fix}_{{\it lp}G_{2,1}}(0)$,
 \ \ \ and  

\smallskip

  \ \ \ \ \ ${\sf Fix}_{{\it lp}G_{2,1}}(0) \cdot (g)_0 \ = \ $
${\sf Fix}_{{\it lp}G_{2,1}}(0) \cdot w_{(g,g^{-1})}$ .
\end{lem}
{\bf Proof.} By Theorem \ref{repres_G_circuits_by_G} we have 
 \ $w_{(g,g^{-1})}( 0 \, x) \ = \ 0 \ g(x)$ \ for all $x \in {\rm Dom}(g)$.
So, $w_{(g,g^{-1})}$ and $(g)_0$ act in the same way on $0 \, \{0,1\}^*$. 
Also, both $w_{(g,g^{-1})}$ and $(g)_0$ map $0 \, \{0,1\}^*$ into 
$0 \, \{0,1\}^*$, and both map $1 \, \{0,1\}^*$ into $1 \, \{0,1\}^*$.
The Lemma follows from this.
 \ \ \ $\Box$

\bigskip

We abbreviate $\{ \tau_{i,j} : 0 < i < j\}$ by $\tau$.
The element $w_{(g,g^{-1})}$ of $G_{2,1}$, represented by the word
$W_{(g,g^{-1})}$, belongs to ${\sf Stab}_{{\it lp}G_{2,1}}(0,1)$ as we saw 
in Theorem \ref{repres_G_circuits_by_G}. However, the word $W_{(g,g^{-1})}$ 
itself is a sequence over the generating set 
$\Gamma_{G_{2,1}} \cup \tau$ of $G_{2,1}$.
Therefore, in order to follow the action of $W_{(g,g^{-1})}$ and of its 
prefixes we need to take ${\sf Fix}(0)$ as a subgroup of $G_{2,1}$.
This leads us to the {\bf Schreier left coset graph} of 
${\sf Fix}_{G_{2,1}}(0)$ within $G_{2,1}$, over the generating set 
$\Gamma_{G_{2,1}} \cup \tau$. 
By definition this Schreier graph has vertex set 
$G_{2,1}/{\sf Fix}_{G_{2,1}}(0)$, i.e., the left cosets, of the form 
$g \cdot {\sf Fix}_{G_{2,1}}(0)$ with $g \in G_{2,1}$.
And it has directed edges of the form 
 \ $g \cdot {\sf Fix}_{G_{2,1}}(0) \ \stackrel{\gamma}{\longrightarrow}$
$ \ \gamma g \cdot {\sf Fix}_{G_{2,1}}(0)$ \ for $g \in G_{2,1}$, 
 \ $\gamma \in \Gamma_{G_{2,1}} \cup \tau$. 
Lemma \ref{relation_w_g} implies that for all
$g \in {\it lp}G_{2,1}$, 

\smallskip

 \ \ \ \ \  $(g)_0 \cdot {\sf Fix}_{G_{2,1}}(0) \ = \ $
$w_{(g,g^{-1})} \cdot {\sf Fix}_{G_{2,1}}(0)$. 

\smallskip

\noindent 
We assume that $\Gamma_{G_{2,1}} = \Gamma_{G_{2,1}}^{-1}$, so the Schreier
graph is symmetric, and hence it has a distance function based on path 
length; we denote this distance by 

\smallskip

 \ \ \ \ \ $d_{G/F}(.,.): \ \ $
$G_{2,1}/{\sf Fix}_{G_{2,1}}(0) \times G_{2,1}/{\sf Fix}_{G_{2,1}}(0) $
$\ \longrightarrow \ {\mathbb N}$.

%\smallskip
%
%Since ${\sf Fix}_{G_{2,1}}(0)$ is an isomorphic copy of $G_{2,1}$ within
%itself, and 
%$\Lambda_1: g \in G_{2,1} \longmapsto (g)_0 \in {\sf Fix}_{G_{2,1}}(0)$ is an
%isomorphic self-embedding of $G_{2,1}$, we call the Schreier coset graph of
%$G_{2,1}/{\sf Fix}_{G_{2,1}}(0)$ the {\bf self-embedding Schreier graph}
%of $G_{2,1}$.

\begin{lem} \label{embed_lpG_in_Schreier} \    
There are injective morphisms

\medskip

$g \in {\it lp}G_{2,1} \ \hookrightarrow \ g \in G_{2,1} \ $ 
$\stackrel{\simeq}{\longrightarrow} \ $
$(g)_0 \in {\sf Fix}_{G_{2,1}}(1) \ \stackrel{\simeq}{\longrightarrow}$
 \ $(g)_0 \cdot {\sf Fix}_{G_{2,1}}(0) \ \in \ $
  ${\sf Stab}_{G_{2,1}}(0,1)/{\sf Fix}_{G_{2,1}}(0)$,

\medskip

\noindent and an inclusion map

\medskip

$(g)_0 \cdot {\sf Fix}_{G_{2,1}}(0) \ \in \ $
  ${\sf Stab}_{G_{2,1}}(0,1)/{\sf Fix}_{G_{2,1}}(0) \ \ $
$ \hookrightarrow \ \ $
$(g)_0 \cdot {\sf Fix}_{G_{2,1}}(0) \ \in \ G_{2,1}/{\sf Fix}_{G_{2,1}}(0)$.

\medskip

\noindent In particular, 

\medskip

 $g \in G_{2,1} \ \ \longmapsto \ \ $
$(g)_0 \cdot {\sf Fix}_{G_{2,1}}(0) \ \in \ $
  $G_{2,1}/{\sf Fix}_{G_{2,1}}(0)$  

\medskip

\noindent is an embedding of $G_{2,1}$, as a set, into the vertex set 
$G_{2,1}/{\sf Fix}_{G_{2,1}}(0)$ of the Schreier graph. 
\end{lem}
{\bf Proof.} 
Recall that the map \ $\Lambda_0:$ 
$g \in G_{2,1} \longmapsto (g)_0 \in {\sf Fix}_{G_{2,1}}(1)$ \ is a 
bijective morphism (Lemma \ref{isom_of_G}). Also, the map 
 \ $u \in {\sf Fix}_{G_{2,1}}(1) \ \longmapsto \ $
$u \cdot {\sf Fix}_{G_{2,1}}(0) \ \in \ G_{2,1}/{\sf Fix}_{G_{2,1}}(0)$
 \ is injective; indeed, if \
$u \cdot {\sf Fix}_{G_{2,1}}(0) = v \cdot {\sf Fix}_{G_{2,1}}(0)$ \
with \ $u, v \in {\sf Fix}_{G_{2,1}}(1)$ \ then \
$v^{-1}u \in {\sf Fix}_{G_{2,1}}(0) \cap {\sf Fix}_{G_{2,1}}(1)$
$ \ = \ \{ {\bf 1}\}$.

The map \     
$g \in G_{2,1} \ \longmapsto \ (g)_0 \cdot {\sf Fix}_{G_{2,1}}(0) $
  $ \ \in \ {\sf Stab}_{G_{2,1}}(0,1)/{\sf Fix}_{G_{2,1}}(0)$ \
is a surjective group homomorphism since ${\sf Fix}_{G_{2,1}}(0)$ is a 
normal subgroup of ${\sf Stab}_{G_{2,1}}(0,1)$.
Since \ ${\sf Fix}_{G_{2,1}}(0) \cap {\sf Fix}_G(1) = \{ {\bf 1}\}$, this
homomorphism is injective from ${\sf Fix}_{G_{2,1}}(1)$ onto
 \ ${\sf Stab}_{G_{2,1}}(0,1)/{\sf Fix}_{G_{2,1}}(0)$.

The combination of these maps provides an isomorphism from $G_{2,1}$ onto
${\sf Stab}_{G_{2,1}}(0,1)/{\sf Fix}_{G_{2,1}}(0)$. Hence we also have an
embedding of $G_{2,1}$, as a set, into the vertex set 
$G_{2,1}/{\sf Fix}_{G_{2,1}}(0)$ of the Schreier graph.
 \ \ \ $\Box$

\medskip

Since by Lemma \ref{embed_lpG_in_Schreier} we can consider $G_{2,1}$ as a 
subset of the vertex set $G_{2,1}/{\sf Fix}_{G_{2,1}}(0)$ of the Schreier 
graph, the path-distance $d_{G/F}(.,.)$ on $G_{2,1}/{\sf Fix}_{G_{2,1}}(0)$
leads to a distance on $G_{2,1}$, inherited from $d_{G/F}(.,.):$

\begin{defn} \label{Schreier_distance_inherited} \  
For all $g, g' \in G_{2,1}$ the {\em Schreier graph distance inherited by} 
$G_{2,1}$ is 

\smallskip

 \ \ \ \ \ $D(g, g') \ = \ d_{G/F}\big((g)_0 \cdot {\sf Fix}_{G_{2,1}}(0),$
  $ \ (g')_0 \cdot {\sf Fix}_{G_{2,1}}(0) \big)$.
\end{defn}

\noindent
The comparison of the Schreier graph distance $D(.,.)$ on ${\it lp}G_{2,1}$ 
with the word-length that ${\it lp}G_{2,1}$ inherits from its embedding into 
${\it lep}M_{2,1}$ leads to the following {\it distortion of} 
${\it lp}G_{2,1}$:

\begin{defn} \label{distor_Delta} \  
In ${\it lp}G_{2,1}$ we consider the {\em distortion}

\smallskip

  \ \ \ \ \ $\Delta(n) \ = \ {\rm max}\{ D({\bf 1},g) : \, $
 $|g|_{{\it lep}M_{2,1}} \leq n , \  \ g \in {\it lp}G_{2,1}\}$.
\end{defn}

\noindent
We now state and prove the main theorem relating $\Delta(.)$ and $\alpha$. 
Recall that $\alpha(.)$ is the computational asymmetry function of boolean
permutations, defined in terms of circuit size.

%% Main Theorem

\begin{thm} {\bf (Computational asymmetry vs.\ distortion).} 
 \label{polyn_rel_alpha_delta} \
The computational asymmetry function $\alpha(.)$ and the distortion 
$\Delta(.)$ of ${\it lp}G_{2,1}$ are polynomially related.
More precisely, for all $n \in {\mathbb N}:$

\smallskip
 
 \ \ \ \ \ $\big(\alpha(n)\big)^{1/2} \ \ \leq \ \ $
$ c' \cdot \Delta(n)$
$ \ \ \leq \ \ c \ n^4 +  c \cdot \big( \alpha(c \, n) \big)^4$

\smallskip

\noindent where $c \geq c' \geq 1$ are constants.
\end{thm}
{\bf Proof.} The Theorem follows immediately from Lemmas
\ref{distortionLEQalpha4} and \ref{dist_upperB_on_alpha}.
 \ \ \ $\Box$

\begin{lem} \hspace{-2mm}. \label{distortionLEQalpha4} \
There is a constant $c \geq 1$ such that for all $n \in {\mathbb N}:$
 \ \ \ $\Delta(n) \leq c \ n^4 +  c \cdot \big( \alpha(c \, n) \big)^4$.
\end{lem}
{\bf Proof.}
By Lemma \ref{relation_w_g}, 
 \ $(g)_0 \cdot {\sf Fix}_{G_{2,1}}(0) \ = \ $
 $ w_{(g,g^{-1})} \cdot {\sf Fix}_{G_{2,1}}(0)$, hence 

\smallskip

$d\big({\sf Fix}_{G_{2,1}}(0), \ (g)_0 \cdot {\sf Fix}_{G_{2,1}}(0)\big) $
 $ \ = \ $
$d\big({\sf Fix}_{G_{2,1}}(0), \ $
  $w_{(g,g^{-1})} \cdot {\sf Fix}_{G_{2,1}}(0)\big)$. 

\smallskip

\noindent 
Since the word $W_{(g,g^{-1})}$ and the Schreier graph use the same 
generating set, namely $\Gamma_{G_{2,1}} \cup \, \tau$, we have 

\smallskip

$d\big({\sf Fix}_{G_{2,1}}(0), \ $
  $w_{(g,g^{-1})} \cdot {\sf Fix}_{G_{2,1}}(0)\big)$
$ \ \ \leq \ \ |W_{(g,g^{-1})}|$. 

\smallskip

\noindent By Theorem \ref{repres_G_circuits_by_G}, \ 
$|W_{(g,g^{-1})}| \ \leq \ O(|C_g|^4 + |C_{g^{-1}}|^4)$. 
And by the definition of the computational asymmetry function,
 \ $|C_{g^{-1}}| \leq \alpha(|C_g|)$. Hence

\smallskip

$d\big({\sf Fix}_{G_{2,1}}(0), \ (g)_0 \cdot {\sf Fix}_{G_{2,1}}(0)\big)$
 $ \ \leq \ O(|C_g|^4 + |C_{g^{-1}}|^4)$
$ \ \leq \ O\big(|C_g|^4 + \alpha(|C_g|)^4\big)$.

\smallskip

\noindent By Proposition \ref{circ_lepM}, 
 \ $|C_g| = O(|g|_{{\it lep}M_{2,1}})$. Hence, for some constants 
$c'', c' \geq 1$, 

\smallskip

$d\big({\sf Fix}_{G_{2,1}}(0), \ (g)_0 \cdot {\sf Fix}_{G_{2,1}}(0)\big)$
$ \ \leq \ c' \cdot |g|_{{\it lep}M_{2,1}}^4 \ + \ $
$ c' \cdot \alpha(c'' \cdot |g|_{{\it lep}M_{2,1}})^4$.

\smallskip

\noindent Thus,

${\rm max}\big\{  d\big({\sf Fix}_{G_{2,1}}(0), \ $
               $(g)_0 \cdot {\sf Fix}_{G_{2,1}}(0)\big) : \ $
 $ |g|_{{\it lep}M_{2,1}} \leq n, \ g \in {\it lp}G_{2,1} \big\}$
 $ \ \ \leq \ \ c' \, n^4 + c' \, \alpha(c'' \, n)^4$.

\smallskip

\noindent By Definition \ref{distor_Delta} of the distortion function 
$\Delta$ we have therefore 

\smallskip

$\Delta(n) \ \leq \ c' \, n^4 + c' \ \alpha(c'' \, n)^4$.

\smallskip

\noindent This proves the Lemma. 
 \ \ \ $\Box$

%%%
\begin{lem} \label{dist_upperB_on_alpha} \
There is a constant $c \geq 1$ such that for all $n \in {\mathbb N}:$
 \ \ $\alpha(n) \ \leq \ c \cdot \Delta(c \, n)^2$.
\end{lem}
{\bf Proof.}  We first prove the following.

\smallskip

\noindent {\sf Claim:}  For every $g \in {\it lp}G_{2,1}$, 
the inverse permutation $g^{-1}$ can be computed by a circuit $C_{g^{-1}}$
of size \ \ $|C_{g^{-1}}| \leq c \cdot d \big({\sf Fix}_{G_{2,1}}(0), \, $
                         $(g)_0 \cdot {\sf Fix}_{G_{2,1}}(0) \big)^2 $, 
 \ for some constant $c \geq 1$.

\smallskip

\noindent Proof of the Claim: 
There is a word $W'$  of length \, $|W'| \, = \, $
$d\big({\sf Fix}_{G_{2,1}}(0), \, (g)_0 \cdot {\sf Fix}_{G_{2,1}}(0)\big)$ 
 \, over $\Gamma_{G_{2,1}} \cup \, \tau$ that labels a shortest
path from ${\sf Fix}_{G_{2,1}}(0)$ to $(g)_0 \cdot {\sf Fix}_{G_{2,1}}(0)$ 
in the Schreier graph of $G_{2,1}/{\sf Fix}_{G_{2,1}}(0)$. 
Let $W = (W')^{-1}$ (the free-group inverse of $W'$), so $|W| = |W'|$. Let 
$w$ be the element of $G_{2,1}$ represented by $W$. Then $W$ labels a 
shortest path from ${\sf Fix}_{G_{2,1}}(0)$ to 
$(g^{-1})_0 \cdot {\sf Fix}_{G_{2,1}}(0)$ in the Schreier graph of 
$G_{2,1}/{\sf Fix}_{G_{2,1}}(0)$; this path has length \ $|W| = |W'| = \ $
$d\big({\sf Fix}_{G_{2,1}}(0), \, (g)_0 \cdot {\sf Fix}_{G_{2,1}}(0)\big)$
$ \ = \ d\big({\sf Fix}_{G_{2,1}}(0), \ $
$(g^{-1})_0 \cdot {\sf Fix}_{G_{2,1}}(0)\big)$. 
 
We have \,  
$w \cdot {\sf Fix}_{G_{2,1}}(0) = (g^{-1})_0 \cdot {\sf Fix}_{G_{2,1}}(0)$, 
thus for all $x \in \{0,1\}^*:$ \ $w(0 \, x) = 0 \ g^{-1}(x)$. We now 
take the word $VWU$ over the generating set 
$\Gamma_{M_{2,1}} \cup \, \tau$ of the monoid $M_{2,1}$, where we choose
the words $U$ and $V$ to be
 \ $U = ({\sf and}, \, {\sf not}, \, {\sf fork}, \, {\sf fork})$, 
 \ and \ $V = ({\sf or})$. The functions {\sf and, not, fork, or} were defined
in Subsection 1.1. 
Then for all $x = x_1 \ldots x_n \in \{0,1\}^*$, with 
$x_1, \ldots, x_n \in \{0,1\}$, we have 

\smallskip

$x_1 \ldots x_n \ \ \stackrel{{\sf fork}}{\longrightarrow} \ \ $
$x_1 \, x_1 \ldots x_n \ \ $
$\stackrel{{\sf fork}}{\longrightarrow} \, $
$\stackrel{{\sf not}}{\longrightarrow} \ \ $
$\overline{x_1} \, x_1 \, x_1 \ldots x_n \ \ $
$\stackrel{{\sf and}}{\longrightarrow} \ \ 0 \, x_1 \ldots x_n $
$ \ \ = \ \ 0 \, x $

\medskip

$\stackrel{W}{\longrightarrow} \ \ 0 \ g^{-1}(x) \ \ $ 
$\stackrel{{\sf or}}{\longrightarrow} \ \ g^{-1}(x)$.

\smallskip

\noindent The last {\sf or} combines $0$ and the first bit of $g^{-1}(x)$, 
and this makes $0$ disappear.
Thus overall, \ $VWU(x) = g^{-1}(x)$. The length is \, $|VWU| = |W| + 5$. 

Since $g^{-1} \in {\it lp}G_{2,1} \subset {\it lep}M_{2,1}$, Theorem 
\ref{distor_lepM_M} implies that there exists a word $Z$ over the generators
 \, $\Gamma_{{\it lep}M_{2,1}} \cup \, \tau$ \, of ${\it lep}M_{2,1}$ such 
that 

\smallskip

(1) \ \ $|Z| \leq c_1 \cdot |VWU|^2$, \, for some constant $c_1 \geq 1$, and 

(2) \ \ $Z$ represents the same element of ${\it lep}M_{2,1}$ as $VWU$, 
       namely $g^{-1}$. 

\smallskip

\noindent Moreover, by Prop.\ \ref{circ_lepM}, the word $Z$ can be 
transformed into a circuit of size $\leq c_2 \cdot |Z|$ (for some constant 
$c_2 \geq 1$).
This proves that there is a circuit $C_{g^{-1}}$ for $g^{-1}$ of size \,    
$|C_{g^{-1}}| \leq c \cdot |W|^2$ \, (for some constant $c \geq 1$).
Since we saw that \ $|W| = $
$d_{G/F}({\sf Fix}_{G_{2,1}}(0), \, (g)_0 \cdot {\sf Fix}_{G_{2,1}}(0))$,
 \ the Claim follows.
 \ \ \ {\sf [End, Proof of the Claim.]}

\bigskip

\noindent By definition, \ $D({\bf 1}, g) = $   
$d_{G/F}({\sf Fix}_{G_{2,1}}(0), \, (g)_0 \cdot {\sf Fix}_{G_{2,1}}(0))$.
Hence, by the Claim above: 

\smallskip
 
 \ $|C_{g^{-1}}| \leq c \cdot \big( D({\bf 1}, g) \big)^2$.

\smallskip

\noindent
By Prop.\ \ref{circ_lepM} the word-length in ${\it lep}M_{2,1}$ and the
circuit size are linearly related; hence 
 \, $|g|_{{\it lep}M_{2,1}} \leq c_0 \, |C_g|$, for some constant 
$c_0 \geq 1$. Therefore, 

\medskip

$\alpha(n) \ = \ {\rm max}\{ |C_{g^{-1}}| : \ |C_g| \leq n, \ $
                 $ g \in {\it lp}G_{2,1} \}$

\smallskip

$\leq \ \ $
${\rm max}\{ |C_{g^{-1}}| : \ |g|_{{\it lep}M_{2,1}} \leq c_0 \, n, \ $
                 $ g \in {\it lp}G_{2,1} \} $ 

\smallskip

$\leq \ \ {\rm max}\big\{ c \cdot \big( D({\bf 1}, g) \big)^2 : \, $
   $|g|_{{\it lep}M_{2,1}} \leq c_0 \, n, \ $
                 $ g \in {\it lp}G_{2,1} \big\}$

\smallskip

$\leq \ \ c \cdot \big( \Delta(c_0 \, n) \big)^2$.

\smallskip

\noindent This proves the Lemma.  \ \ \ $\Box$

%%%%%%%%%%%%%%%%%%%

\section{Other bounds and distortions}

\subsection{Other distortions in the Thompson groups and monoids}

\noindent 
The next proposition gives more upper bounds on the computational 
asymmetry function $\alpha$.

\begin{pro}\hspace{-2mm}. \label{alphaLEQdistortion} \
Assume $\Gamma_{{\it lep}G_{2,1}} \subset \Gamma_{{\it lep}M_{2,1}} $
$\subset \Gamma_{M_{2,1}}$.  
Let $\delta_{{\it lp}G, {\it lep}M} = $  
$\delta \big[ \, |.|_{\Gamma_{{\it lp}G_{2,1}} \cup \tau}, \ $
  $|.|_{\Gamma_{{\it lep}M_{2,1}} \cup \tau} \, \big]$ 
be the distortion function of ${\it lp}G_{2,1}$ in the Thompson monoid 
${\it lep}M_{2,1}$, based on word-length. Similarly, let 
$\delta_{{\it lp}G, M} = $
$\delta \big[ \, |.|_{\Gamma_{{\it lp}G_{2,1}} \cup \tau}, \ $
 $|.|_{\Gamma_{M_{2,1} \cup \tau}} \, \big]$ be the distortion
function of ${\it lp}G_{2,1}$ in the Thompson monoid $M_{2,1}$.
Then for some constant $c \geq 1$ and for all $n \in {\mathbb N}$,

\smallskip

\hspace{1.5in}  
$\alpha(n) \ \leq \ c \cdot \delta_{{\it lp}G, {\it lep}M}(c \, n)$
$ \ \leq \ c \cdot \delta_{{\it lp}G, M}(c \, n)$.
\end{pro}
{\bf Proof.} We first prove that \   
$\delta_{{\it lp}G, {\it lep}M}(n) \ \leq \ \delta_{{\it lp}G, M}(n)$.
Recall that by definition, \
$\delta_{{\it lp}G, {\it lep}M}(n) \ = \ $
${\rm max}\{ |g|_{{\it lp}G_{2,1}} : g \in {\it lp}G_{2,1}, \ $
           $ |g|_{{\it lep}M_{2,1}} \leq n\}$, and similarly for
$\delta_{{\it lp}G, M}(n)$.
Since $\Gamma_{{\it lep}M_{2,1}} \subset \Gamma_{M_{2,1}}$ we have
 \ $|x|_{{\it lep}M_{2,1}} \leq |x|_{M_{2,1}}$. Hence, \
$\{ |g|_{{\it lp}G_{2,1}} : g \in {\it lp}G_{2,1}, \ $
  $ |g|_{{\it lep}M_{2,1}} \leq n\} \ \subseteq \ $
$\{ |g|_{{\it lp}G_{2,1}} : g \in {\it lp}G_{2,1}, \ $
  $ |g|_{M_{2,1}} \leq n\}$. By taking max over each of these two sets it
follows that \ $\delta_{{\it lp}G, {\it lep}M}(n) \ \leq \ $
$\delta_{{\it lp}G, M}(n)$.

\smallskip

Next we prove that \ $\alpha(n) \ \leq \ $
$c \cdot \delta_{{\it lp}G, {\it lep}M}(c \, n)$.
For any $g \in {\it lp}G_{2,1}$ we have
$C(g^{-1}) \leq O(|g^{-1}|_{{\it lep}M_{2,1}})$, by Prop.\
\ref{wl_asymM_vs_wl_asymLepM}. Moreover,
$|g^{-1}|_{{\it lep}M_{2,1}} \leq |g^{-1}|_{{\it lp}G_{2,1}}$ since
${\it lp}G_{2,1}$ is a subgroup of ${\it lep}M_{2,1}$, and since the
generating set used for ${\it lp}G_{2,1}$ (including all $\tau_{i,j}$)
is a subset of the generating set used for ${\it lep}M_{2,1}$.
For any group with generating set closed under inverse we have
$|g^{-1}|_G = |g|_G$.  And by the definition of the distortion
$\delta_{{\it lp}G, {\it lep}M}$ we have \ $|g|_{{\it lp}G_{2,1}} \leq $
$\delta_{{\it lp}G, {\it lep}M}(|g|_{{\it lep}M_{2,1}})$.
And again, by Prop.\ \ref{wl_asymM_vs_wl_asymLepM},
$|g|_{{\it lep}M_{2,1}} \leq O(C(g))$.
Putting all this together we have

\medskip

$C(g^{-1}) \ \leq \ c_1 \cdot |g^{-1}|_{{\it lep}M_{2,1}} \ \leq \ $
$c_1 \cdot |g^{-1}|_{{\it lp}G_{2,1}} \ = \ c_1 \cdot |g|_{{\it lp}G_{2,1}}$

\medskip

$ \ \leq \ $
$c_1 \cdot \delta_{{\it lp}G, {\it lep}M}(|g|_{{\it lep}M_{2,1}}) \ \leq \ $
$c_1 \cdot \delta_{{\it lp}G, {\it lep}M}(c_2 \ C(g))$.

\medskip

\noindent Thus, $c_1 \cdot \delta_{{\it lp}G, {\it lep}M}(c_2 \ C(g))$ is 
an upper bound on $C(g^{-1})$. Since, by definition, $\alpha(C(g))$ is the 
smallest upper bound on $C(g^{-1})$, it follows that 
$\alpha(C(g)) \leq c_1 \cdot \delta_{{\it lp}G, {\it lep}M}(c_2 \ C(g))$.
 \ \ \ $\Box$

\bigskip

Recall that in the definition \ref{distor_Delta} of the distortion $\Delta$ 
we compared $D(.,.)$ with the word-length in ${\it lep}M_{2,1}$.
If, instead, we compare $D(.,.)$ with the word-length in $M_{2,1}$ we obtain
the following distortion of ${\it lp}G_{2,1}$ : 

$$\delta(n) \ = \ {\rm max}\{ D({\bf 1},g) : \, 
 |g|_{M_{2,1}} \leq n , \  \ g \in {\it lp}G_{2,1}\} . $$

\begin{pro} \label{dist_in_M-lepM} \
The distortion functions $\Delta(.)$ and $\delta(.)$ are polynomially 
related. More precisely, there are constants $c', c_1, c_2 \geq 1$ such that
for all $n \in {\mathbb N}$: \ \
$\Delta(n) \ \leq \ c_1 \ \delta(n)$
$ \ \leq \ c_2 \ \Delta(c' \, n^2)$.
\end{pro}
{\bf Proof.} Let's assume first that 
$\Gamma_{{\it lep}M_{2,1}} \subseteq \Gamma_{M_{2,1}}$, from which it follows
that $|g|_{M_{2,1}} \leq |g|_{M_{2,1}}$. Therefore,
 \ $\{D({\bf 1},g) : |g|_{{\it lep}M_{2,1}} \leq n\} \ \subseteq \ $
$\{D({\bf 1},g) : |g|_{M_{2,1}} \leq n\}$.
Hence, $\Delta(n) \leq \delta(n)$.

By Theorem \ref{distor_lepM_M}, 
$|g|_{{\it lep}M_{2,1}} \leq c \cdot |g|_{M_{2,1}}^2$. So,
 \ $\{D({\bf 1},g) : |g|_{M_{2,1}} \leq n\} \ \subseteq \ $
$\{D({\bf 1},g) : |g|_{{\it lep}M_{2,1}} \leq c \, n^2\}$. Hence,
$\delta(n) \leq \Delta(c \, n^2)$. 

When we do not have 
$\Gamma_{{\it lep}M_{2,1}} \subseteq \Gamma_{M_{2,1}}$, the constants in the
theorem change, but the statement remains the same.
 \ \ \ $\Box$

%%%%%%%%%%%%%%

\subsection{ Monotone boolean functions and distortion }

On $\{0,1\}^*$ we can define the {\em product order}, also called 
``bit-wise order''. It is a partial order (and in fact, a lattice order), 
denoted by ``$\preceq$'', and defined as follows. First, $0 \prec 1$; 
next, for any $u,v \in \{0,1\}^*$ we have $u \preceq v$ \ iff 
 \ $|u| = |v|$ and $u_i \preceq v_i$ for all $i = 1, \ldots, |u|$,
 where $u_i$ (or $v_i$) denotes the $i$th bit of $u$ (respectively $v$).  

By definition, a partial function $f: \{0,1\}^* \to \{0,1\}^*$ is 
{\em monotone} (also called ``product-order preserving'') iff for all 
$u,v \in {\rm Dom}(f): \ u \preceq v$ implies $f(u) \preceq f(v)$.

The following fact is well known (see e.g., \cite{Wegener} Section 4.5):
A function $f: \{0,1\}^m \to \{0,1\}^n$ is monotone iff $f$ can be computed
by a combinational circuit that only uses gates of type {\sf and}, 
{\sf or}, {\sf fork}, and wire-swappings; i.e., {\sf not} is absent. 
A circuit of this restricted type is called a {\em monotone circuit}.

Razborov \cite{RazCLIQ} proved super-polynomial lower bounds for the size of
monotone circuits that solve the clique problem, and in \cite{RazMATCH} he
proved super-polynomial lower bounds for the size of
monotone circuits that solve the perfect matching problem for bipartite 
graphs; the latter problem is in {\sf P}. Tardos \cite{Tardos}, based on work
by Alon and Boppana \cite{AlonBopp}, gave an 
exponential lower bound for the size of monotone circuits that solve  a 
problem in {\sf P}; see also \cite{Handb} (Chapter 14 by Boppana and Sipser).
Thus, there exist problems that can be solved by
polynomial-size circuits but for which monotone circuits must have 
exponential size. In particular (for some constants $b>1, c>0$), there 
are infinitely many monotone functions $f_n: \{0,1\}^n \to \{0,1\}^n$ 
such that $f_n$ has a combinational circuit of size $\leq n^c$, but $f_n$ 
has no monotone circuit of size $\leq b^n$.

\smallskip

Based on an alphabet $A = \{a_1, \ldots, a_k\}$ with 
$a_1 \prec a_2 \prec \ \ldots \ \prec a_k$ we define a partial function 
$f: A^* \to A^*$ to be {\em monotone} iff $f$ preserves the product order of 
$A^*$. The monotone functions enable us to define the following submonoid of
the Thompson-Higman monoid ${\it lep}M_{k,1}:$

\medskip

${\it mon}M_{k,1} \ = \ \{ \varphi \in {\it lep}M_{k,1} : \varphi$
 can be represented by a monotone function $P \to Q$, 

\hspace{2in} where $P$ and $Q$ are prefix codes, with $P$ maximal $\}$.

\medskip

\noindent An essential extension or restriction of an element of 
${\it mon}M_{k,1}$ is again in ${\it mon}M_{k,1}$, so this set is 
well-defined as a subset of ${\it lep}M_{k,1}$. It is easily seen to be 
closed under composition, so ${\it mon}M_{k,1}$ is a submonoid of 
${\it lep}M_{k,1}$.

We saw that all monotone finite functions have circuits made from gates of 
type {\sf and}, {\sf or}, {\sf fork}.  Hence ${\it mon}M_{2,1}$ has the 
following generating set: 

\smallskip

 \ \ \ \ \    
$\{ {\sf and}, {\sf or}, {\sf fork} \} \cup \{\tau_{i,j} : j > i \geq 1\}$. 

\smallskip

\noindent 
The results about monotone circuit size imply the following distortion 
result. Again, ``exponential'' refers to a function with a lower
bound of the form \ $n \in {\mathbb N} \longmapsto \exp(\sqrt[c]{c' \, n})$,
for some constants $c' > 0$ and $c \geq 1$.

\begin{pro} \label{distor_monotone} \     
Consider the monoid ${\it mon}M_{2,1}$ over the generating set 
$\{ {\sf and}, {\sf or}, {\sf fork} \} \cup \{\tau_{i,j} : j>i \geq 1\}$,
and the monoid ${\it lep}M_{2,1}$ over the generating set 
$\Gamma_{{\it lep}M_{2,1}} \cup \{\tau_{i,j} : j>i \geq 1\}$, where
$\Gamma_{{\it lep}M_{2,1}}$ is finite.
Then ${\it mon}M_{2,1}$ has exponential word-length distortion in 
${\it lep}M_{2,1}$. 
\end{pro}
{\bf Proof.} Let $\Gamma_{\it mon} = \{ {\sf and}, {\sf or}, {\sf fork} \}$.
By Prop.\ \ref{circ_lepM} we have \   
$|f|_{\Gamma_{{\it lep}M_{2,1}} \cup \tau} = |C_f|$,
where $|C_f|$ denotes the ordinary circuit size of $f$. 
By a similar argument we obtain: \  
$|f|_{\Gamma_{\it mon} \cup \tau} = |{\it mon}C_f|$, where $|{\it mon}C_f|$ 
denotes the monotone circuit size of $f$.  We saw that as a consequence of the
work of Razborov, Alon, Boppana, and Tardos, there exists an infinite set of 
monotone functions that have polynomial-size circuits but whose monotone 
circuit-size is exponential. The exponential distortion follows.
 \ \ \ $\Box$

\smallskip

Since ${\it lep}M_{2,1}$ has quadratic distortion in $M_{2,1}$, 
${\it mon}M_{2,1}$ also has exponential word-length distortion in 
$M_{2,1}$.

%%%%%%%%%%%%%%%%%%%%%%%%%%%%%%%%%%%%%%%%%%%%%%%%%%%%%%%%%

\bigskip

\bigskip

%%%%%%%%%%%%%%%%%%%%%%%%%%%%%%%%%%%%%%%%%%%%%%%%%%%%%%%%%%%%%%%%%%%%%%%%%%%

%%%%%%%%%%%%%%%%%%%%%%%%%%%%%

\bigskip

\bigskip

\noindent {\bf Jean-Camille Birget} \\
Dept.\ of Computer Science \\
Rutgers University at Camden \\
Camden, NJ 08102, USA \\
{\tt birget@camden.rutgers.edu}


\begin{thebibliography}{99}

{\small

\bibitem{AlonBopp} N.\ Alon, R.\ Boppana, ``The monotone circuit complexity
  of boolean functions'', {\it Combinatorica} 7 (1987) 1-23.

\bibitem{ArzhGubaSapir} G.\ Arzhantseva, V.\ Guba, M.\ Sapir,
  ``Metrics on diagram groups and uniform embeddings in a Hilbert space'',
   Mathematics Arxiv (2004)  http://arxiv.org/abs/math.GR/0411605.

\bibitem{ArzhDrutuSapir} G.\ Arzhantseva, C.\ Dru\c{t}u, M.\ Sapir,
  ``Compression functions of uniform embeddings of groups into Hilbert and
  Banach spaces'', Mathematics Arxiv (2006)
  http://arxiv.org/abs/math.GR/0612378

\bibitem{Bennett73} C.H.~Bennett, ``Logical reversibility of computation'',
  {\it IBM J.\ of Research and Development} 17 (Nov.\ 1973) 525-532.

\bibitem{Bennett89} C.H.~Bennett, ``Time/Space tradeoffs for reversible
computation'', {\it SIAM J.\ of Computing} 18 (1989) 766-776.

\bibitem{BiThomps} J.C.\ Birget, ``The groups of Richard Thompson and
  complexity'', {\it International J. of Algebra and Computation} 14(5,6)
  (Dec.\ 2004) 569-626 (Mathematics ArXiv: math.GR/0204292, Apr.\ 2002).

\bibitem{BiCoNP} J.C.~Birget, ``Circuits, coNP-completeness, and the groups
  of Richard Thompson'', {\it International J.~of Algebra and Computation}
  16(1) (Feb.\ 2006) 35-90
 (Mathematics ArXiv: http://arXiv.org/abs/math.GR/0310335, Oct.\ 2003).

\bibitem{BiFact} J.C.\ Birget, ``Factorizations of the Thompson-Higman
  groups, and circuit complexity'',
  Mathematics ArXiv: math.GR/0607349, July 2006.

\bibitem{BiThomMon} J.C.~Birget, ``Monoid generalizations of the Richard
  Thompson groups'', Mathematics ArXiv: math.GR/0704.0189, 
  2 Apr.\ 2007.

\bibitem{Bourgain} J.\ Bourgain, ``On Lipschitz embedding of finite metric
  spaces in Hilbert space'', {\it Israel J.\ of Mathematics} 52 (1985) 46-52.

\bibitem{BoppLag} R.B.~Boppana, J.C.~Lagarias, ``One-way functions and 
  circuit complexity'', {\it Information and Computation} 74 (1987) 
  226-240.

\bibitem{CFP} J.\ W.\ Cannon, W.\ J.\ Floyd, W.\ R.\ Parry,
``Introductory notes on Richard Thompson's groups'',
{\it L'Enseignement Math\'ematique} 42 (1996) 215-256.

\bibitem{DiffieHellman} W.\ Diffie, M.E.\ Hellman, ``New directions in
  cryptography'', {\it IEEE Transactions in Information Theory} 22
  (1976) 644-655.

\bibitem{Farb} B.\ Farb, ``The extrinsic geometry of subgroups and the
  generalized word problem'', {\it Proc.\ London Mathematical Society}
  (3) 68 (1994) 577-593.

\bibitem{FredToff} E.~Fredkin, T.~Toffoli, ``Conservative logic'',
 {\it International J.\ Theoretical Physics} 21 (1982) 219-253.

\bibitem{GrollSel} J.\ Grollman, A.\ Selman, ``Complexity measure for
  public-key cryptosystems'', {\it SIAM J.\ on Computing} 17 (1988) 309-335.

\bibitem{Gromov} M.\ Gromov, ``Asymptotic invariants of infinite groups'',
 in {\it Geometric Group Theory} (G.\ Niblo, M.\ Roller, editors),
 London Mathematical Society Lecture Notes Series 182, Cambridge Univ.\
 Press (1993).

\bibitem{Hig74} G.\ Higman, ``Finitely presented infinite simple groups'',
  Notes on Pure Mathematics 8, The Australian National University,
  Canberra (1974).

\bibitem{HiltgenDiss} A.P.L.\ Hiltgen, ``Cryptographically relevant
  contributions to combinatorial complexity'', Dissertation, ETH-Z\"urich.
  Hartung-Gorre Verlag, Konstanz (1994). 

\bibitem{HiltgenAUSC} A.P.L.\ Hiltgen, ``Construction of feebly-one-way
  families of permutations'', in {\it Advances in Cryptology -- 
  AUSCRYPT'92}, Lecture Notes in Computer Science 718 (1993) 422-434.

\bibitem{IndykMatousek} P.\ Indyk, J.\ Matousek, ``Low-distortion embeddings 
  of finite metric spaces'', Chapter 8 in {\it Handbook of Discrete and 
  Computational Geometry} (J.E.\ Goodman and J.\ O'Rourke, editors), 
  CRC Press LLC, Boca Raton, FL; Second Edition (2004).

\bibitem{Lec} Y.\ Lecerf, ``Machines de Turing r\'{e}versibles ...'',
 {\em Comptes Rendus de l'Acad\'{e}mie des Sciences, Paris} 257 No.\ 18 
 (Oct.\ 1963) 2597 - 2600.

\bibitem{Lupanov} O.B.\ Lupanov, ``A method of circuit synthesis'',
  {\it Izv.\ V.U.Z.\ Radiofiz.} 1 (1958) 120-140.

\bibitem{MargMeakSunik} S.\ Margolis, J.\ Meakin, Z.\ \v{S}uni\'k,
  ``Distortion functions and the membership problem for submonoids of 
  groups and monoids'', {\it Contemporary Mathematics}, AMS, 
  372 (2005) 109-129.

\bibitem{Massey} J.L.\  Massey, ``The Difficulty with Difficulty'',
  IACR Distinguished Lecture delivered at EUROCRYPT '96, July 17, 1996,
  Saragossa, Spain. \ \ ( http://www.iacr.org/publications/dl/ )

\bibitem{McKTh} R.\ McKenzie, R.J.\ Thompson,
``An elementary construction of unsolvable word problems in group theory'',
in {\it Word Problems}, (W.\ W.\ Boone, F.\ B.\ Cannonito, R.\ C.\ Lyndon,
editors), North-Holland (1973) pp.\ 457-478.

\bibitem{Nekr_book} V.\ Nekrashevych, {\it Self-similar groups}, 
  Mathematical Surveys and Monographs vol.\ 117 (2005), American 
  Mathematical Society.

\bibitem{Olsh} A.Y.\ Ol'shanskii, ``On subgroup distortion in finitely
presented groups'', {\it Matematicheskii Sbornik} 188 (1997) 51-98.

\bibitem{OlSap} A.Y.\ Ol'shanskii, M.V.\ Sapir, ``Length and area functions
on groups and quasi-metric Higman embedding'', {\it International J.\ of
Algebra and Computation} 11 (2001) 137-170.

\bibitem{RazCLIQ} A.A.\ Razborov, ``Lower bounds for the monotone complexity 
  of some boolean functions'', {\it Doklady Akademii Nauk SSSR} 281(4)
  (1985) 798-801. (English transl.: {\it Soviet Mathematical Doklady} 31 
  (1985) 354-357.)

\bibitem{RazMATCH} A.A.\ Razborov, ``Lower bounds of monotone complexity 
  of the logical permanent function'', {\it Matematicheskie Zametki} 37(6)
  (1985) 887-900. (English transl.: {\it Mathematical Notes of the Academy
  of Sciences of the USSR} 37 (1985) 485-493.) 

\bibitem{JESavage} J.E.\ Savage, {\it Models of Computation}, Addison-Wesley
(1998).

\bibitem{Scott} E.A.\ Scott, ``A construction which can be used
to produce finitely presented infinite simple groups'',
{\it J. of Algebra} 90 (1984) 294-322.

\bibitem{Selman} A.\ Selman, ``A survey of one-way functions in complexity
theory'', {\it Mathematical Systems Theory} 25 (1992) 203-221.

\bibitem{Shannon} C.E.\ Shannon, ``The synthesis of two-terminal switching
  circuits'', {\it Bell System Technical J.} 28 (1949) 59-98. 

\bibitem{Shende} V.~Shende, A.~Prasad, I.~Markov, J.~Hayes,
  ``Synthesis of reversible logic circuits'', {\it IEEE Transactions on
  Computer-Aided Design of Integrated Circuits and Systems} 22(6) (2003)
  710-722.

\bibitem{Tardos} E.\ Tardos, ``The gap between monotone and non-monotone
  circuit complexity is exponential'', {\it Combinatorica} 7(4) (1987)
  141-142.

\bibitem{Th0} Richard J.\ Thompson, Manuscript (1960s).

\bibitem{Th} R.J.\ Thompson, ``Embeddings into finitely generated
  simple groups which preserve the word problem'',
  in {\it Word Problems II}, (S.\ Adian, W.\ Boone, G.\ Higman, editors),
  North-Holland (1980) pp.\ 401-441.

\bibitem{Toff80Memo} T.~Toffoli, ``Reversible computing'', MIT Laboratory 
  for Computer Science, Technical Memo MIT/LCS/TM-151 (1980). 

\bibitem{Toff80Conf} T.~Toffoli, ``Reversible computing'', {\it Automata, 
  Languages and Programming} (7th Colloquium), Lecture Notes in Computer 
  Science 85 (July 1980) 623-644.  (Abridged version of \cite{Toff80Memo}.)

\bibitem{Handb} J.\ van Leeuwen (editor), {\it Handbook of Theoretical
  Computer Science}, volume {\bf A}, MIT Press and Elsevier (1990).

\bibitem{Wegener} I.\ Wegener, {\it The complexity of boolean functions},
  Wiley/Teubner (1987).

} %End \small

\end{thebibliography}
\end{document}